\documentclass[
]{elsarticle}

\journal{Applied Mathematics and Computation}

\usepackage{hycolor}
\usepackage{hyperref}

\usepackage[utf8]{inputenc}
\usepackage{xcolor}
\newcommand{\td}[2]{\cfrac{\textup{d}#1}{\textup{d}#2}\, }

\newcommand{\e}{\operatorname{e}}
\renewcommand{\u}{\widehat{u}}

\newcommand\figwidth{85mm}

\makeatletter
\@namedef{u8:\detokenize{ｄ}}{\ensuremath{{\,\mathrm{d}}}}
\@namedef{u8:\detokenize{Δ}}{\ensuremath{\Delta}}
\@namedef{u8:\detokenize{ω}}{\ensuremath{\omega}}
\@namedef{u8:\detokenize{λ}}{\lambda}
\@namedef{u8:\detokenize{τ}}{\tau}
\@namedef{u8:\detokenize{·}}{\cdot}
\@namedef{u8:\detokenize{α}}{\alpha}
\@namedef{u8:\detokenize{ε}}{\epsilon}
\@namedef{u8:\detokenize{÷}}{\frac}
\@namedef{u8:\detokenize{⦅}}{\left(}
\@namedef{u8:\detokenize{⦆}}{\right)}
\makeatother
\usepackage{lineno}
\modulolinenumbers[5]
\usepackage{amsmath, amsfonts} 
\usepackage{multicol}

\usepackage{subcaption}
\usepackage{boldline}
\newcommand\dg\mathfrak

\usepackage{graphicx}

\usepackage{algorithm2e}

\begin{document}

\begin{frontmatter}

	\title{Performance of Borel-Padé-Laplace integrator for the solution of stiff and non-stiff problems}

\author{Ahmad Deeb}
\author{Aziz Hamdouni}


\author{Dina Razafindralandy\texorpdfstring{\corref{mycorrespondingauthor}}{}}
\cortext[mycorrespondingauthor]{Corresponding author}
\ead{drazafin@univ-lr.fr}

\address{Laboratoire des Sciences de l'Ingénieur pour l'Environnement \\ UMR 7356 La Rochelle Université -- CNRS\\ Avenue Michel Crépeau, 17042 La Rochelle Cedex 1, France}

\begin{abstract}
	A stability analysis of the Borel-Padé-Laplace series summation technique, used as explicit time integrator, is carried out. Its numerical performance on stiff and non-stiff problems is analyzed. Applications to ordinary and partial differential equations are presented. The results are compared with those of many popular schemes designed for stiff and non-stiff equations.
\end{abstract}

\begin{keyword}
	Borel-Laplace summation, divergent series, time integrator, stiff equations
\end{keyword}

\end{frontmatter}


\section{Introduction}

Stiff problems occur in many areas of engineering science, such as mechanics, electrical and chemical engineering (see for instance \cite{book:hairer2,isereles96,butcher03,lambert1991,willoughby74}). However, their resolution has remained a challenge for numerical analysts. The reason is that many numerical methods designed for general ordinary differential equations exhibit a high instability when solving stiff problems, unless an excessively small time step is used. As a consequence, numerical schemes with better stability properties have been developed especially for stiff problems.

One method used to estimate the largest time step allowed by a given numerical scheme without breaking its stability is the analysis of the linear stability domain. The scheme is called $A$-stable if this domain contains the half complex plane with negative real part, meaning that the scheme is stable in some sense however large is the time step, for the numerical solution of a 1D linear equation. The notion of a linear stability domain will be recalled later. See also \cite{book:hairer2,lambert1991,hackbusch14} for different notions of stability. Of course, even if a scheme is $A$-stable, the time step is limited in pratice due to precision requirements.

Among the most widely used numerical schemes for stiff equations, one can mention implicit linear multistep methods based on backward difference formulas (BDF) \cite{curtiss52,isereles96}. Their stability is limited to low orders. Indeed, only the first order (implicit Euler) and second order schemes are $A$-stable. BDF schemes of order 3 to 6 exhibit a weaker stability property called $A(\alpha)$-stability, and the formulas of order greater than 6 are unstable. A generalization of BDF which uses a second derivative permits to obtain implicit $A(\alpha)$-stable schemes up to order 10. See \cite{book:hairer2} for instance.

Another important family of numerical schemes for differential equations are Runge-Kutta methods (RK) \cite{runge95,butcher63,butcher96}, which are one-step schemes. Explicit RK schemes are not $A$-stable and not suitable for stiff equations. However, compared to multistep methods, it is easier to find stable implicit Runge-Kutta schemes. For example, Gauss, Radau IA and IIA, and Lobatto IIIA, IIIB and IIIC are $A$-stable \cite{book:hairer2,butcher16}.

Of course, there exist some other schemes suitable for stiff problems. 
A common feature of all these methods is their implicit character. However, the cost of an implicit scheme may be very high. This is particularly true for long-time dynamics problems (celestial mechanics, molecular dynamics, \ldots) where the use of implicit methods is hardly conceivable. The development of explicit, yet with a good enough stability property, numerical schemes is desirable.

An approach which has been used to this aim is to build stabilized RK schemes \cite{book:hairer2,verwer96}. These schemes are not $A$-stable like the implicit RK schemes but have a larger stability domain than standard explicit RK schemes.

Other semi-explicit schemes which are built for stiff problems are exponential time differencing (ETD) integrators. They are based on an exact, exponential type, resolution of the linear part of the equation. In doing so, the stiff part of the solution is correctly captured if it is an exponentially decaying term. The complete solution, the expression of which can be found by the variation of constants method, is then computed numerically. Various schemes have been proposed for this task \cite{Friedli_1978,Norsett_1969,vanderhoven_1974,certaine1960,hochbruck2005323,hochbruck2010exponential}. One of the most popular exponential integrators is the exponential time differencing associated to an explicit 4-th order Runge-Kutta method (ETDRK4) developed by Cox and Matthews \cite{cox_2002}. The algorithm is not completely explicit since it requires the (pseudo-)inversion of a matrix. Moreover, they generally need the approximation of the action of a matrix exponential, which is numerically expensive.

In the present article, we examine the performance of a Borel-Laplace integrator (BL) in solving stiff and non-stiff systems. BL is an entirely explicit, arbitrary high-order scheme. It is based on a decomposition of the solution into its time Taylor series, followed by a Borel-Laplace summation procedure to accelerate the convergence, or in the case of a divergent series, to obtain an asymptotical solution. The first goal of the article is the study of the stability of BL.  We will see that, although not $A$-stable (typical for explicit methods), BL admits a stability region which grows very fast with the order of the scheme. This enables large time steps compared to many popular explicit and even implicit schemes in practice. The second goal of the article is to show that BL is suited to the resolution of stiff equations and to high-dimensional problems.

At its origin, the Borel-Laplace summation method was intended to define the asymptotic sum of a Gevrey series
 \cite{borel}. It has recently gained more interest when authors showed that the heat equation and many equations in mechanics (Burgers and Navier-Stokes equations, \ldots), quantum physics or astronomy have divergent but Gevrey Taylor series \cite{lutz1999,lysik2009,Costin_2006,dyson_1952,suslov05,kontopolous_2002}. The Borel-Laplace summation method has been transformed into numerical algorithm \cite{thomann00} and used for the first time as a time integrator by Razafindralandy and Hamdouni \cite{jcp13}. Since then, a number of features of the Borel-Laplace integrator was studied. For example, it generally allows much larger time steps than other explicit methods for the resolution of many problems \cite{jcp13}. Its ability to cross some types of singularities, its high-order symplecticity, or its high-order iso-spectrality in solving a Lax pair problem have been stated in \cite{esaim14}. Another advantage of BL is that decreasing or increasing the approximation order is as simple as changing the value of a parameter in the code. However, nowhere in the cited works on the Borel-Laplace integrator stiffness has been addressed. One aim of the present article is, as mentioned, to fill this gap.

 The Borel-Laplace algorithm that will be discussed here results from the representation of the Borel sum as a Laplace integral. It makes use of a Padé approximation, as will be seen later, and is accordingly named Borel-Padé-Laplace algorithm (BPL). A representation of the Borel sum as an inverse factorial series also leads to an efficient algorithm \cite{dcds16} but will not be used.

This paper is organized as follows. In section 2, the Borel-Padé-Laplace algorithm is briefly recalled. In section 3, a linear stability analysis is carried out. The stability regions, corresponding to different values of parameters, are plotted. In section 4, numerical performance on stiff and non-stiff ODE problems as well as on a PDE is analyzed.

\section{Borel-Padé-Laplace integrator}

Consider an ordinary differential equation or a semi-discretized partial differential equation :
\begin{equation}
	\begin{cases}
		\td ut=F(t,u),\\[5pt] u(0)=u_0,
	\end{cases}
	\label{equation}
\end{equation}
where
\[ u:\begin{array}{ccc} [0,T]&\longrightarrow& \mathbb R^n\\t&\longmapsto &u(t)\end{array}\]
is the unknown, $n\in \mathbb N^*$ is the dimension of the system and $F$ is a non-linear operator
\[ F:\begin{array}{ccc} \mathbb R\times\mathbb R^n&\longrightarrow& \mathbb R^n\\(t,v)&\longmapsto &F(t,v).\end{array}\]
Borel-Laplace integrator is based on approximating the solution to (\ref{equation}) via a (convergent or divergent) time series 
\begin{equation}
	\u(t)=\sum_{k=0}^{\infty}u_kt^k\in(\mathbb C[[t]])^n
	\label{u_series}
\end{equation}
and performing a Borel-Laplace summation procedure on this series. In  equation (\ref{u_series}), $\mathbb C[[t]]$ stands for the ring of formal power series in $t$, with complex coefficients. To simplify, assume that $n=1$. The terms $u_k$ are obtained by inserting directly the series expansion (\ref{u_series}) into equation (\ref{equation}). This leads to explicit relations of the form
\begin{equation}u_{k+1}=\frac1{k+1}F_k(u_0,\cdots,u_k)\label{ukp1}\end{equation}
where $F_k$ is the $k-$th Taylor coefficient of $F(t,u(t))$ at $t=0$. It is generally a non-linear function of $u_0,u_1,\cdots,u_k$. The expression of $F_k$ will be explicitely given for each equation we will be dealing with.

Regarding the validity domain of series (\ref{u_series}), there are two possible scenarios. The first one is that the (exact or numerical) radius of convergence of series is zero. In this case, a summation procedure is required. A summation procedure consists in finding an analytic representation (as an integral, a rational function, a continued fraction, \dots) of the exact solution from the series. The one chosen here is the Borel-Laplace summation in which the solution is represented by a Laplace transform of a rational function (see section \ref{theoretic}).

The other possibility is that the series is convergent with a finite or theoretically infinite radius of convergence. However, if the series converges slowly, the partial sum 
\begin{equation}
	u^K(t)=\sum_{k=0}^Ku_kt^k
	\label{partial_sum}
\end{equation}
may give an acceptable approximation only for values of $t$ much smaller than the radius of convergence. In this case, an acceleration of convergence is optional but advisable.

In the present paper, the Borel-Laplace summation will systematically be applied to the series. If the series is convergent, it will act as a convergence acceleration technique. As will be seen, it not only extends the stability region but also increases the speed of the method.


\subsection{Theoretical setting of Borel-Padé-Laplace summation\label{theoretic}}

The theory behind Borel-Laplace summation can be found in many papers \cite{borel,ramis_poincare_1,ramis_poincare_2,costin_2008_book} and shall not be reproduced here. Only the computational aspects are presented. Let us assume that series (\ref{u_series}) is a $p$-Gevrey series in a neighborhood of the origin, that is, 
\begin{equation}
	|u_k|\leq CA^k(k!)^p,\quad\quad \forall k\geq 0
	\label{gevrey}
\end{equation}
for some positive real numbers $A$ and $C$. In fact, it is known that most of series arising in engineering problems are $p$-Gevrey series for some positive rational number $p$. In the sequel, we consider only the case $p=1$.

The numerical summation is done in three stages. First, the Borel transform
\begin{equation}
	\mathcal B\u(\xi)=\sum_{k=0}^{\infty}\,B_k\xi^k\in\mathbb C[[\xi]]
	\label{b_series}
\end{equation}
of series (\ref{u_series}) is considered. In this expression,
\begin{equation}
	B_k=\cfrac{u_{k+1}}{k!},\quad\quad k\geq0.
	\label{bk}
\end{equation}
Series (\ref{b_series}) is convergent at the origin. Next, $\mathcal B\u(\xi)$ is prolonged analytically into a function $P(\xi)$ in the vicinity of a semi-line $\ell$ of the complex plane, linking 0 to $\infty$. Lastly, the Laplace transform (at $1/t$), which is the formal inverse of the Borel transform, is applied to the prolonged function. At the end of the procedure, one gets an integral representation 
\[
\mathcal{S}\u(t)=u_0+\displaystyle
\int_{\ell}P(\xi) e^{-\xi/t}\operatorname d\xi
\]
of the solution. The function $\mathcal{S}\u(t)$ is called the Borel sum of series (\ref{u_series}). The Borel-Laplace summation procedure is summarized in Table \ref{borel_laplace}. 

\begin{table}[ht]\centering  {\sf\small
\begin{tabular}{ccccc}
	$\displaystyle\u(t)=\sum_{k=0}^{+\infty}u_kt^k$& $\sim$   &\hspace{-1cm}$\mathcal S\u(t)=u_0+\displaystyle
\int_{\ell}P(\xi) e^{-\xi/t}\operatorname d\xi$ 
\\
 $\left.\begin{array}{c}\\\text{Borel}\\\\\end{array}\right\downarrow$&
 &$\left\uparrow \begin{array}{c}\\\text{Laplace} \\\\\end{array}\right.$
\\
$\displaystyle\mathcal{B}\u(\xi)=
\sum_{k=0}^{+\infty}B_k \xi^k$
&$\overrightarrow{\hspace{.5cm}\text{Prolongation}\hspace{.5cm}}$
&\hspace{-1cm}$P(\xi)$\\\\
\end{tabular}}
\caption{Borel-Laplace summation}\label{borel_laplace}
\end{table}

If the initial series (\ref{u_series}) is convergent at the origin, then the Borel sum $\mathcal S\u(t)$ takes the same value as the original sum (\ref{u_series})  for $t$ inside the disc of convergence. However, the domain of definition of $\mathcal S\u$ is generally larger than the disc of convergence of the series $\u$. If the initial series is a divergent but Gevrey series, $\mathcal S\u$ is a sectorially analytical function, having the series $\u$ as Gevrey asymptotics.

\subsection{Algorithm\label{sec:algorithm}}

Numerically, only a finite number of terms $u_k$ can be computed. The series is then represented by the degree $K$ polynomial $\u(t)\simeq u^K(t)$
\begin{equation}
	\displaystyle\u(t)\simeq u^K(t)=\sum_{k=0}^Ku_kt^k,
	\label{partial_sum2}
\end{equation}
that is the partial sum already been defined in (\ref{partial_sum}).
The Borel transform of $u^K(t)$ is a degree $(K-1)$ polynomial. The prolongation is carried out via a Padé approximation \cite{brezinski79,brezinski94}. Other prolongation techniques exist but none of them has been used as part of a Borel-Laplace based time integrator, to the best of our knowledge. A usual Gauss-Laguerre quadrature permits to compute the Laplace transform \cite{stroud66}. In simulations, the semi-line $\ell$ is the positive real axis. Due to the realization of the prolongation via a Padé aprroximation, the algorithm is called Borel-Padé-Laplace (BPL).

The function $\mathcal Su^K(t)$ obtained with BPL provides an approximate solution to the equation as long as an accuracy criterion is met. When this is no longer the case, the algorithm (computation of $u_k$'s and Borel summation) is restarted using the last acceptable value $\mathcal Su^K(t_f)$ as initial condition. BPL is then a step-by-step method over time. One way of evaluating the accuracy of the approximate solution is to calculate the residue of the equation. This strategy is rather expensive but, as will be seen, is fast enough to compete with all the other numerical schemes under consideration.


The Borel-Padé-Laplace algorithm can be summarized as follows, for a one-dimensional problem, with the residue as quality criteria:
\begin{enumerate}
	\item Start with $t_0=0$ and $u_0=u(t_0)$.
	\item Compute the first $K$ coefficients of the series $\u$:
		$$u_{k+1}=\cfrac1{k+1}\ F_k(u_0,\dots,u_k),\quad k=0,\cdots,K-1$$
		where $F_k$ is the $k-$th Taylor coefficient of $F(t,u(t))$ at $t=t_0$.
	\item Apply the Borel transformation. In other words, compute the first $K$ coefficients of $\mathcal{B}\u$:
		$$B_k=\cfrac{u_{k+1}}{k!},\quad k=0,\cdots,K-1.$$
	\item Compute a $[K_a/K_b]$ Padé approximant $P(\xi)$ of the polynomial with coefficients $B_k$, \textit{i.e.} determine $a_k$ and $b_k$ such that 
		\begin{equation}
			\begin{array}{ll}
				B_0+B_1\xi+\cdots+B_{K-1}\xi^{K-1}+O(\xi^K)=
				\\[5pt]
				\hspace{4cm}\cfrac{a_0+a_1\xi+\cdots+a_{K_a}\xi^{K_a}}{1+b_1\xi+\cdots+b_{K_b}\xi^{K_b}}=:P(\xi).
			\end{array}
			\label{pade}
		\end{equation}
	\item Obtain the approximate Borel sum by computing the Laplace transform with a $N_G$-point Gauss-Laguerre quadrature formula\footnote{Note that $\displaystyle\int_0^{+\infty}P(\xi)\e^{-\xi/t}ｄ\xi=t\int_0^{+\infty}P(t\xi)\e^{-\xi}ｄ\xi$}:
		\begin{equation}
			\mathcal{S}u^K(t)= u_0+ t\sum_{i=1}^{N_G} P(t\xi_i)w_i.
			\label{gausslaguerre}
		\end{equation}
	\item Find $t_f$ such that the relative residue norm is smaller than a tolerance parameter $\epsilon$ for all $t\leq t_f$:
		\begin{equation}
			\left\|\cfrac{\textrm d\mathcal{S}u^K}{\textrm dt}-F\left(t,\mathcal{S}u^K\right)\right\|\ <\ \epsilon\ \left\|\mathcal{S}u^K\right\|
		\label{eps}
		\end{equation}
		
		Take $\mathcal{S}u^K(t)$ as the approximation of $u(t)$ for $t\in(t_0,t_f]$.

	\item Return to step 2 with $t_0=t_f$, $u_0=\mathcal{S}u^K(t_f)$.
\end{enumerate}
In this algorithm, $K_a$ and $K_b$ are any positive integers such that $K_a+K_b=K-1$. Their influence on the stability region will be analyzed in section \ref{sec:stability}.
The reals $\xi_i$ are the roots of the $N_G$-th Gauss-Laguerre polynomials and the $w_i$ are the corresponding weights. 
Note also that a singular value decomposition will be carried out to improve the robustness of the Padé approximation  in numerical tests, following an algorithm discussed in \cite{gonnet11}.

In stage 6, the final time $t_f$ can be determined as follows. First, we evaluate the numerical of convergence $\tau$ of the series (\ref{u_series}) from its $K$ first terms using a simple criteria developped in \cite{cochelin1994}. It writes
\begin{equation}
	\tau=\left(\delta \cfrac{\|u_1\|}{\|u_K\|}\right)^{1/(K-1)}
	\label{convergence_radius}
\end{equation}
where $\delta$ is a small tolerance parameter such that
\[
\cfrac{\|u^K(t)-u^{K-1}(t)\|}{\|u^K(t)-u_0\|}\simeq\cfrac{\|u_Kt^{K}\|}{\|u_1t\|}\leq\delta.
\]
Next, we check if (\ref{eps}) is verified by $t=τ$. If it is, we try further with $2τ$, $4τ$, \dots until condition (\ref{eps}) is not verified any longer. If, on the contrary, the first try $t=τ$ is not successful, we try further with $τ/2$, $τ/4$ and so on. The last successful value will be taken as $t_f$.
The quantity $t_f-t_0$ is considered as the time step of the algorithm. 

The cut-off order $K$ can be thought as the order of the scheme. Note that one advantage of Borel-Laplace integrator is that, in contrast to many schemes such as BDF or Runge-Kutta, increasing the order is very simple. Increasing $K$ is enough; the algorithm does not need to be modified, no coefficient has to be changed.

\vspace{\baselineskip}

In the next section, the linear stability of BPL is analysed. We study in particular the influence of the summation procedure.

\section{Stability analysis\label{sec:stability}}

First, we recall the notion of stability domain for a discrete scheme. Consider the scalar linear model problem
\begin{equation}
	\begin{cases}
		\td ut=\lambda u\\[5pt]
		u(t_0)=u_0
	\end{cases}
	\label{linear}
\end{equation}
where $λ$ is a complex number with a negative real part. The solution of this equation decreases exponentially to zero when $t$ grows. Consider an iterative scheme  with a constant time step $h$, providing discrete approximate solutions $v_h^n\simeq u(t_n)$ of (\ref{linear}) at discrete times $t_n=t_0+nh$ as follows:
\begin{equation}
	v_h^{n+1}=R(λ,h)v_h^n
	\label{iteration}
\end{equation}
for some function $R$ of $λ$ and $h$. The stability domain of this method is defined as the following subset of the complex plane \cite{book:hairer2,isereles96}:
\begin{equation}
	D=\{\ (λh)\in \mathbb C\ :\ |R(λ,h)|<1\ \}.
\end{equation}
When the time step $h$ is such that $λh$ lies in the stability region, the approximate solution decreases to zero, like the exact one, when $n$ grows.

BPL is not a discrete scheme, in the sense that it does not provide a discrete approximation of the solution but a continuous one. It is however relatively easy to adapt to it the notion of stability domain. For simplicity, assume that $t_0=0$. For equation (\ref{linear}), we have:
\begin{equation}
	F(t,u)=λu\quad\quad\text{and}\quad\quad F_k(u_0,\dots,u_k)=λu_k. 
\end{equation}
Equation (\ref{ukp1}) permits to compute the coefficients of the time series:
\begin{equation}
	u_{k+1}=\cfrac{\lambda u_k}{k+1}  .
\end{equation}
When inserted into series (\ref{u_series}), these coefficients lead of course to the Taylor expansion of the exact solution $\e^{\lambda t}$.
It is a convergent series.

\begin{figure}
	\centering
	\begin{subfigure}{39mm}
		\includegraphics[width=\textwidth]{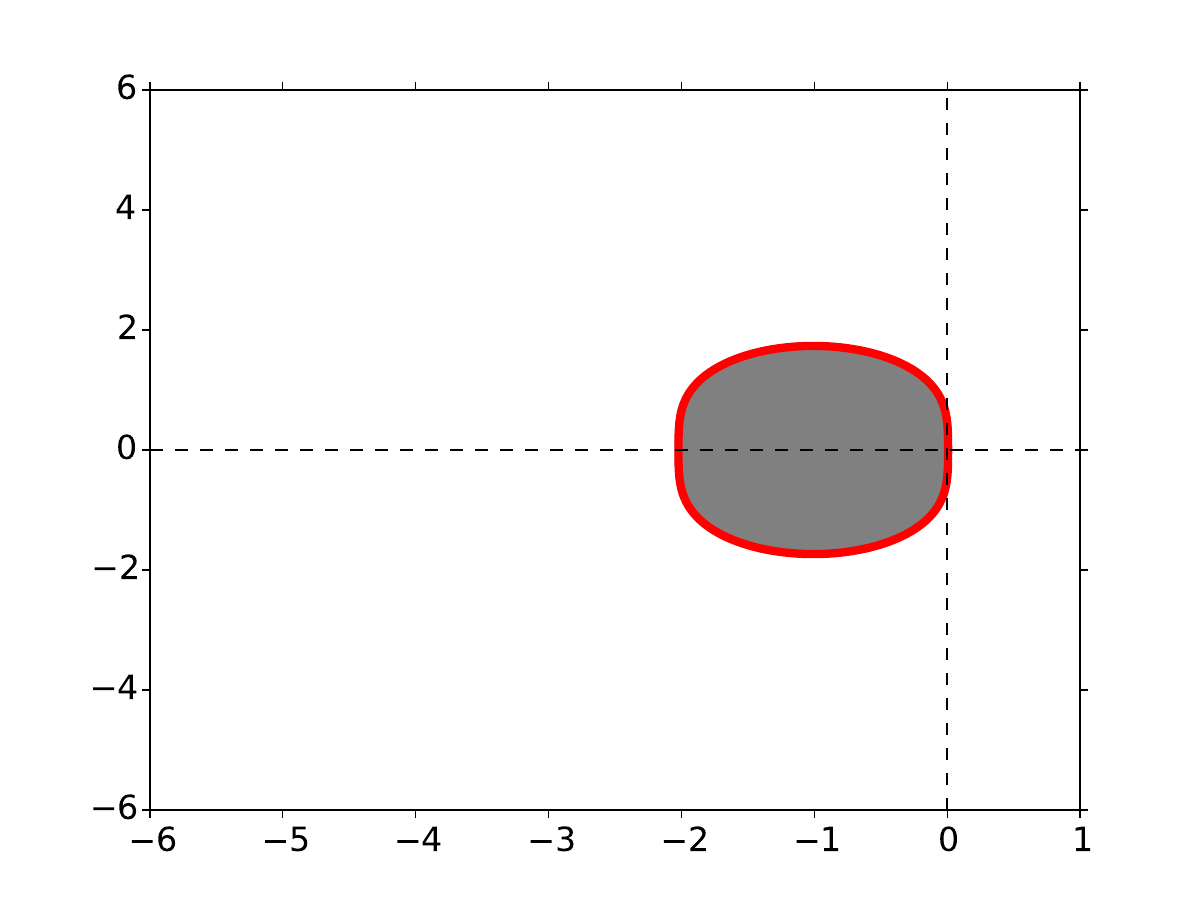}
		\caption{$K=2$}
		\label{fig:region_series_K2}
	\end{subfigure}
	\begin{subfigure}{39mm}
		\includegraphics[width=\textwidth]{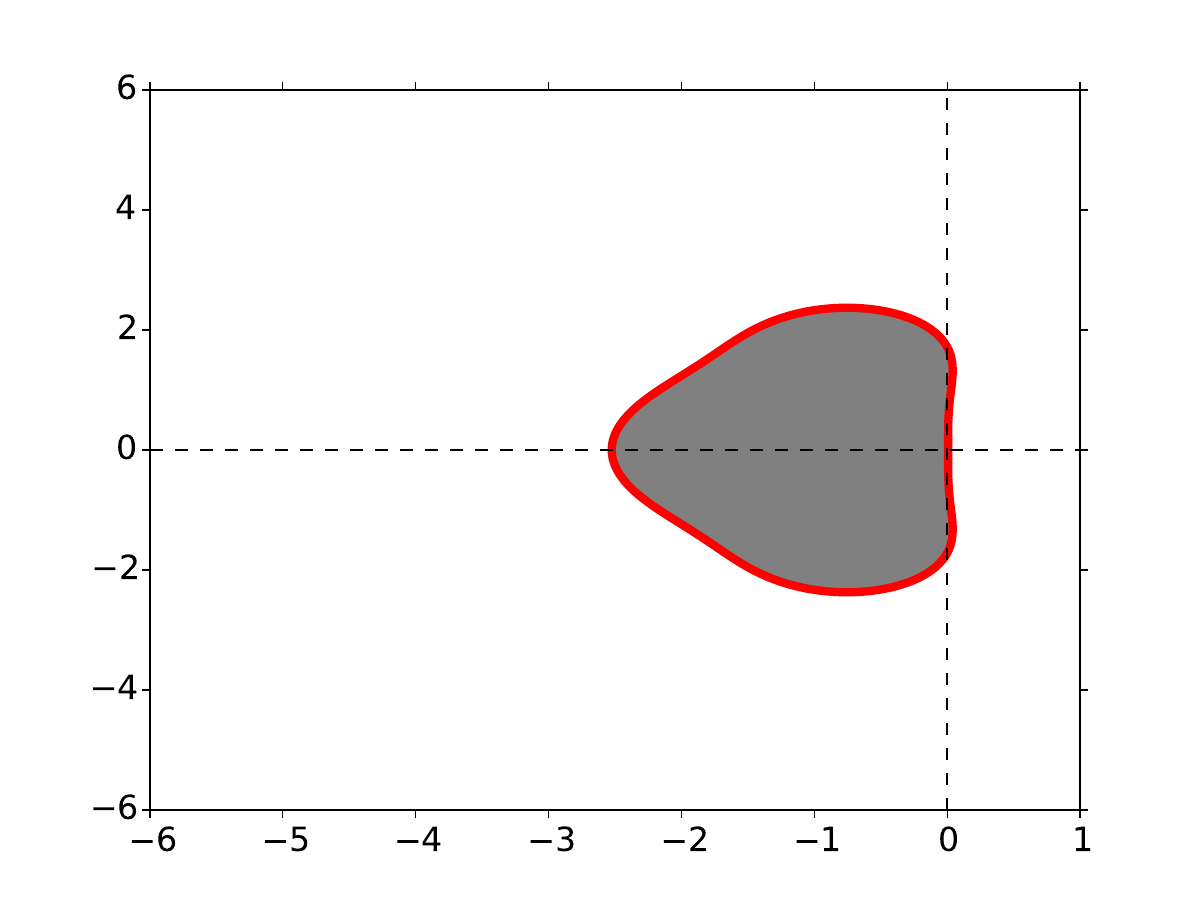}
		\caption{$K=3$}
	\end{subfigure}
	\begin{subfigure}{39mm}
		\includegraphics[width=\textwidth]{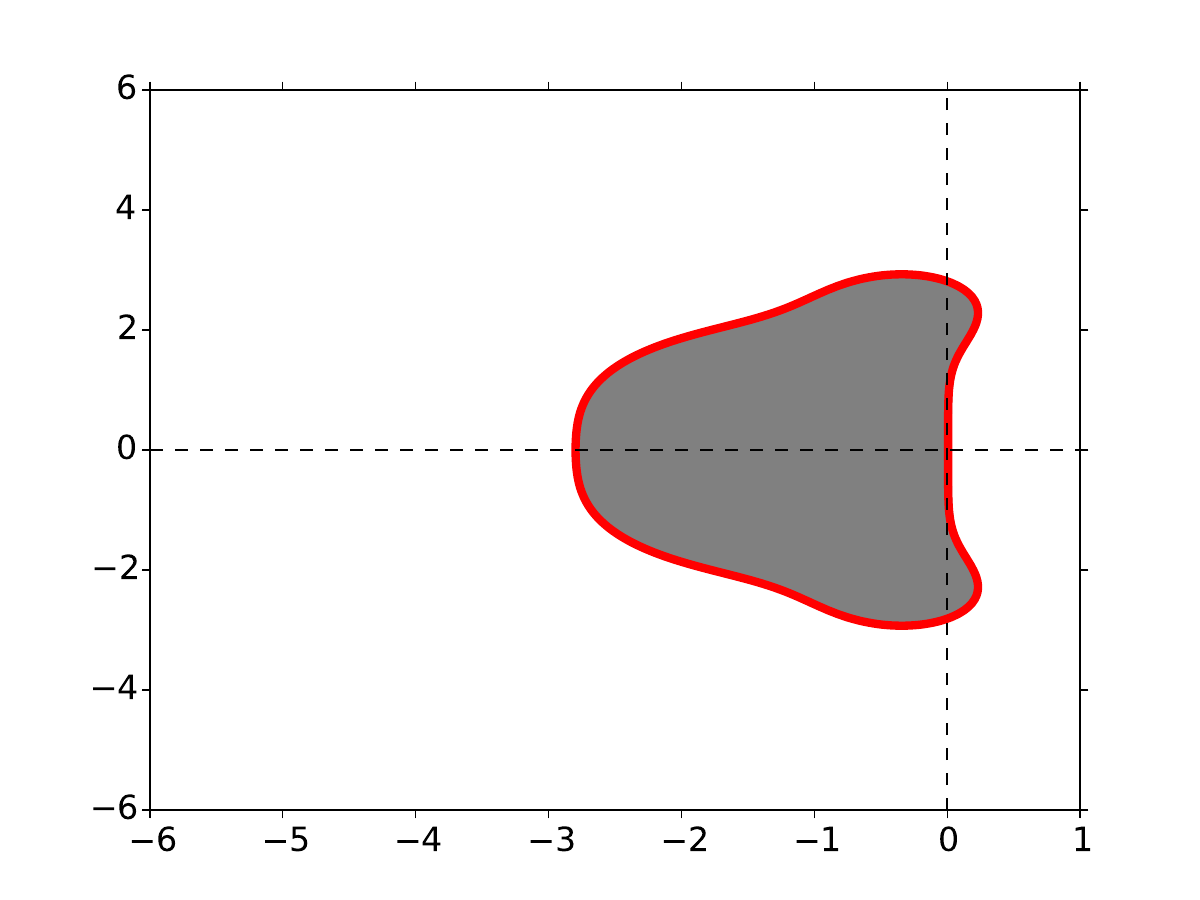}
		\caption{$K=4$}
	\end{subfigure}%

	\vspace{\baselineskip}
	\begin{subfigure}{39mm}
		\includegraphics[width=\textwidth]{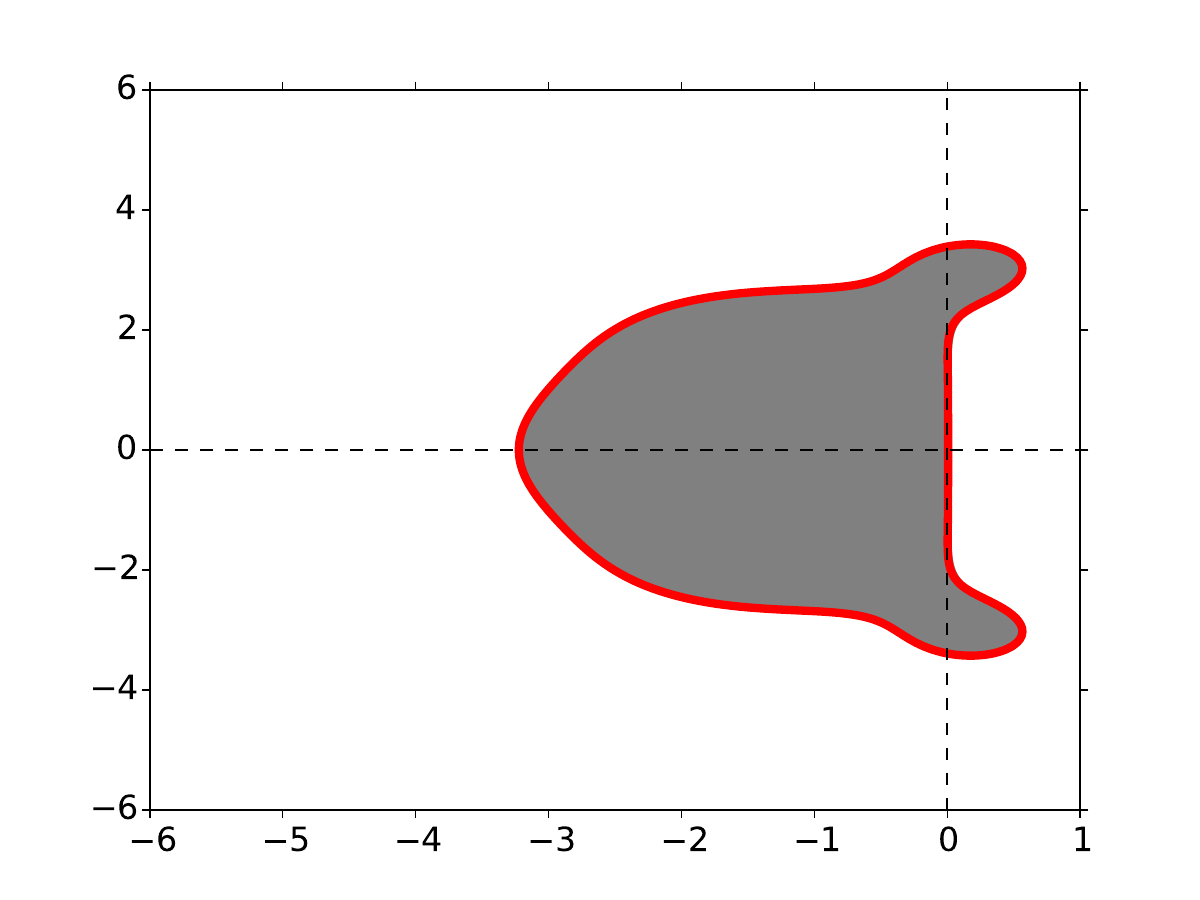}
		\caption{$K=5$}
	\end{subfigure}
	\begin{subfigure}{39mm}
		\includegraphics[width=\textwidth]{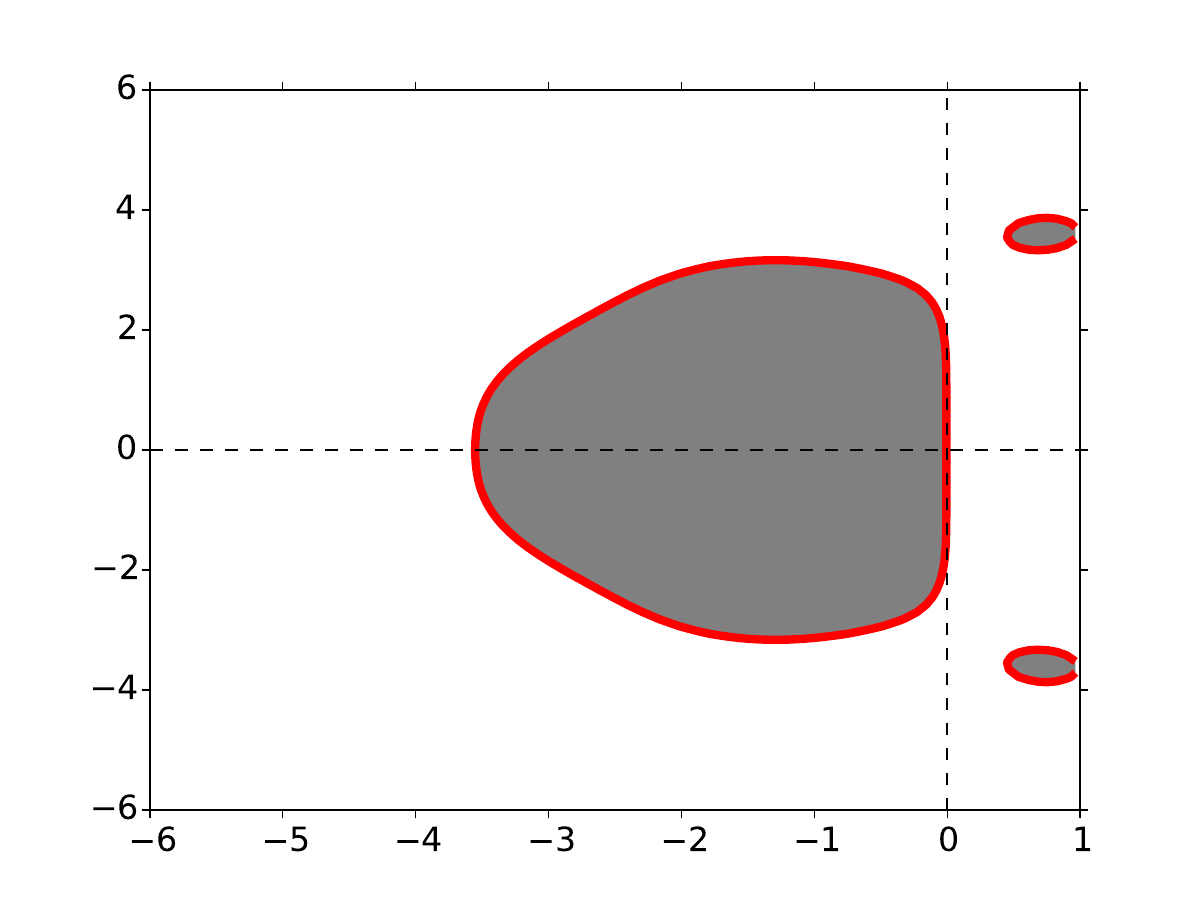}
		\caption{$K=6$}
	\end{subfigure}
	\begin{subfigure}{39mm}
		\includegraphics[width=\textwidth]{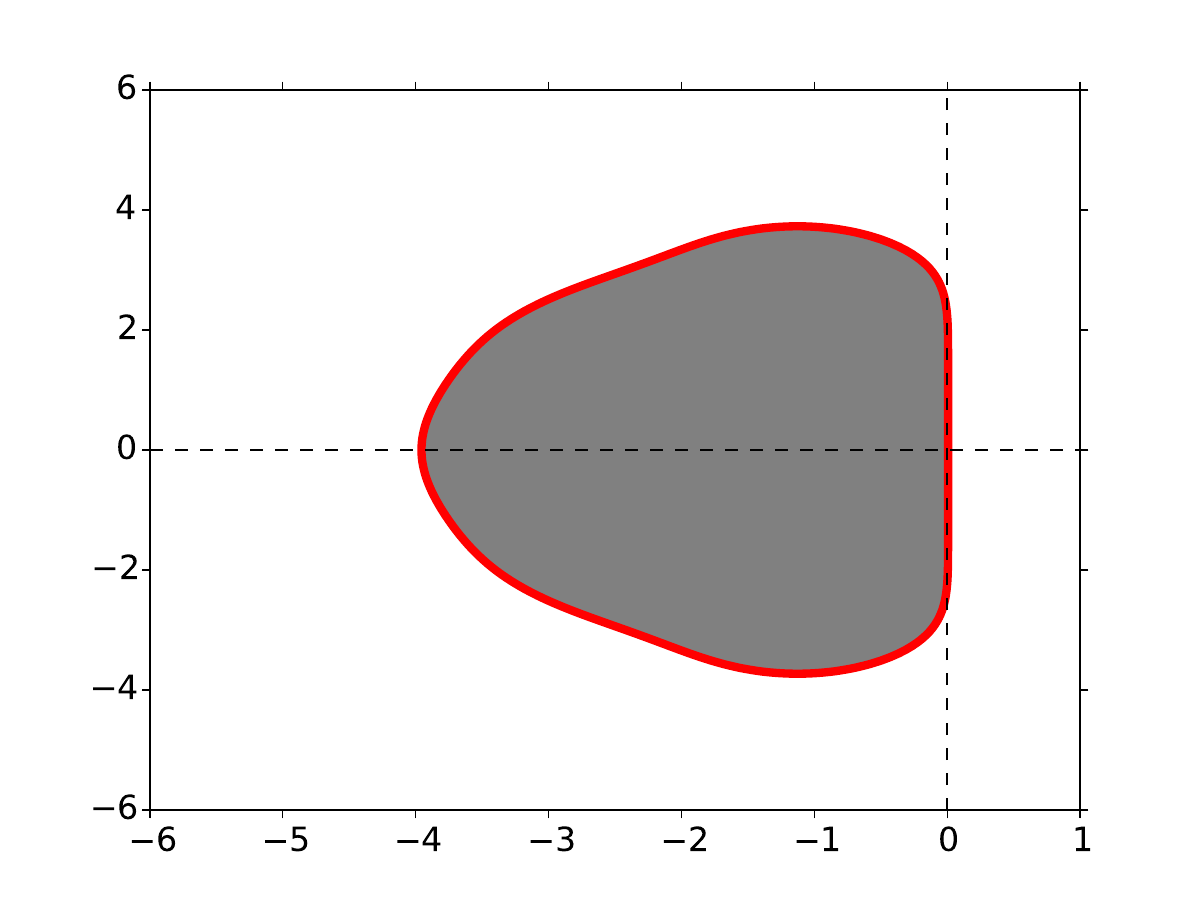}
		\caption{$K=7$}
	\end{subfigure}

	\vspace{\baselineskip}
	\begin{subfigure}{39mm}
		\includegraphics[width=\textwidth]{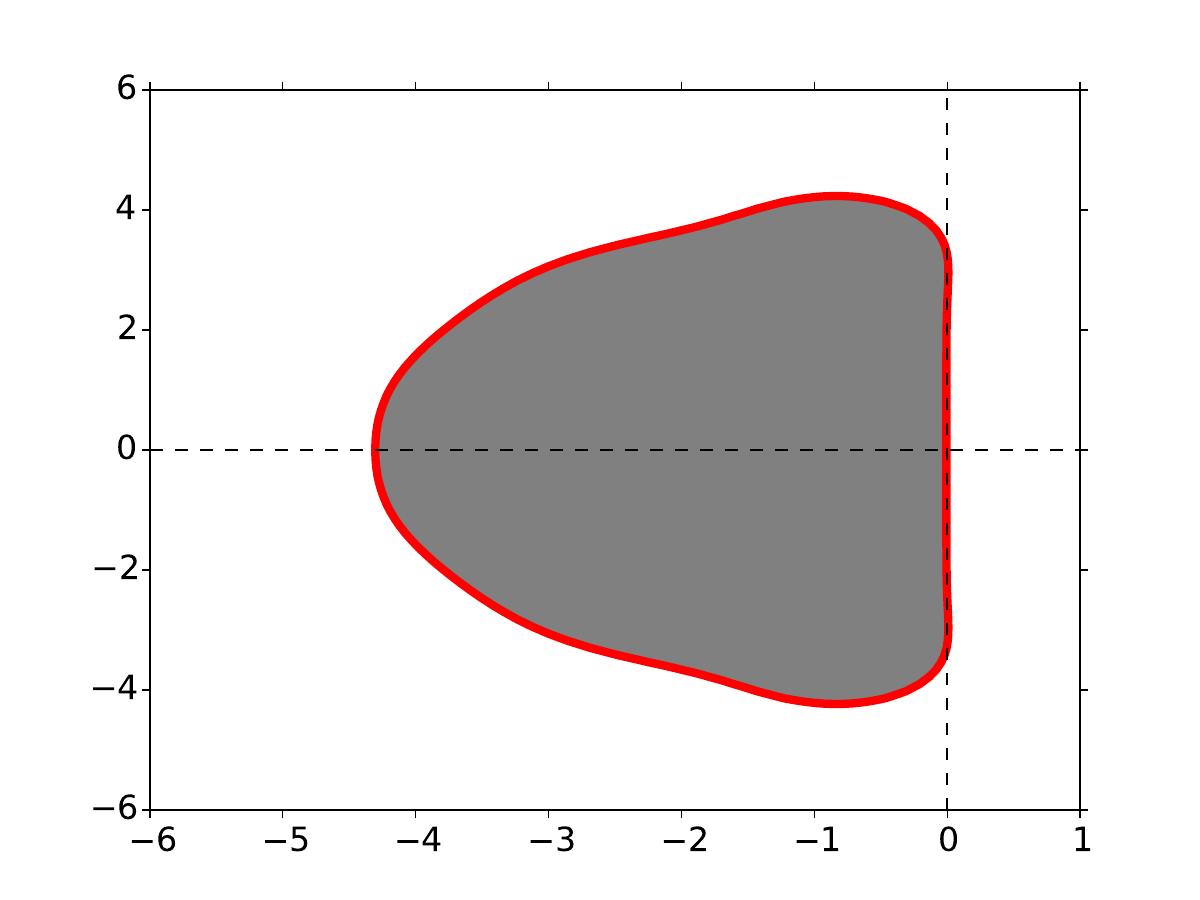}
		\caption{$K=8$}
	\end{subfigure}
	\begin{subfigure}{39mm}
		\includegraphics[width=\textwidth]{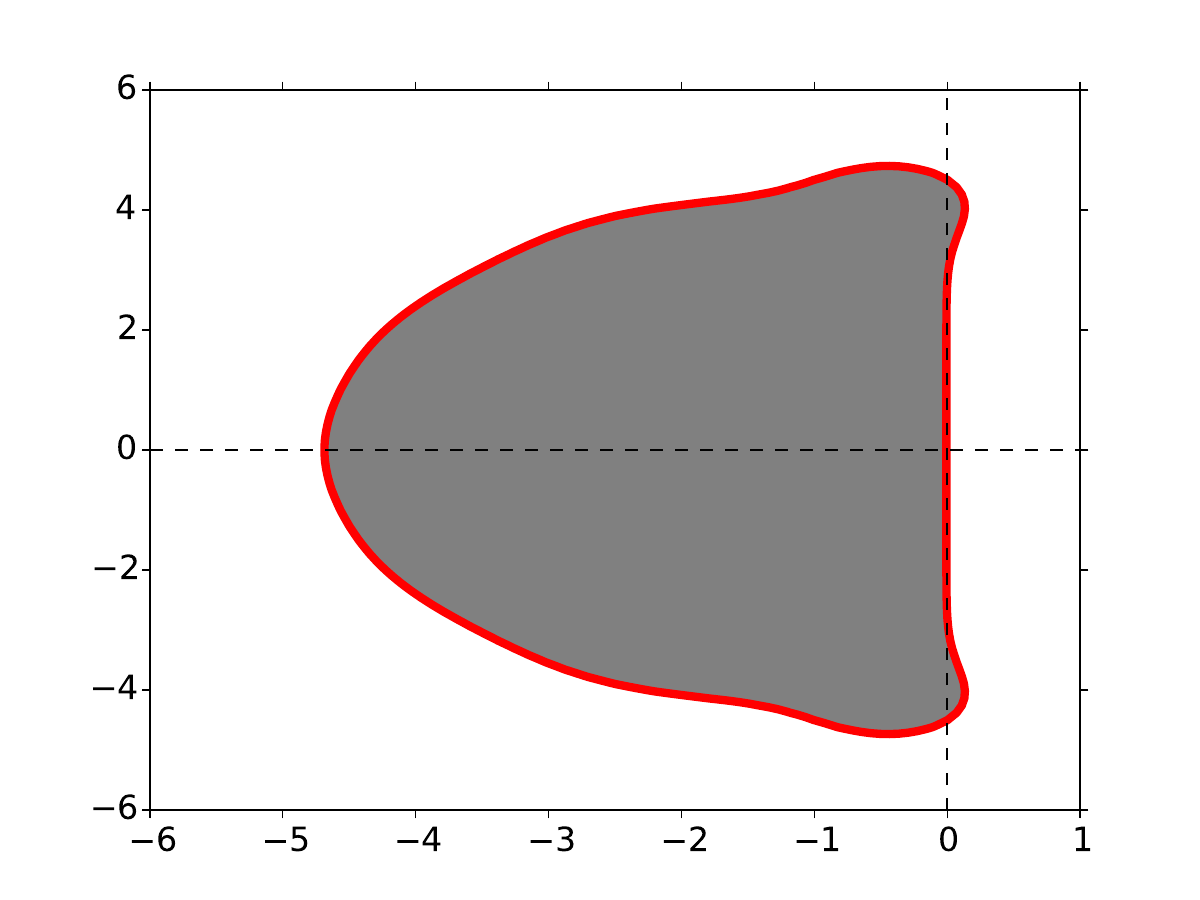}
		\caption{$K=9$}
	\end{subfigure}
	\begin{subfigure}{39mm}
		\includegraphics[width=\textwidth]{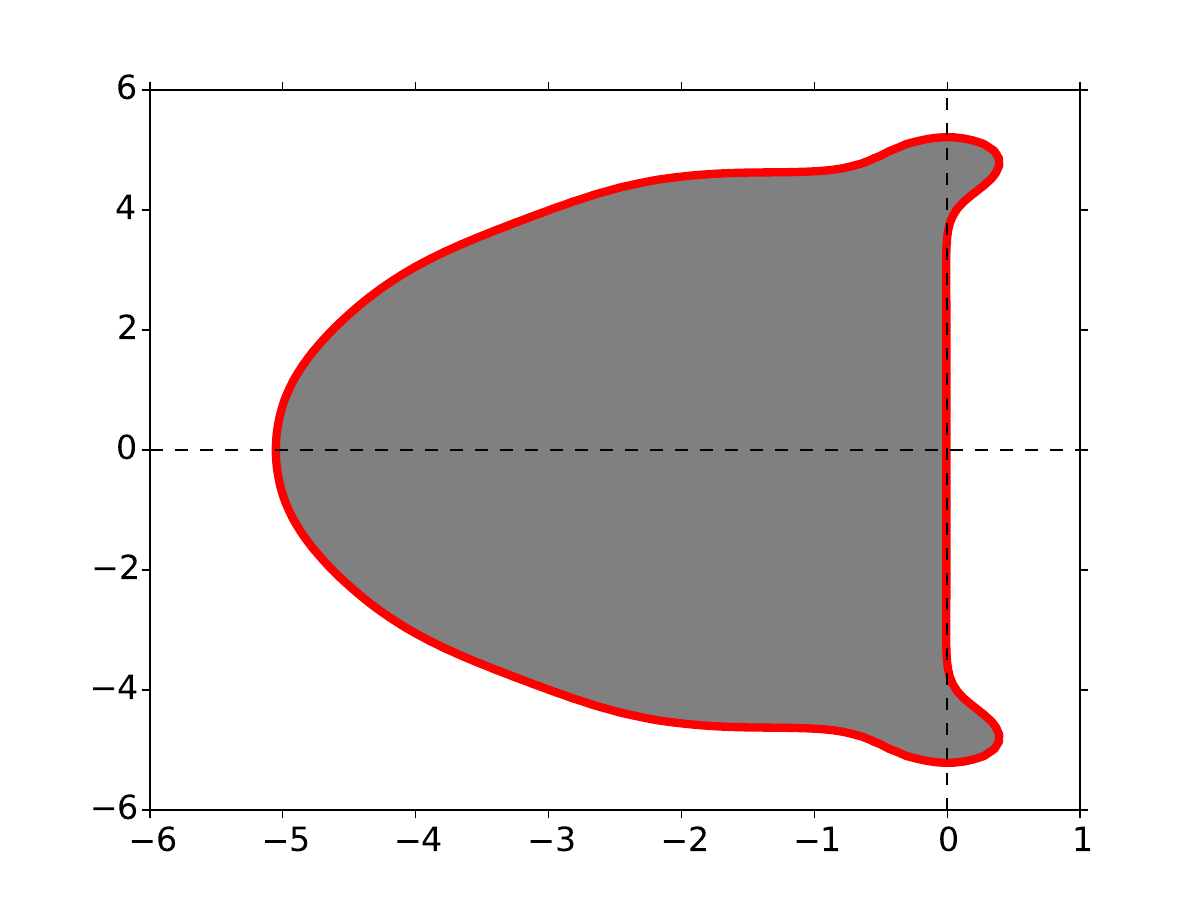}
		\caption{$K=10$}
	\end{subfigure}
	\caption{Linear stability regions of the truncated Taylor series approximation (without Borel summation)}
	\label{fig:region_series}
\end{figure}

We carry out two different linear stability investigations. The first one is when the solution is approximated by the Taylor series truncated at order $K$, without the Borel summation procedure, and the second one is when the solution is approximated with the BPL scheme. When the summation procedure is not applied, the method will be called time Asymptotic Numerical Method (ANM) as in computational solid and fluid mechanics \cite{cochelin07}. With ANM, we have:
\begin{equation}
	u(h)=\left(\sum_{k=0}^K\cfrac{(h\lambda)^k}{k!}\right)u_0.
\end{equation}
Comparing this relation to (\ref{iteration}), we define the domain of linear stability, for a given truncation order $K$, as 
\begin{equation}
	D^K_{ANM}=\left\{z\in\mathbb C\text{ such that } \left|\sum_{k=0}^K\cfrac{z^k}{k!}\right|\leq1\right\}.
\end{equation}
This domain is plotted in Figure \ref{fig:region_series} for $K$ from 2 to $10$. As can be observed, $D^K_{ANM}$ grows with $K$. The growth is however rather slow. Let us use the positive number $|D^K_{ANM}|$ defined as follows as a quantification of the size of $D^K_{ANM}$:
\begin{equation}
	|D^K_{ANM}| =\sup\big\{d\geq 0\ \text{such that}\ [-d,0] \in D^K_{ANM}\big\}.
	\label{size_d}
\end{equation}
This quantity increases almost linearly as can be seen in Figure \ref{fig:region_series}. The slope of the curve is about 0.375.
\begin{figure}
	\centering
	\includegraphics[width=\figwidth]{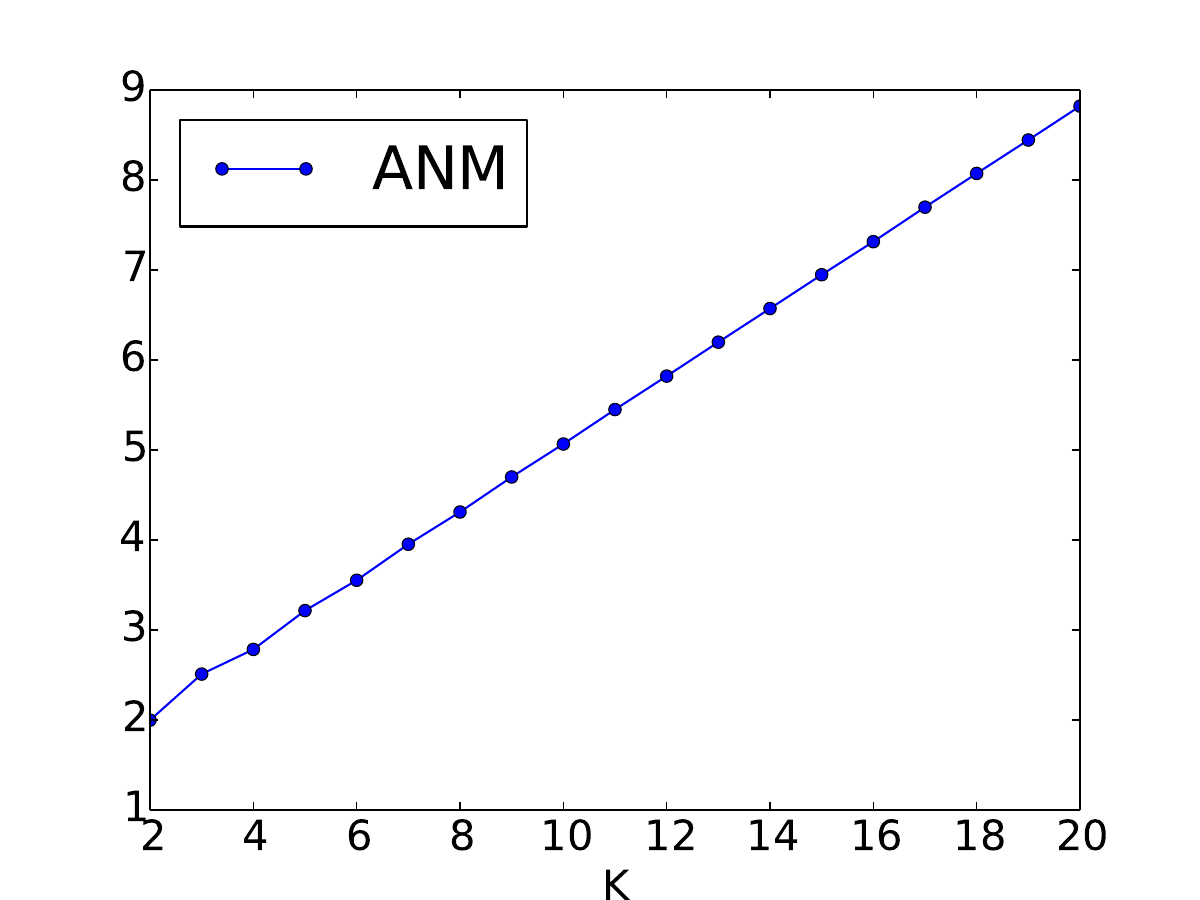}
	\caption{Size $|D_{ANM}^K|$ of the stability region when $K$ grows}
	\label{fig:region_min_anm}
\end{figure}

In fact, the region $D_{ANM}^K$ coincides with the stability region of an explicit $K$-th order Runge-Kutta method. The reason to this is that the stability function of an explicit Runge-Kutta method is a truncated Taylor expansion of the exponential function.

We now apply the summation procedure to the series. To make it more concrete, let us set $u_0=1$ and $K=4$. Like previously, the truncated series solution is
\begin{equation}
	\sum_{k=0}^4\frac{(\lambda t)^k}{k!}=1+\frac{λt}1+\frac {(λt)^2}2+\frac {(λt)^3}6+\frac {(λt)^4}{24}.
\end{equation}
Its Borel transform reads
\begin{equation}
	\sum_{k=0}^3\frac{\lambda ^{k+1}}{(k+1)!}\frac{\xi^k}{k!}=λ\left( 1+\frac{λ\xi}2+\frac {(λ\xi)^2}{12}+\frac {(λ\xi)^3}{144} \right)
	\label{borel_4}
\end{equation}
We take $k_a=1$ and $k_b=2$. The [1/2] Padé approximant of (\ref{borel_4}) is
\begin{equation}
	P(\xi)=λ\, \cfrac{48+14λ\xi}{48-10λ\xi+(λ\xi)^2}.
	\label{pade_exp}
\end{equation}
We then have the following approximate solution of the linear equation (\ref{linear}) with the Borel-Padé-Laplace scheme
\begin{equation}
	u(t)\simeq 1+t\sum_{i=1}^{N_G}P(t\xi_i)ω_i=1+λt\sum_{i=1}^{N_G}Q(λt\xi_i)ω_i
\end{equation}
where $Q$ is the rational function
\begin{equation}
	Q(z)=\cfrac{48+14z}{48-10z+z^2}.
\end{equation}
As already mentioned, $\xi_i$ is the $i$-th root of the $N_G$-th Laguerre polynomial and $ω_i$ is the corresponding weight in Gauss-Laguerre quadrature.
Note that the Padé approximant (\ref{pade_exp}) has no pole on the integration domain (the real positif axis) of the Laplace integral.

\begin{figure}
	\centering
	\begin{subfigure}{39mm}
		\includegraphics[width=\textwidth]{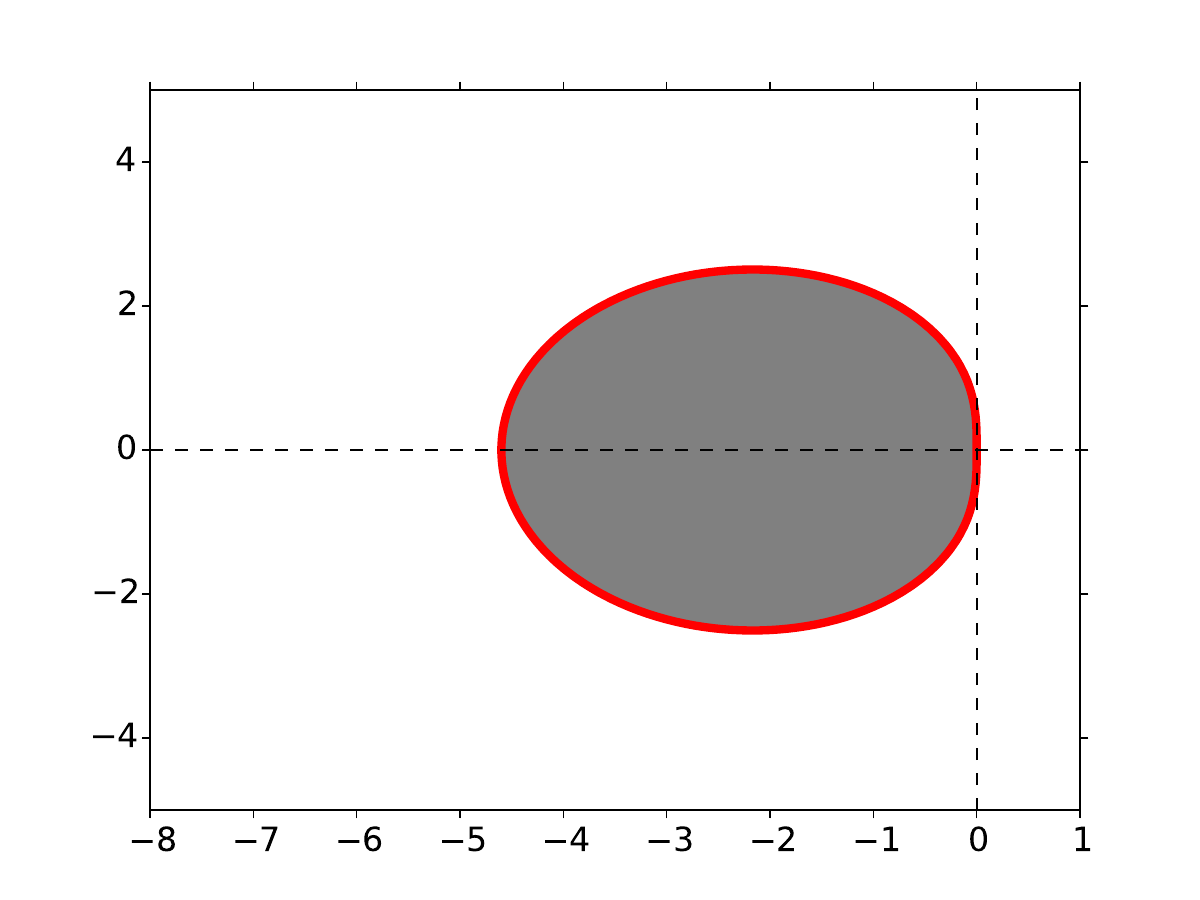}
		\caption{$K=2$, Padé $[0/1]$}\label{region_bpl_2_0}
	\end{subfigure}
	\begin{subfigure}{39mm}
		\includegraphics[width=\textwidth]{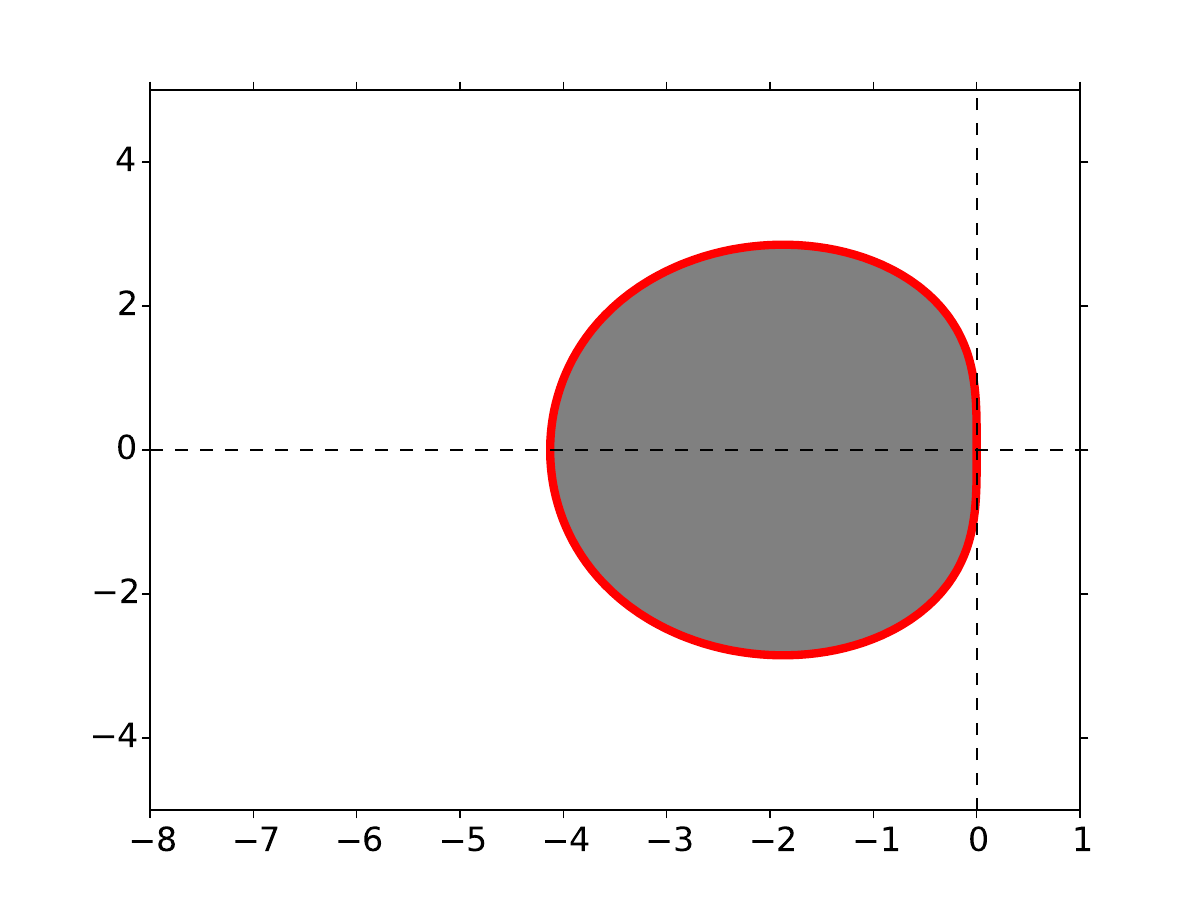}
		\caption{$K=3$, Padé $[1/1]$}
	\end{subfigure}
	\begin{subfigure}{39mm}
		\includegraphics[width=\textwidth]{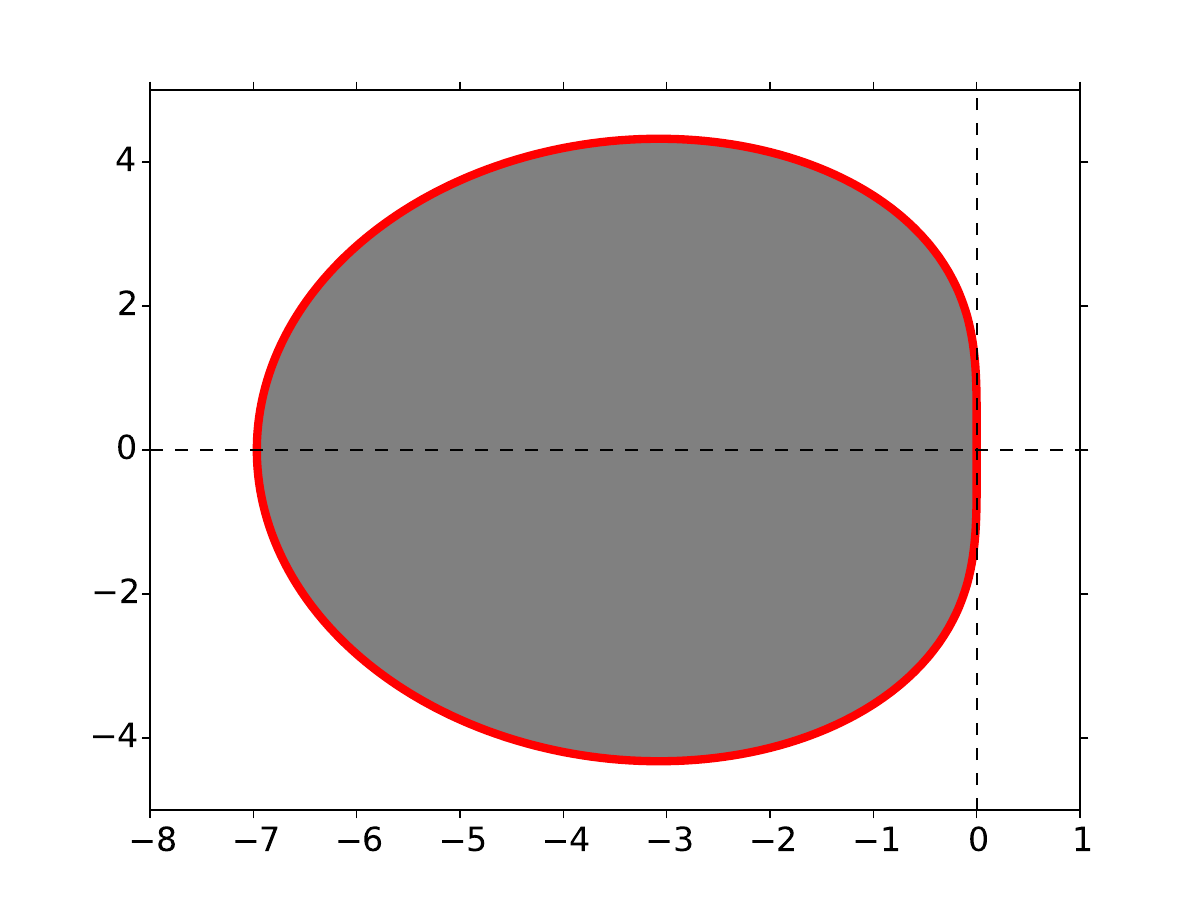}
		\caption{$K=4$, Padé $[1/2]$}
	\end{subfigure}%

	\vspace{\baselineskip}
	\begin{subfigure}{39mm}
		\includegraphics[width=\textwidth]{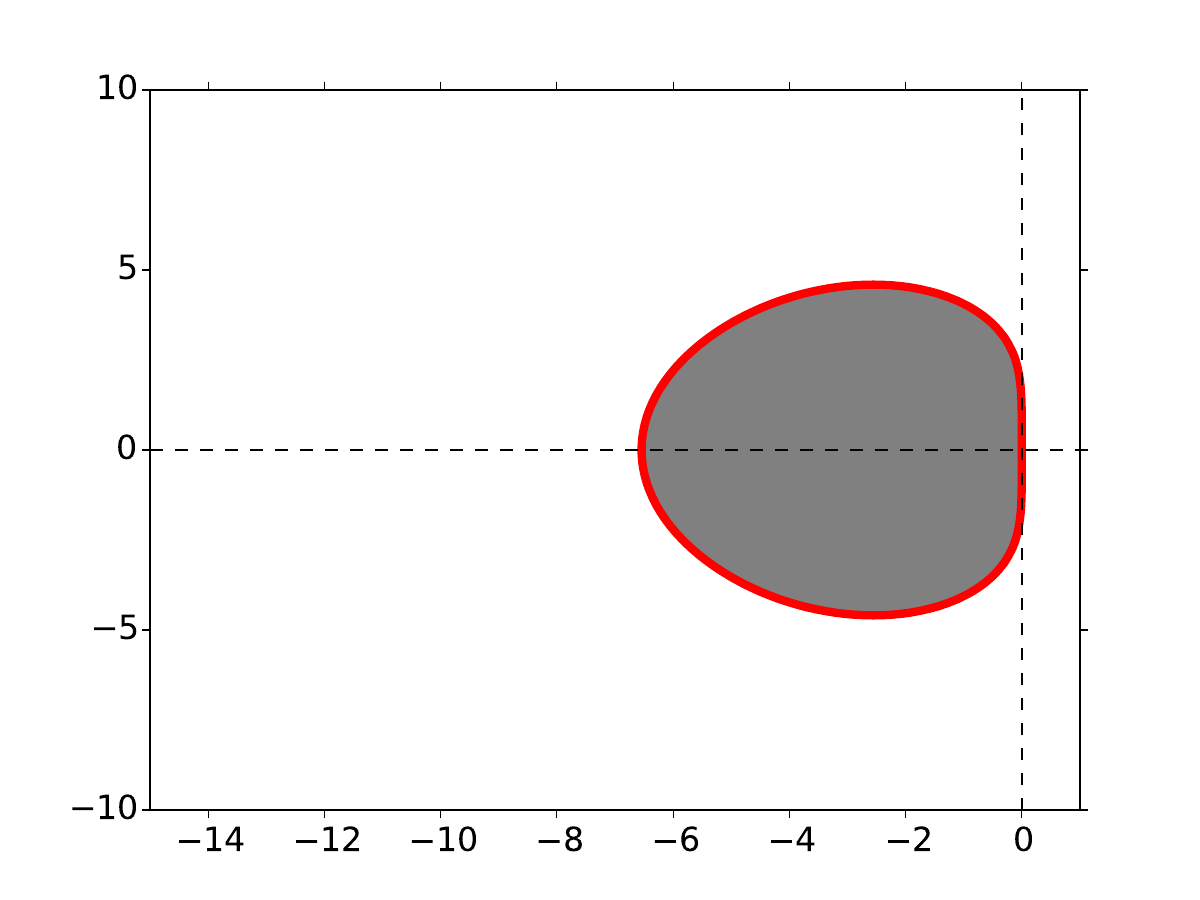}
		\caption{$K=5$, Padé $[2/2]$}
	\end{subfigure}
	\begin{subfigure}{39mm}
		\includegraphics[width=\textwidth]{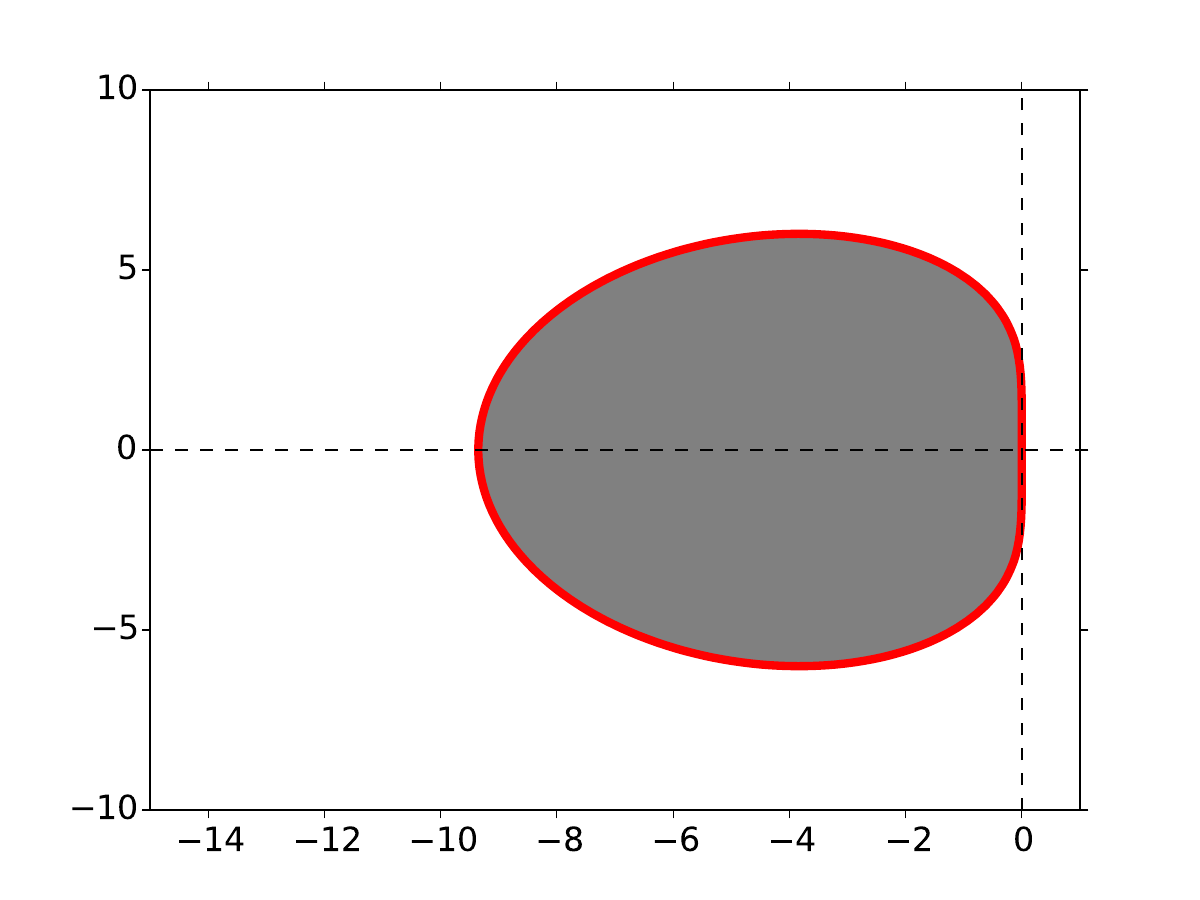}
		\caption{$K=6$, Padé $[2/3]$}
	\end{subfigure}
	\begin{subfigure}{39mm}
		\includegraphics[width=\textwidth]{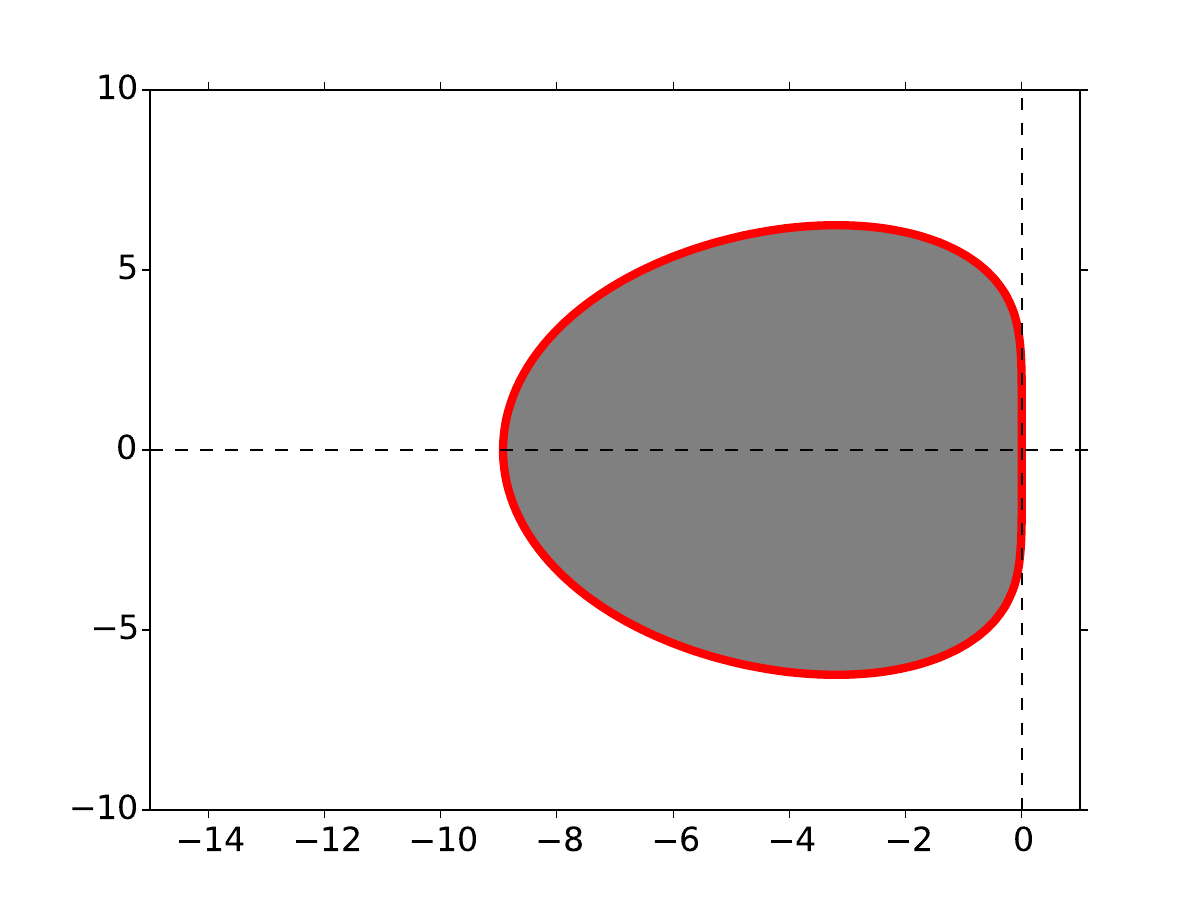}
		\caption{$K=7$, Padé $[3/3]$}
	\end{subfigure}

	\vspace{\baselineskip}
	\begin{subfigure}{39mm}
		\includegraphics[width=\textwidth]{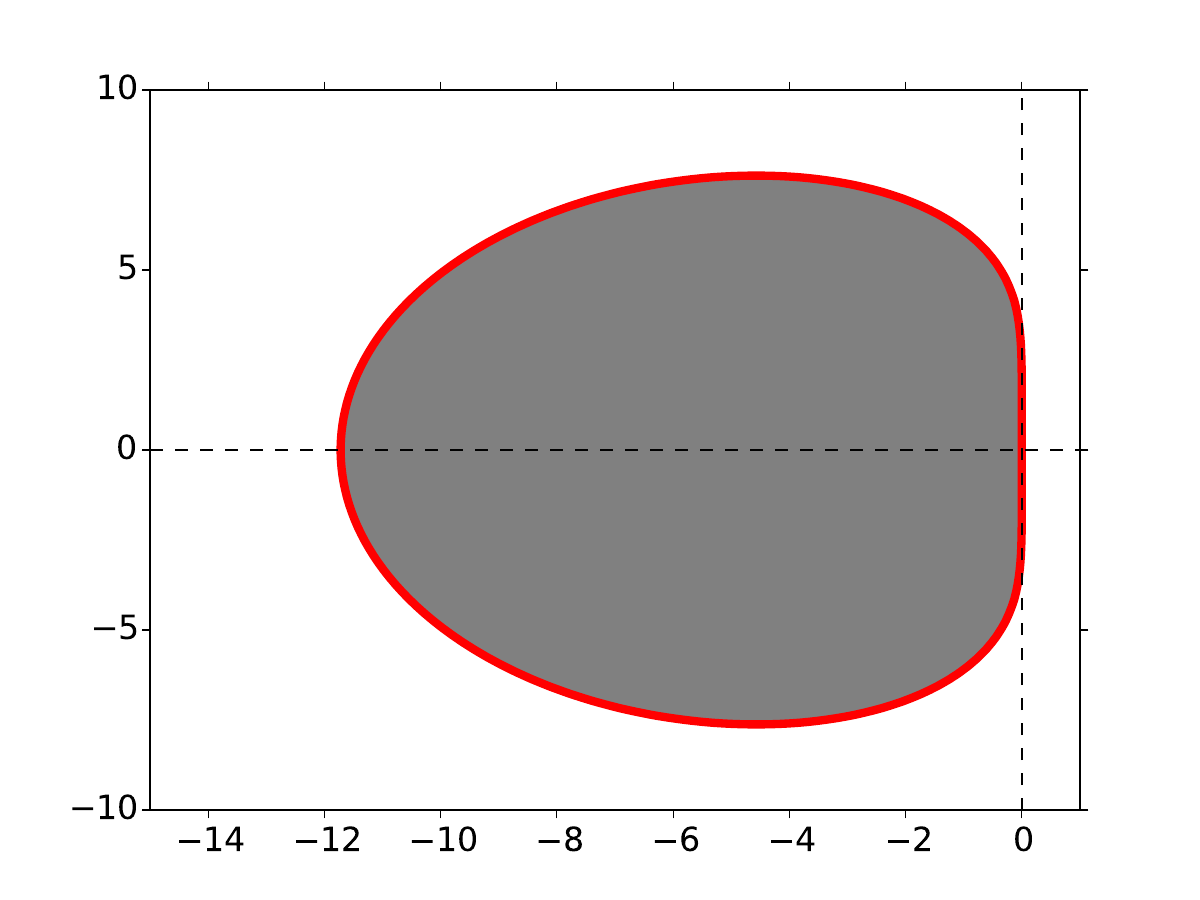}
		\caption{$K=8$, Padé $[3/4]$}
	\end{subfigure}
	\begin{subfigure}{39mm}
		\includegraphics[width=\textwidth]{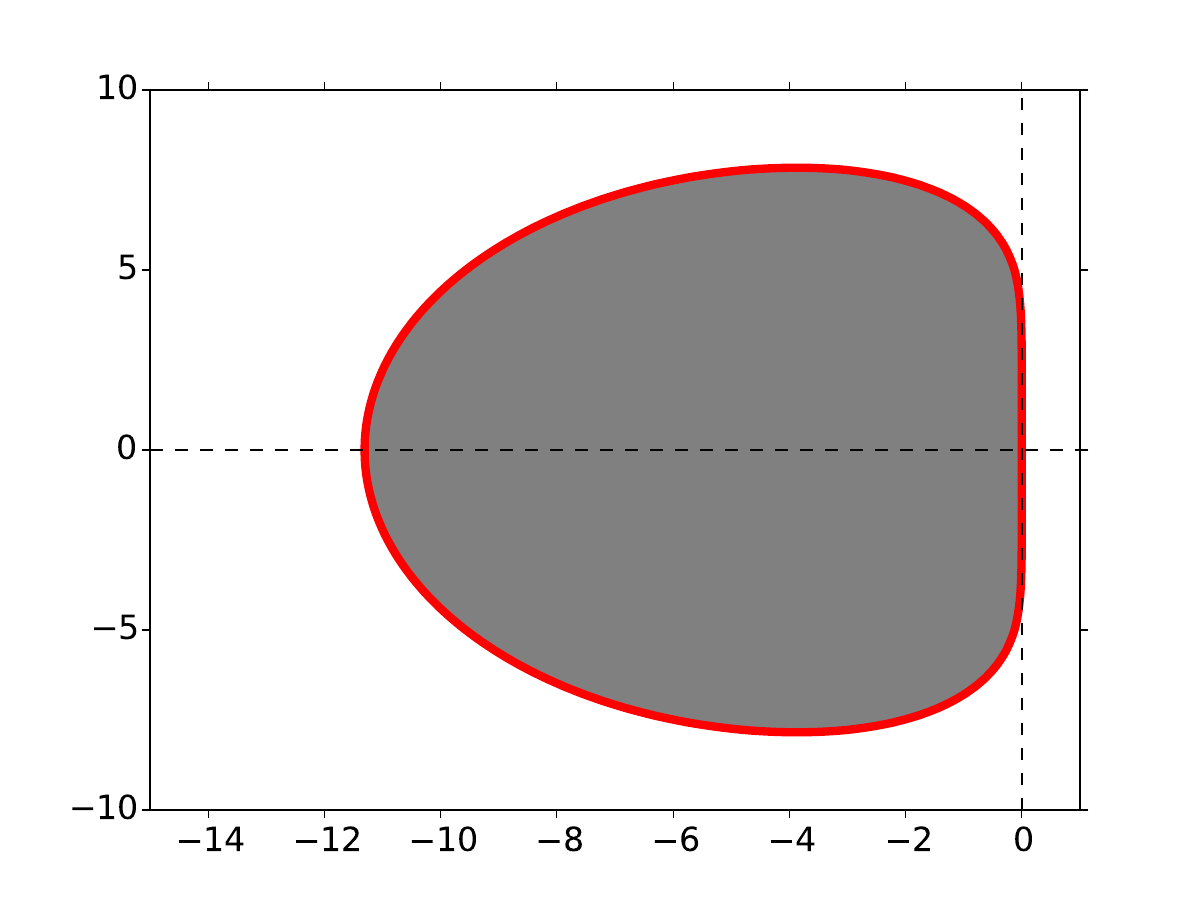}
		\caption{$K=9$, Padé $[4/4]$}
	\end{subfigure}
	\begin{subfigure}{39mm}
		\includegraphics[width=\textwidth]{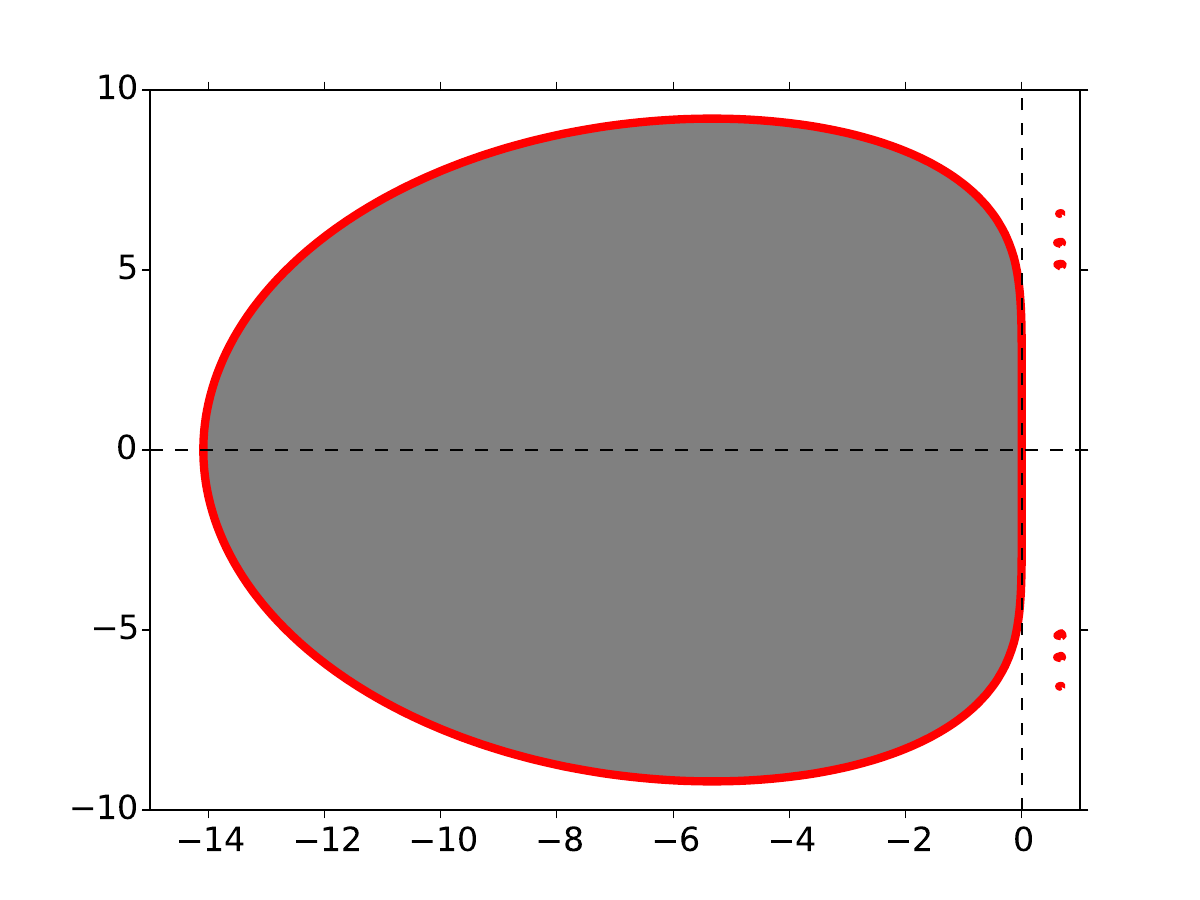}
		\caption{$K=10$, Padé $[4/5]$}
	\end{subfigure}
	\caption{Linear stability regions of Borel-Padé-Laplace integrator with increasing $K$ and close-to-diagonal Padé approximants.}
	\label{fig:region_BPL}
\end{figure}

The stability region of Borel-Padé-Laplace integrator, for $K=4$, is 
\begin{equation}
	D_{BPL}^4=\left\{z\in\mathbb C\text{ such that } \left|1+z\sum_{i=0}^{N_G}Q(z\xi_i)ω_i\right|\leq1\right\}.
\end{equation}
This region is plotted in Figure \ref{fig:region_BPL}, with $N_G=100$ Gauss points, along with the stability regions for other values of $K$, ranging from 2 to 10. In this Figure, the Padé approximant in Borel space is chosen as closed as possible to the diagonal, that is
\begin{equation}
	K_b =K_a\quad\quad\text{or}\quad\quad K_b=K_a+1
	\label{ka}
\end{equation}
depending on the parity of $K$. This choice will be discussed later.
\begin{figure}
	\centering
	\includegraphics[width=\figwidth]{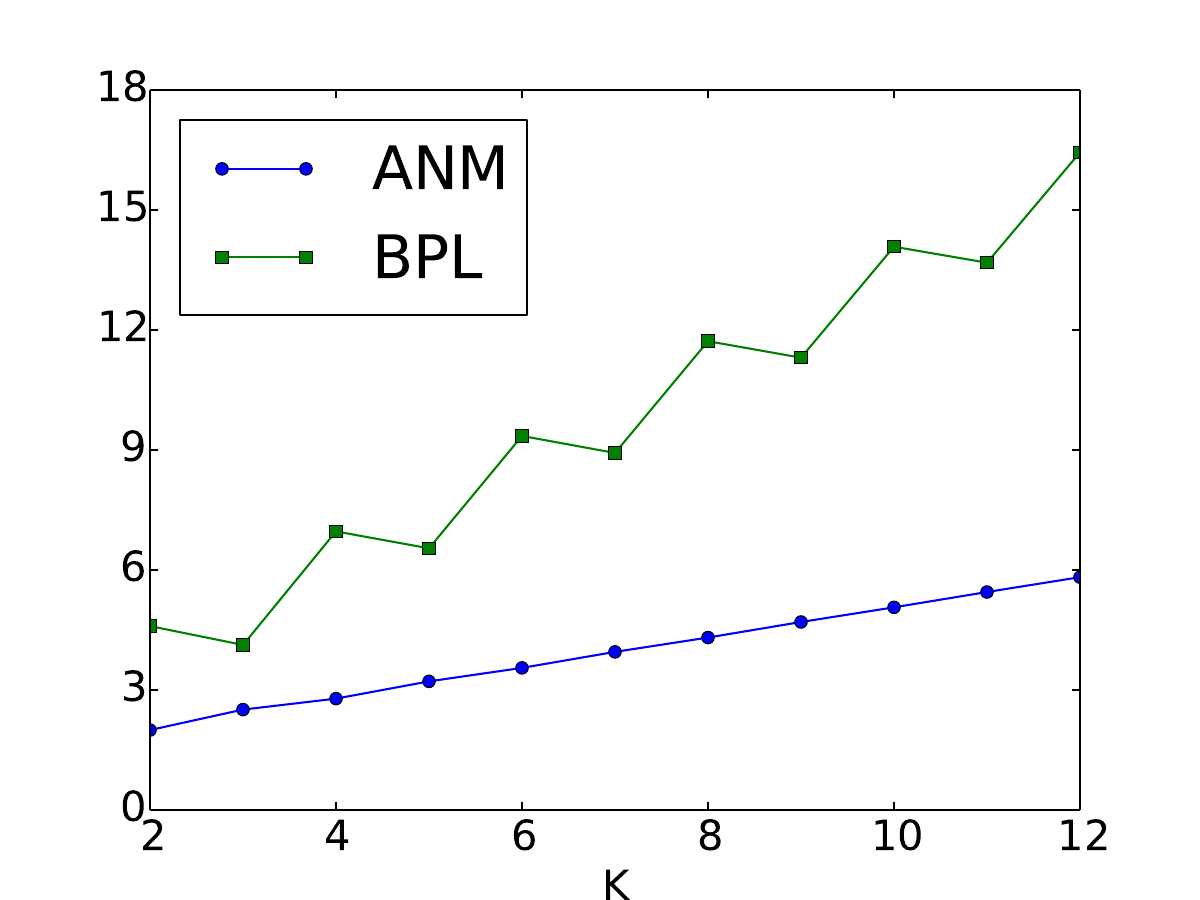}
	\caption{Evolution of $|D_{ANM}^K|$ and $|D_{BPL}^K|$ with $K$}
	\label{fig:region_min}
\end{figure}

As can be seen in Figure \ref{fig:region_BPL}, the regions do not include the half complex plane with negative real part. Indeed, as an explicit scheme, BPL is not $A$-stable. However, this figure clearly shows that, for a fixed $K\geq4$, the stability region of BPL is much wider than that of the simple truncated series scheme or that of an explicit Runge-Kutta scheme. Note that when $K=2$, the Borel summation has no effect, and the stability region is the same as in Figure  \ref{fig:region_series_K2}. 

Another striking point is that the growth of the stability region with $K$ is not as regular as in the case where the Borel summation is not applied. This is due to the choice of (almost) diagonal Padé approximant. However, if we consider either only odd $K$ or only even $K$, the growth is regular again. In all cases, the overall growth rate is higher than in Figure  (\ref{fig:region_series_K2}). Indeed, let us quantify the size of $D^K_{BPL}$ as previously with
\begin{equation}
	|D^K_{BPL}| =\sup\big\{d\geq 0\ \text{such that}\ [-d,0] \in D^K_{BPL}\big\}
\end{equation}
The dependence of $|D^K_{BPL}|$ on $K$ is plotted in Figure \ref{fig:region_min} (along with the evolution of $|D^K_{ANM}|$ for comparison), for $K$ ranging from 2 to 12. The average slope of the curve of $|D^K_{BPL}|$ is about 0.543.

Since we are in the particular situation where the Taylor coefficients $u_k$ of the solution decrease very rapidly with $k$, only few terms of the series are numerically meaningful. More precisely, the coefficients $B_k$ of the Borel transformed series $\mathcal{B}\u$ are below our machine precision (about $2·10^{-16}$) for $k>11$. So, taking $K\geq 12$ does not bring any substantial improvement.

\begin{figure}
	\centering
	\begin{subfigure}{39mm}
		\includegraphics[width=\textwidth]{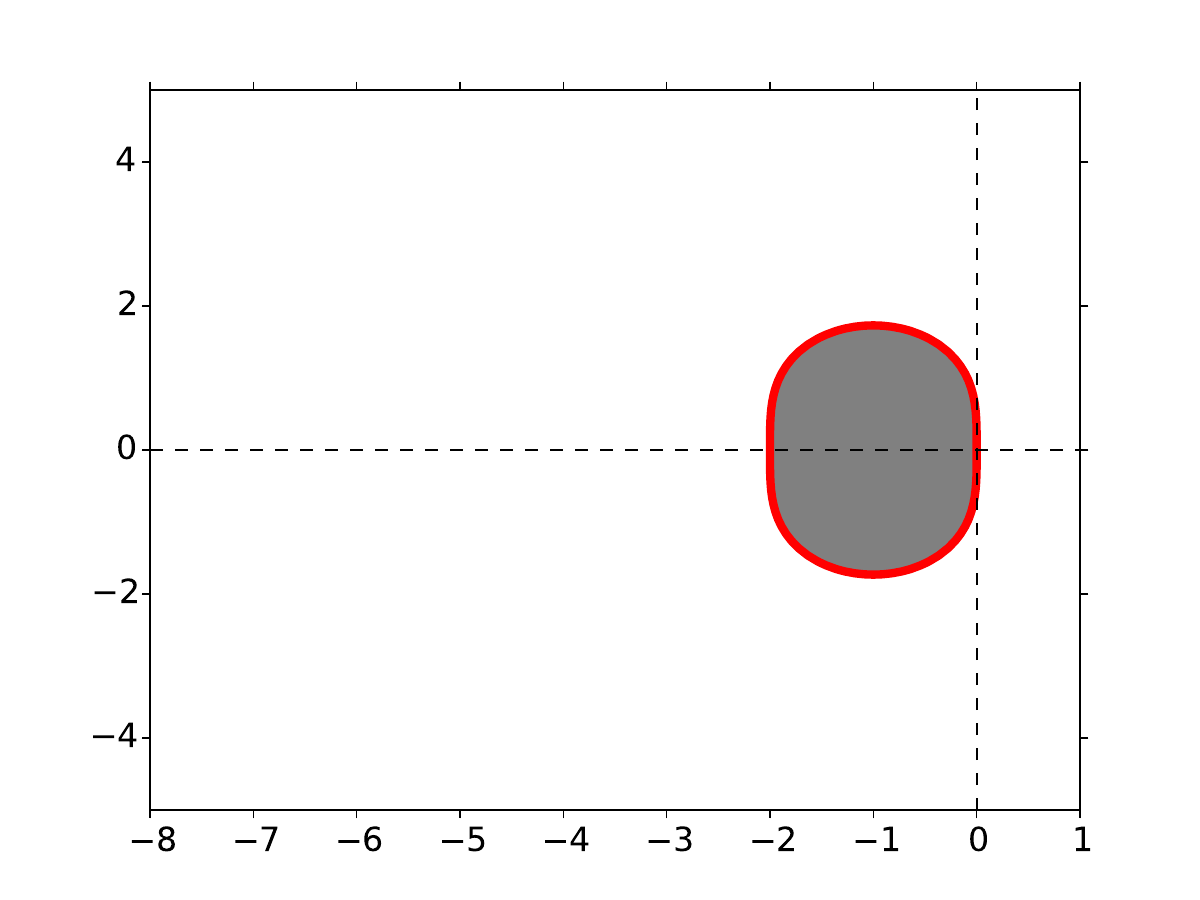}
		\caption{$K=2$, Padé $[1/0]$}
	\end{subfigure}
	\begin{subfigure}{39mm}
		\includegraphics[width=\textwidth]{regions_region_bpl_3_1}
		\caption{$K=3$, Padé $[1/1]$}
	\end{subfigure}
	\begin{subfigure}{39mm}
		\includegraphics[width=\textwidth]{regions_region_bpl_4_1}
		\caption{$K=4$, Padé $[1/2]$}
	\end{subfigure}%

	\vspace{\baselineskip}
	\begin{subfigure}{39mm}
		\includegraphics[width=\textwidth]{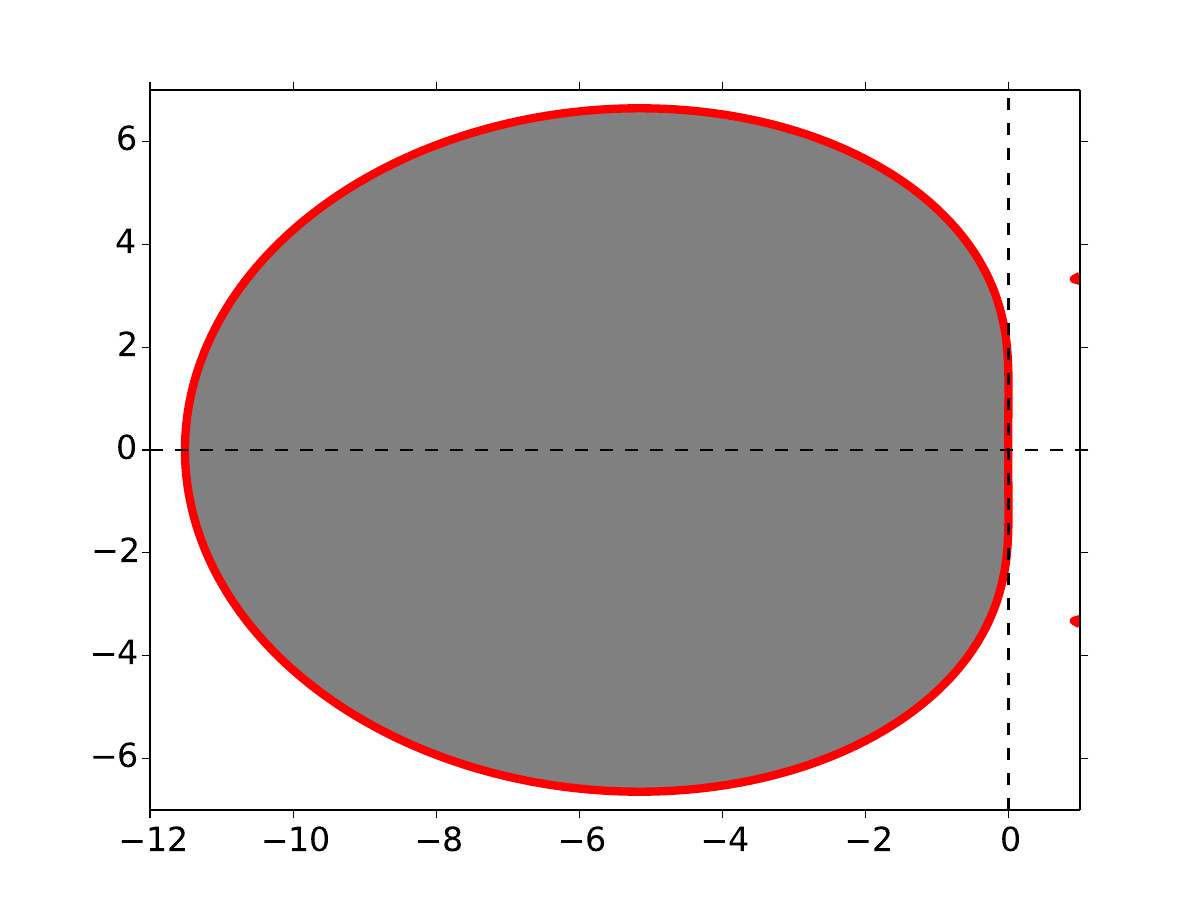}
		\caption{$K=5$, Padé $[1/3]$}
	\end{subfigure}
	\begin{subfigure}{39mm}
		\includegraphics[width=\textwidth]{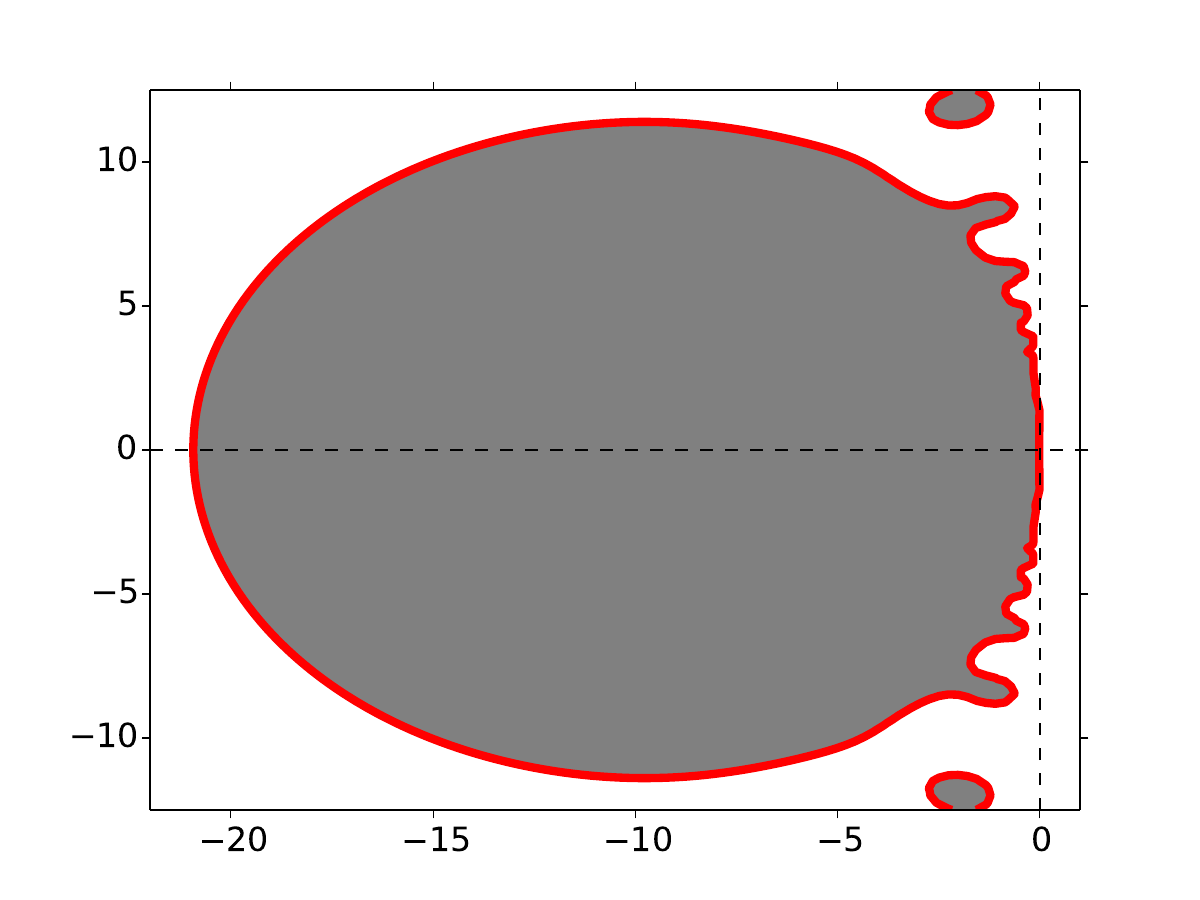}
		\caption{$K=6$, Padé $[1/4]$}
	\end{subfigure}
	\begin{subfigure}{39mm}
		\includegraphics[width=\textwidth]{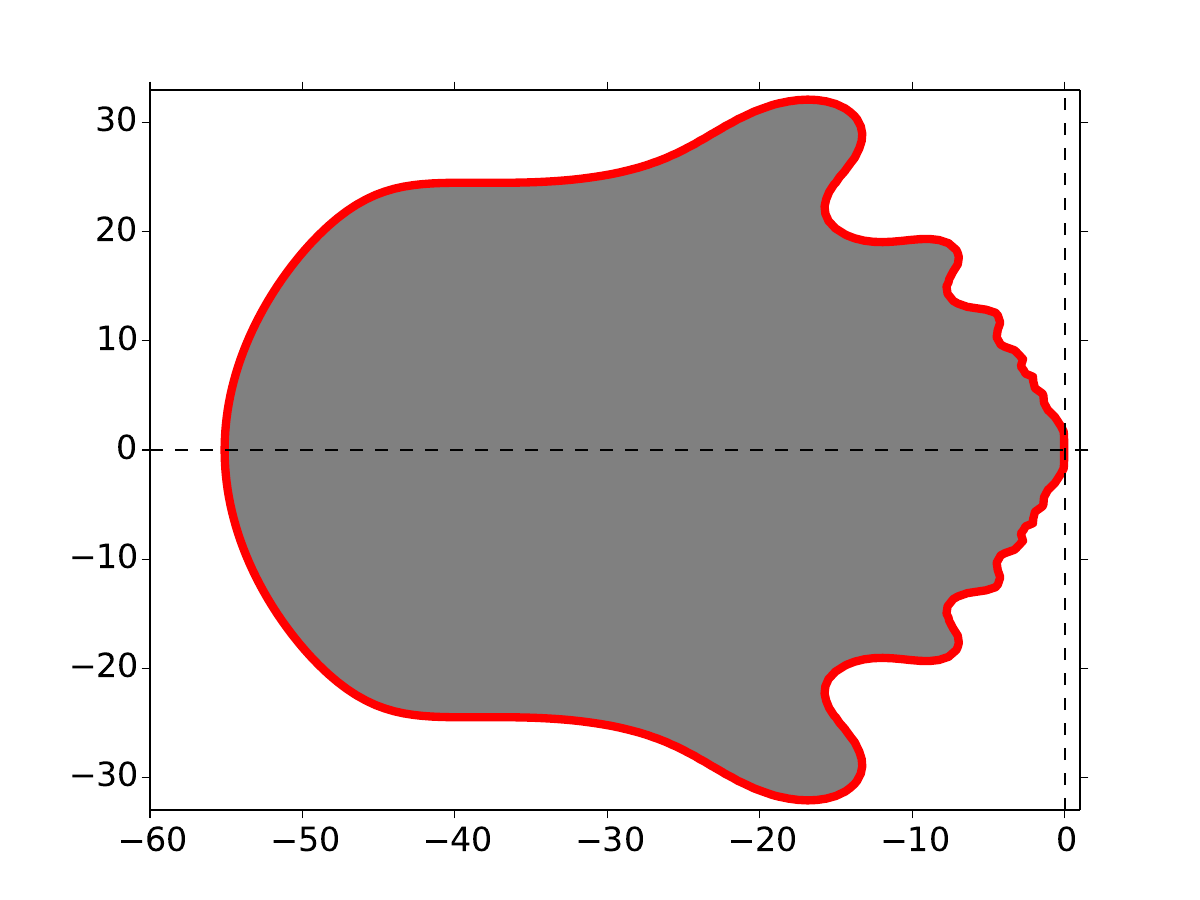}
		\caption{$K=7$, Padé $[1/5]$}
	\end{subfigure}

	\vspace{\baselineskip}
	\begin{subfigure}{39mm}
		\includegraphics[width=\textwidth]{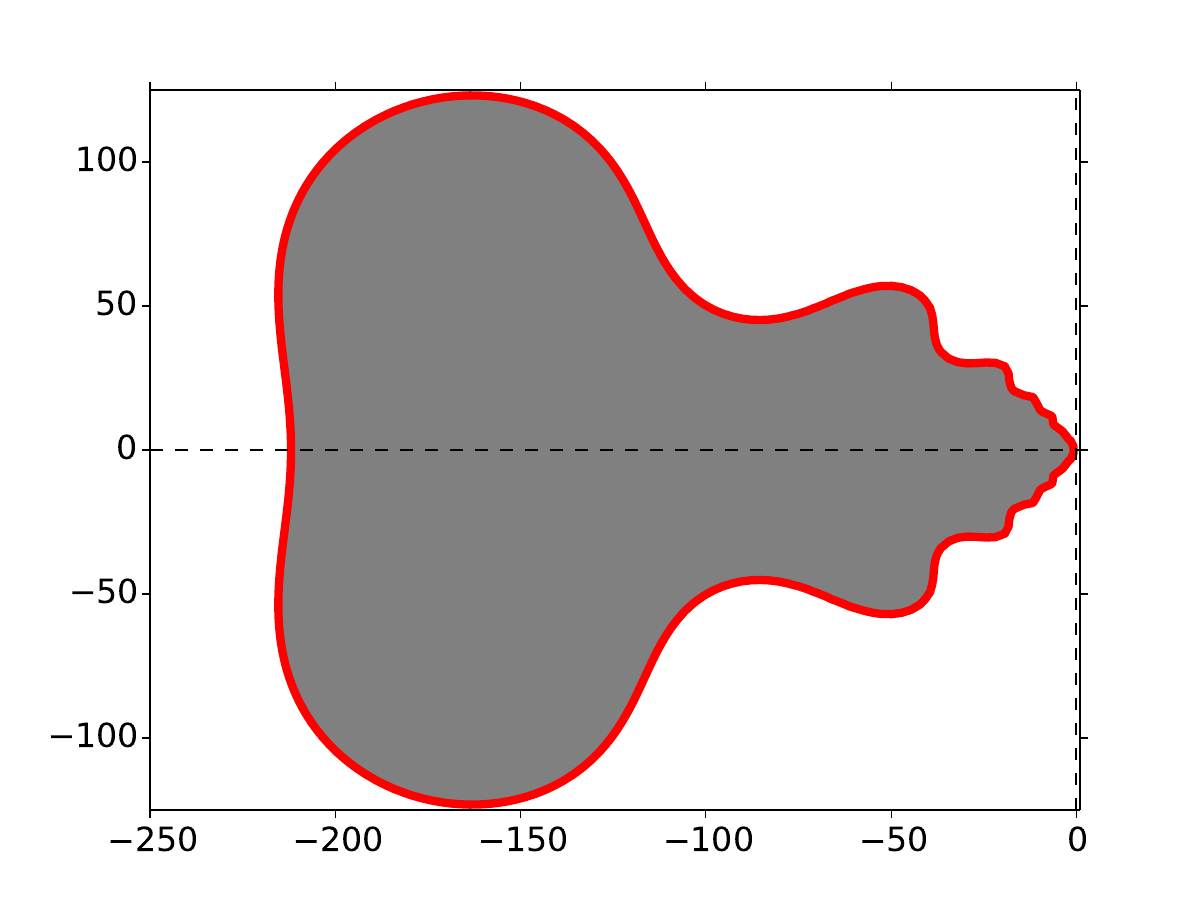}
		\caption{$K=8$, Padé $[1/6]$}
	\end{subfigure}
	\begin{subfigure}{39mm}
		\includegraphics[width=\textwidth]{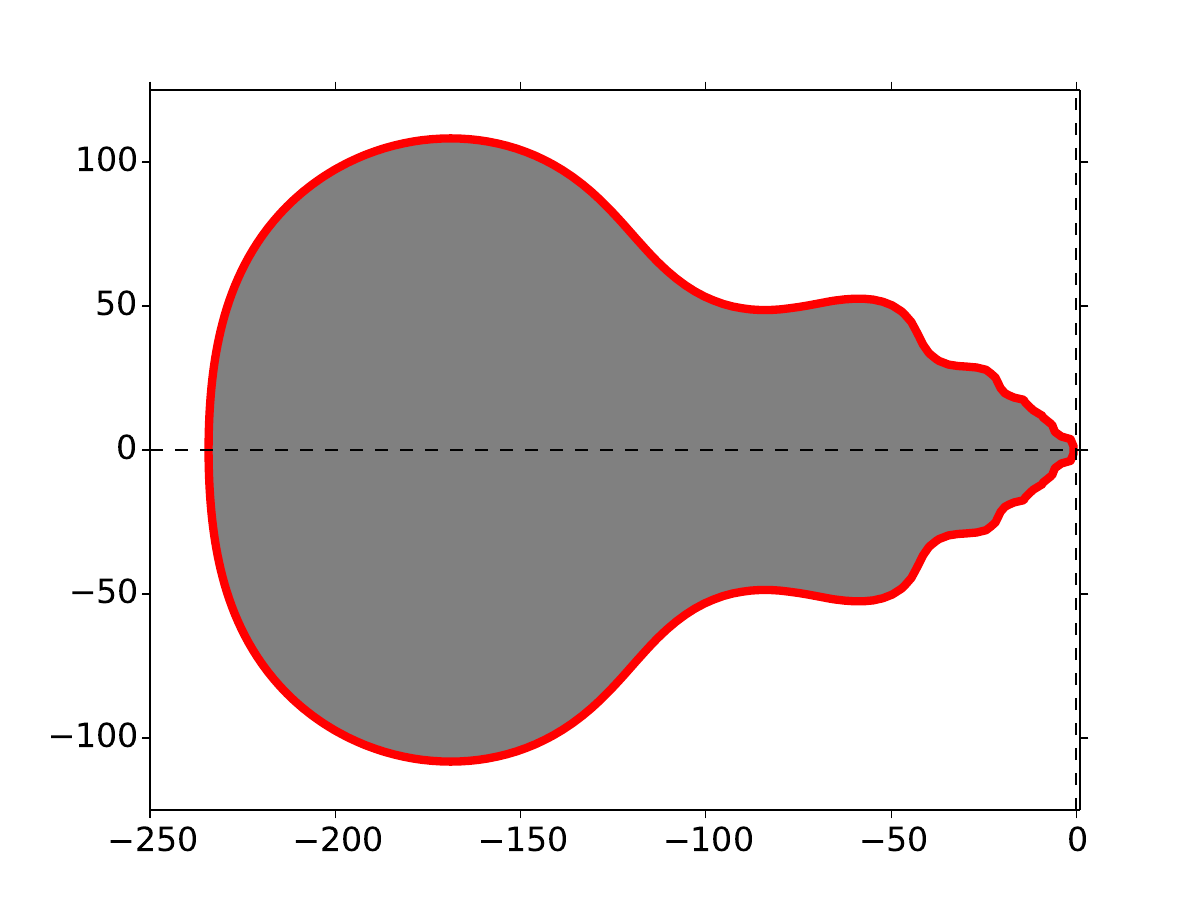}
		\caption{$K=9$, Padé $[1/7]$}
	\end{subfigure}
	\begin{subfigure}{39mm}
		\includegraphics[width=\textwidth]{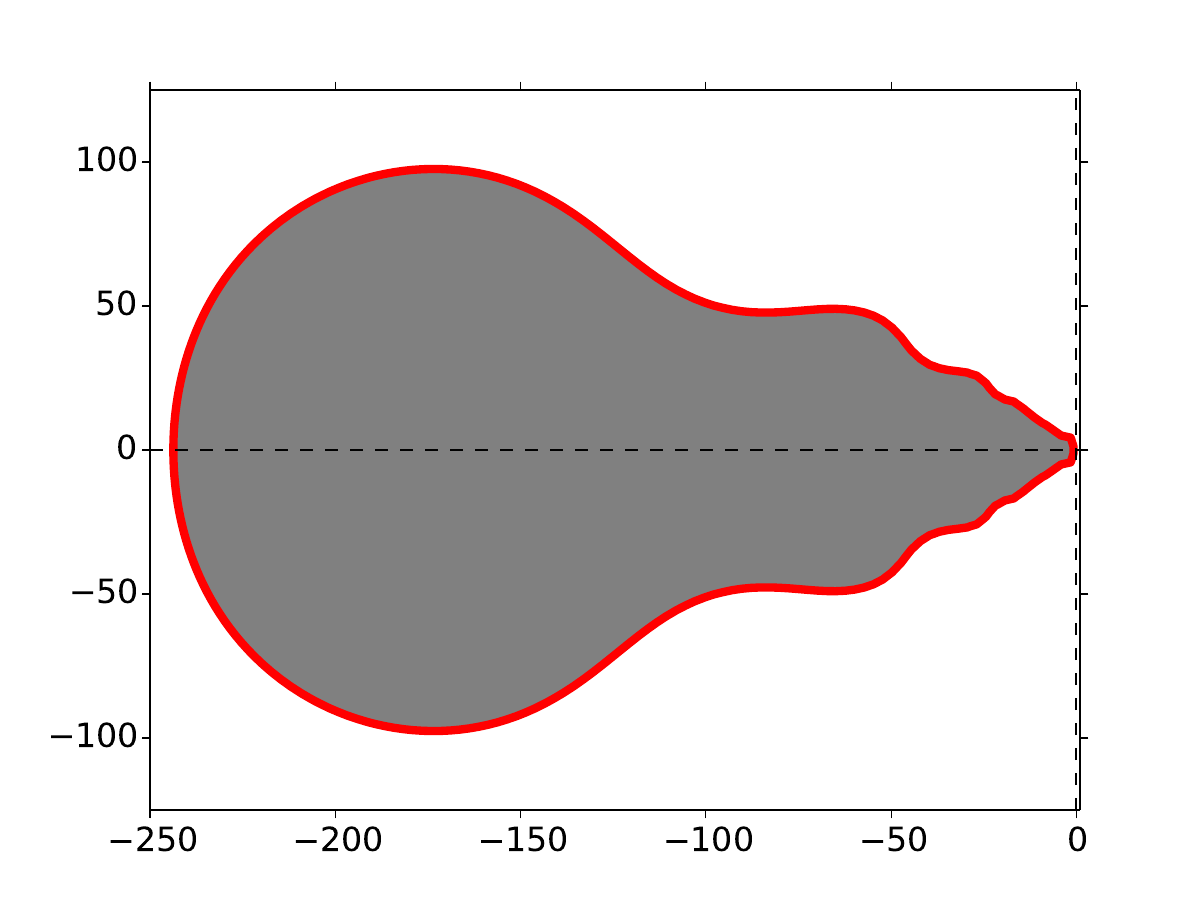}
		\caption{$K=10$, Padé $[1/9]$}
	\end{subfigure}
	\caption{Linear stability regions of Borel-Padé-Laplace integrator with increasing $K$ and fixed $K_a=1$}
	\label{fig:region_BPL_b}
\end{figure}

These observations indicate the importance of the Borel summation procedure, even in a situation where it has not been developed for. Indeed, the summation procedure enlarges very significantly the stability region, even when the Taylor series of the solution is convergent, with an infinite radius of convergence. This stability region can even be larger if we play with the parameters of the Padé approximants in Borel space. For example, we plot in Figure \ref{fig:region_BPL_b} the evolution of $D^K_{BPL}$ with $K$ when $K_a$ is fixed to 1. As can be observed, the stability domain grows with $K$ and is almost always far larger compared to Figures \ref{fig:region_series} and \ref{fig:region_BPL}. The growth rate is also far higher.

\begin{figure}
	\centering
	\begin{subfigure}{39mm}
		\includegraphics[width=\textwidth]{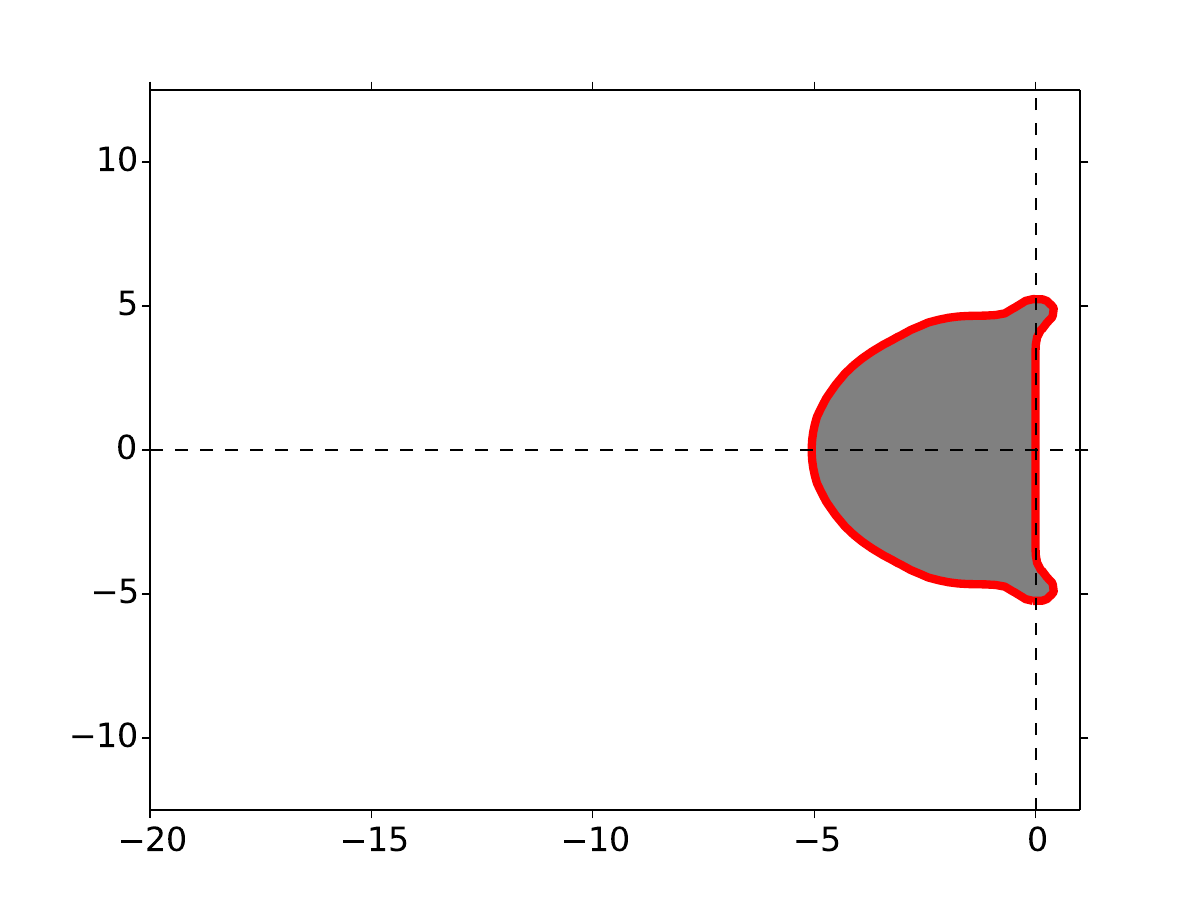}
		\caption{$K=10$, Padé $[9/0]$}
	\end{subfigure}
	\begin{subfigure}{39mm}
		\includegraphics[width=\textwidth]{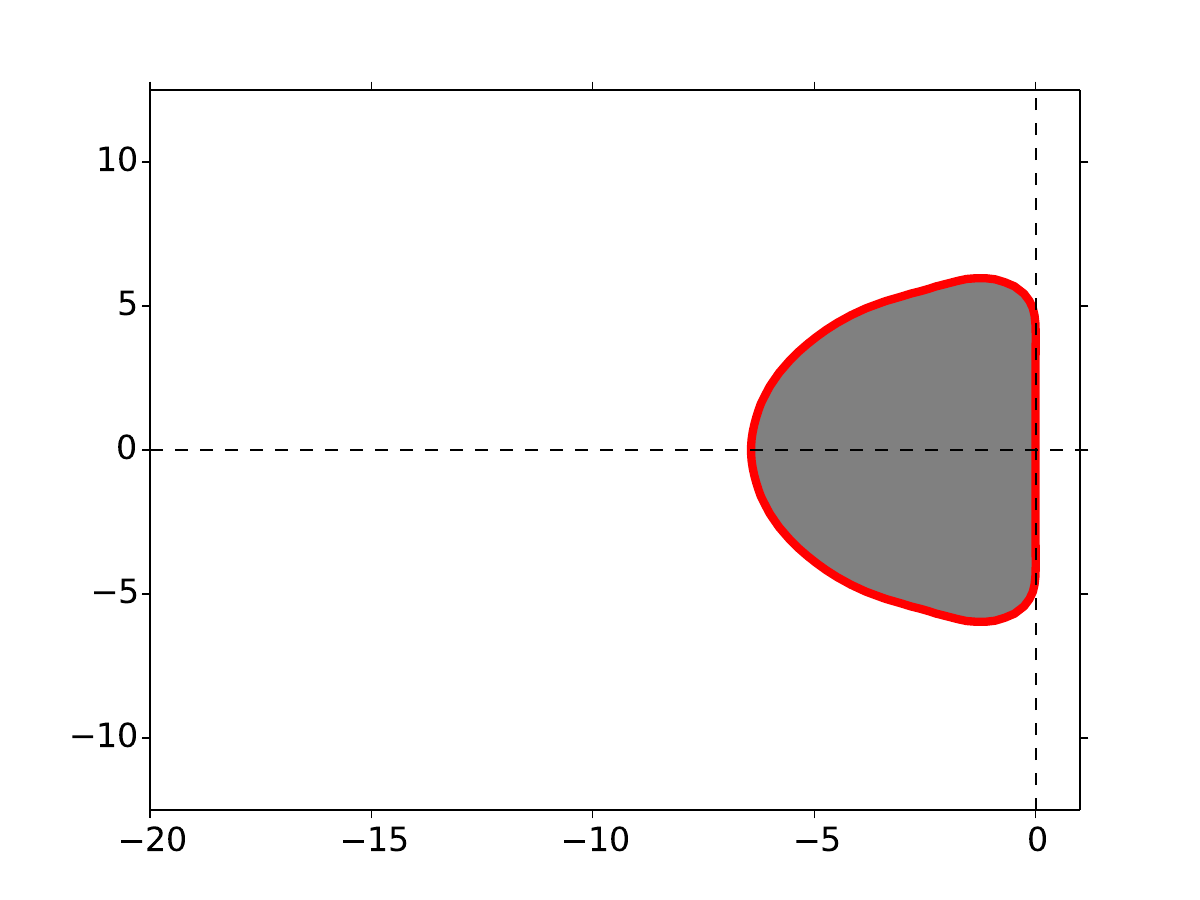}
		\caption{$K=10$, Padé $[8/1]$}
	\end{subfigure}
	\begin{subfigure}{39mm}
		\includegraphics[width=\textwidth]{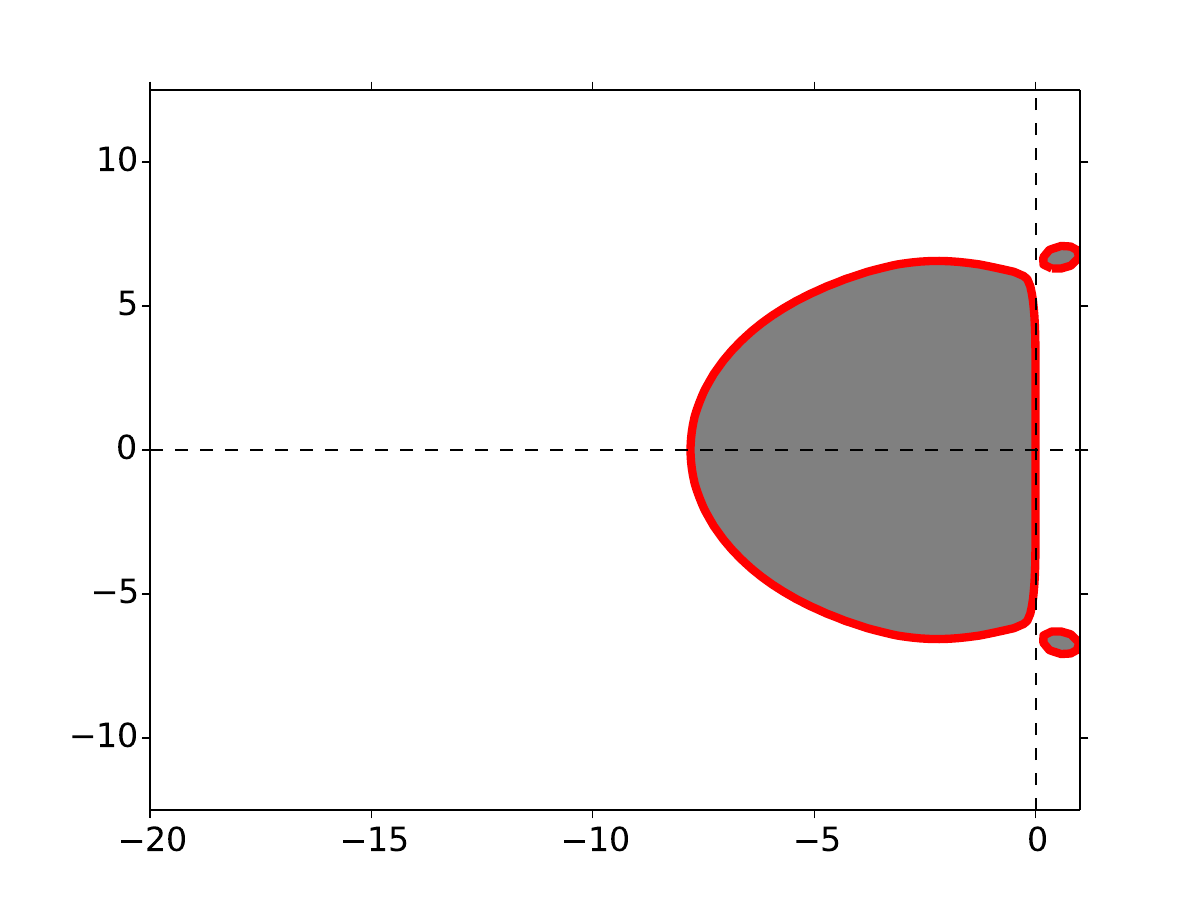}
		\caption{$K=10$, Padé $[7/2]$}
	\end{subfigure}%

	\vspace{\baselineskip}
	\begin{subfigure}{39mm}
		\includegraphics[width=\textwidth]{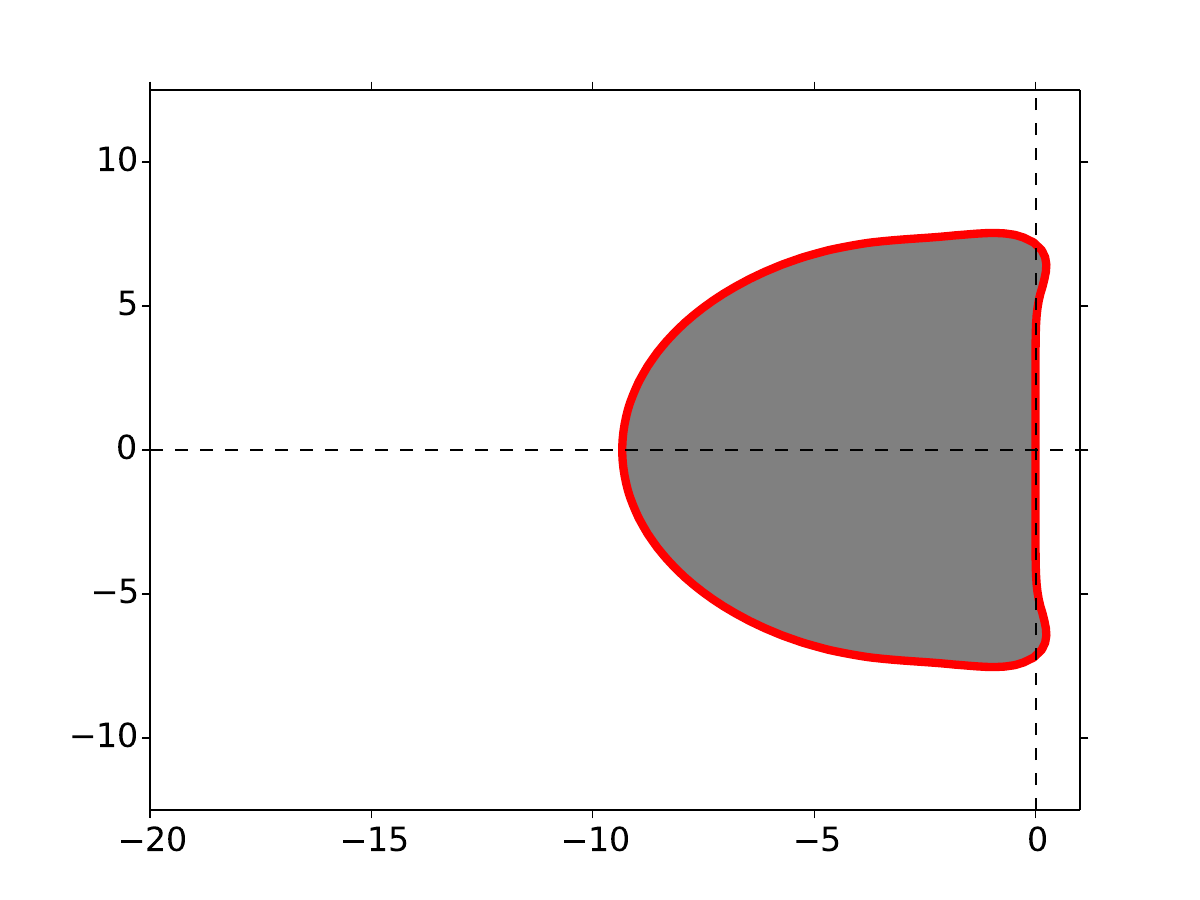}
		\caption{$K=10$, Padé $[6/3]$}
	\end{subfigure}
	\begin{subfigure}{39mm}
		\includegraphics[width=\textwidth]{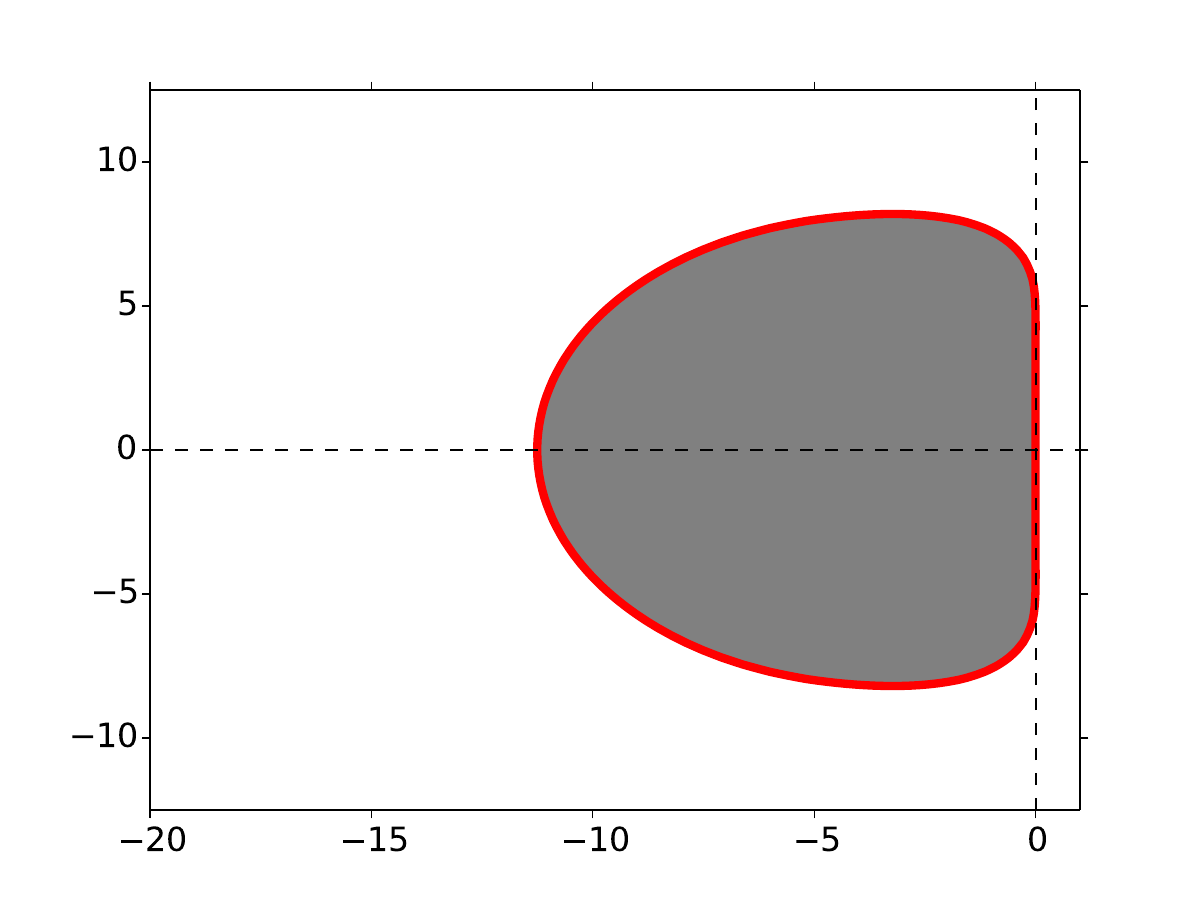}
		\caption{$K=10$, Padé $[5/4]$}
	\end{subfigure}
	\begin{subfigure}{39mm}
		\includegraphics[width=\textwidth]{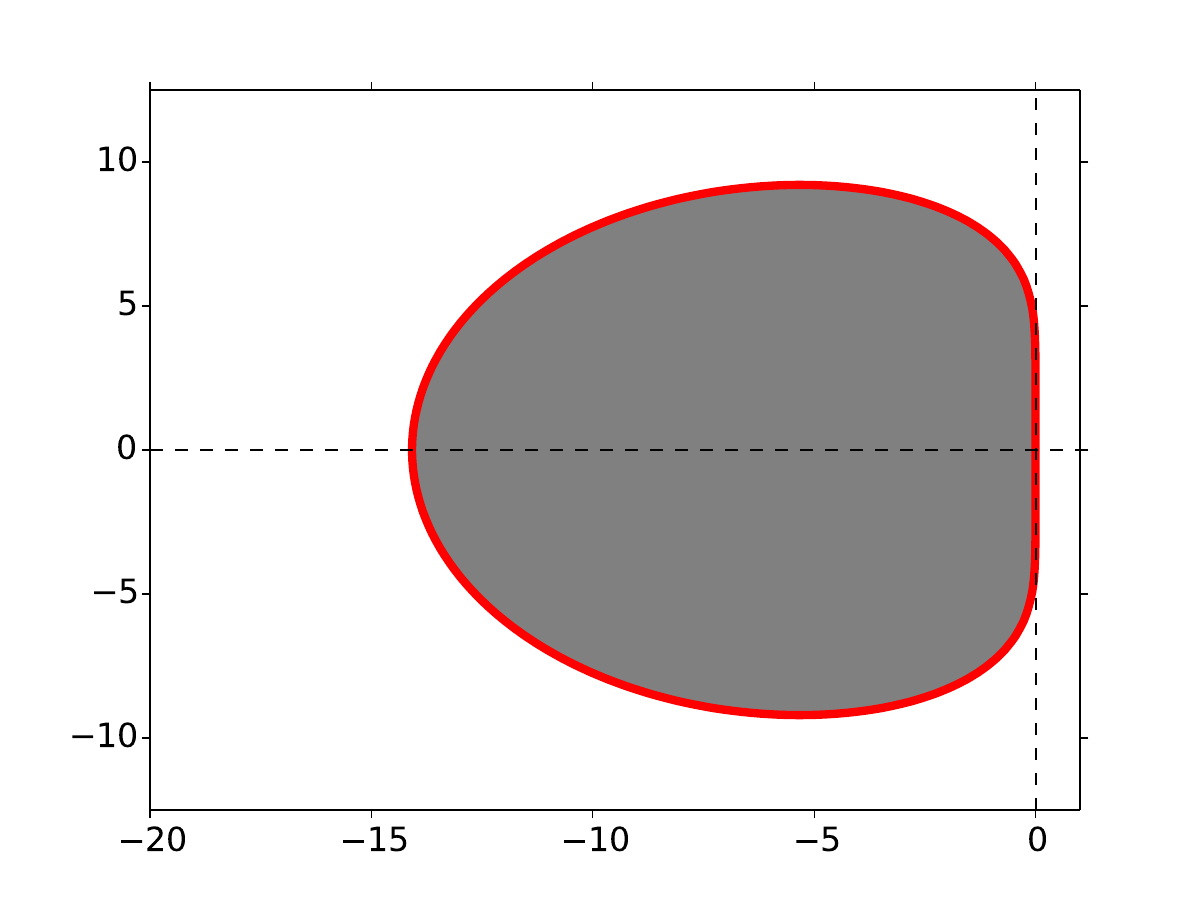}
		\caption{$K=10$, Padé $[4/5]$}
	\end{subfigure}

	\vspace{\baselineskip}
	\begin{subfigure}{39mm}
		\includegraphics[width=\textwidth]{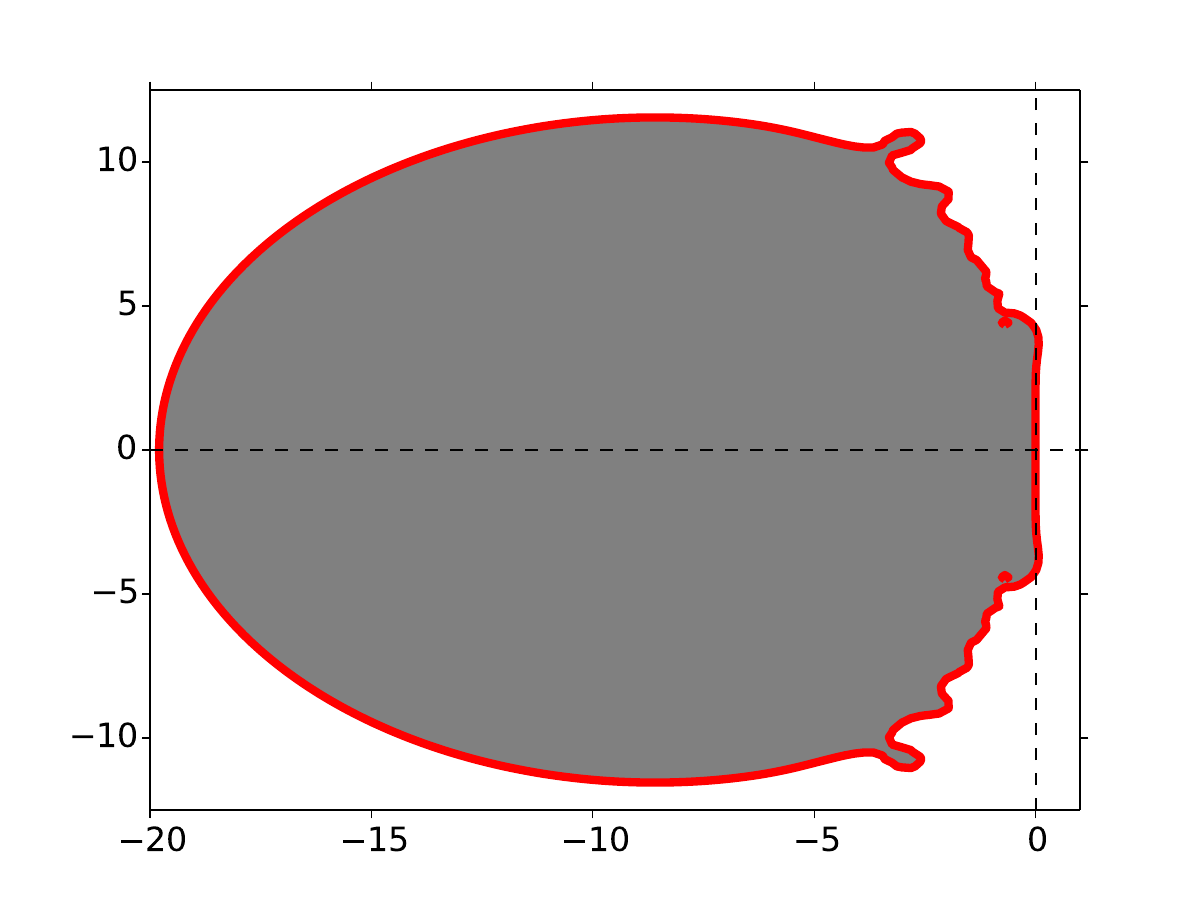}
		\caption{$K=10$, Padé $[3/6]$}
	\end{subfigure}
	\begin{subfigure}{39mm}
		\includegraphics[width=\textwidth]{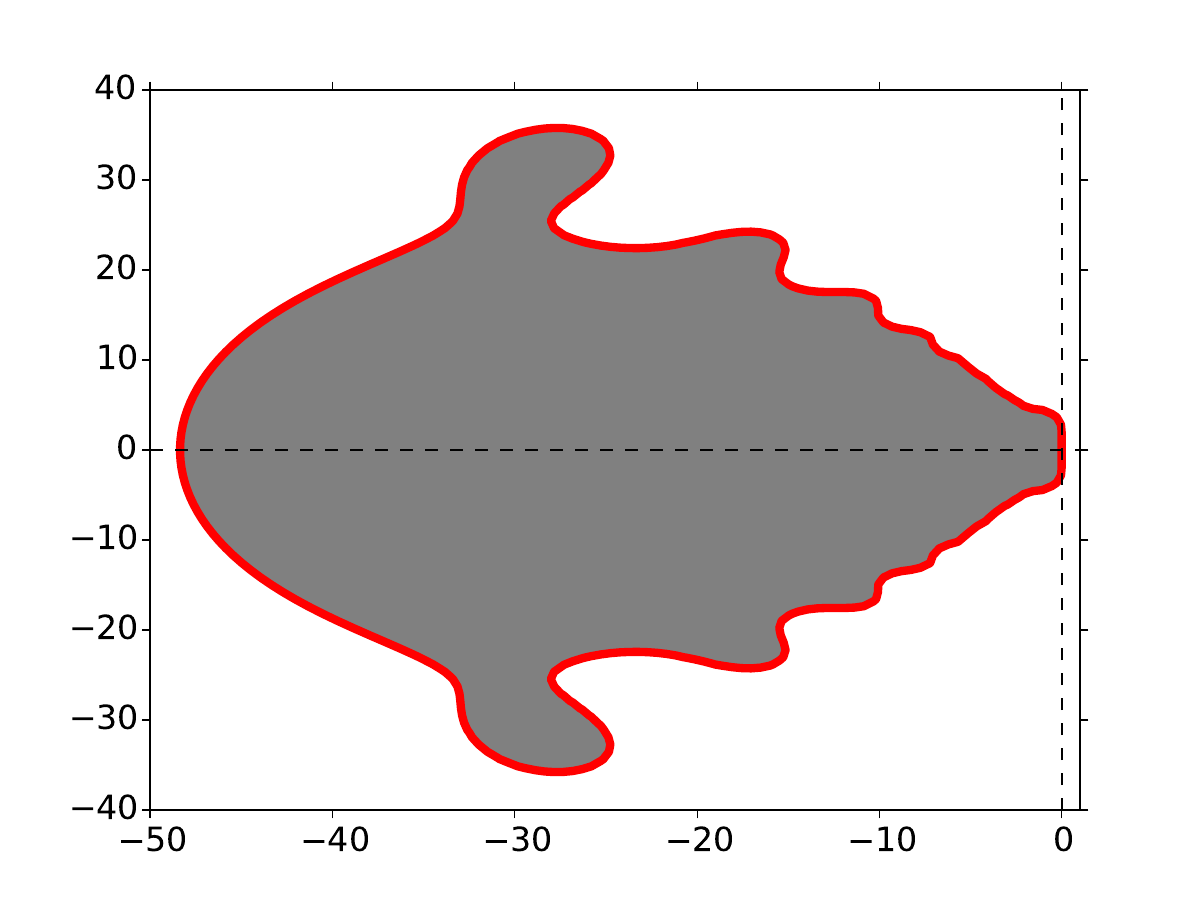}
		\caption{$K=10$, Padé $[2/7]$}
	\end{subfigure}
	\begin{subfigure}{39mm}
		\includegraphics[width=\textwidth]{regions_region_bpl_10_1}
		\caption{$K=10$, Padé $[1/8]$}
	\end{subfigure}
	
	\vspace{\baselineskip}
	\begin{subfigure}{39mm}
		\includegraphics[width=\textwidth]{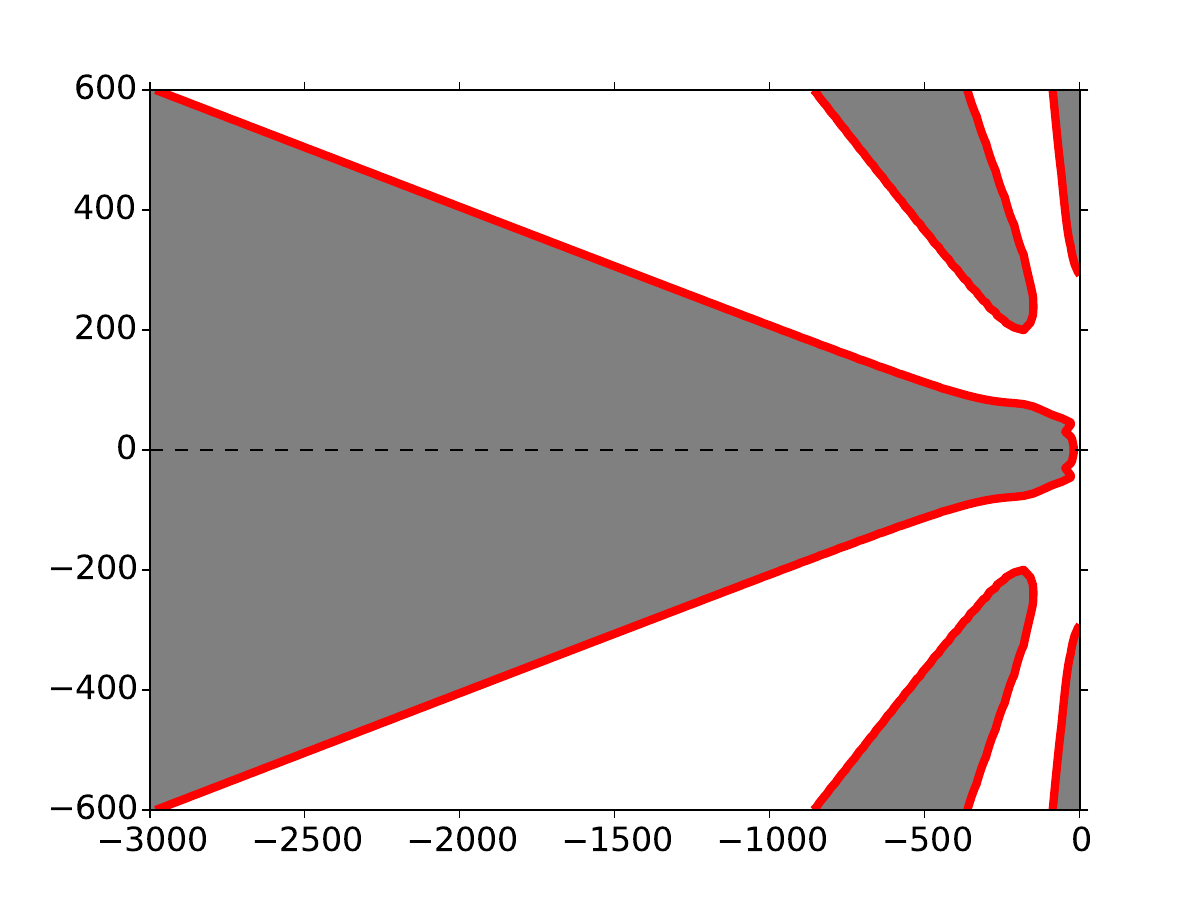}
		\caption{$K=10$, Padé $[0/9]$}\label{region_BPL_c0}
	\end{subfigure}
	\caption{Linear stability regions of Borel-Padé-Laplace integrator for fixed $K=10$ and decreasing $K_a$}
	\label{fig:region_BPL_c}
\end{figure}

To end up, we would like to analyse graphically the influence of $K_a$ (or $K_b$) when $K$ is fixed. The stability regions corresponding to $K=10$ and different Padé degrees are plotted in Figure \ref{fig:region_BPL_c}. It can be observed that the stability region grows with the degree $K_b$ of the Padé denominator. When $K_a=0$ (Figure \ref{region_BPL_c0}), BPL tends to be $A(\alpha)$-stable for some angle $\alpha$. A theoretical study on the optimal choice of $K_a$ and $K_b$, which take into account the stability and the precision, would be very interesting but has not been carried out yet.
In the sequel, a (almost-) diagonal Padé approximant satisfying relation (\ref{ka}) will be chosen. This is motivated by some good properties of diagonal Padé (convergence \cite{baker61,nuttall_1970,baker00}, invariance under linear fractional transformation \cite{baker75}, $A$-stability of diagonal Padé approximants to the exponential function \cite{book:hairer2}, \dots). As seen, it may not correspond to an optimal choice but it will be shown that it is good enough to obtain a very competitive performance in terms of computation time.



Note that some of the plotted stability regions are not complete. Indeed, the whole stability regions may contain other parts in the complex plane, but these parts have been excluded from the graphics.

\vspace{\baselineskip}

In the next section, the performance of BPL in solving stiff and non-stiff equations is analysed.  
As mentioned, all the previous figures were plotted with $N_G=100$ Gauss points. For $N_G=20$ and for $N_G=200$, only small changes have been recorded for Figure \ref{fig:region_BPL}. So, for the upcoming numerical tests, $N_G$ is set to 20.

\section{Numerical performance}

Unless otherwise stated, the order $K$ of BPL is set to 10. The degrees of the numerator and the denominator of the Padé approximant (\ref{pade}) are $K_a=4$ and $K_b=5$.

We compare BPL with some classical numerical schemes. Some of them are popular choices for solving stiff equations. These schemes have either a forth or a tenth consistency order.

\subsection{Classical schemes}

The following schemes are considered. 

\begin{itemize}
	\item The 4-stage 4-th order explicit Runge-Kutta algorithm with a Fehlberg adaptive time step \cite{fehlberg70,book:hairer2}, called RK4 hereinafter. 
	\item The 5-stage 10-th order implicit adaptive Gauss-Legendre method which is a Runge-Kutta scheme combined with a Gauss-Legendre quadrature \cite{book:hairer2}, refered as GAU. 
	\item The 4-step 4-th order implicit backward differentiation formula \cite{curtiss52,book:hairer2}, initialized with RK4, and named BDF in this article.		
	\item The 4-th order exponential time differencing method combined with the adaptive Runge-Kutta-Fehlberg method \cite{cox_2002}. This method is generally called ETDRK4, but will simply be shortened to ETD.
\end{itemize}
RK4 has been chosen for its popularity and speed in solving non-stiff equations, GAU for its order 10 (the same order as that set for BPL), BDF for its popularity in solving stiff equations and ETD because it is a relatively recent integrator for stiff equations. All of these methods are adaptive. The step size is updated at each time iteration with a formula
\[
h_{n+1}=0.9\,h_n\left( \cfrac{\tau}{e_{n+1}}\right)^{\frac1{k+1}}
\]
where $k$ is the order of the scheme and $\tau$ is a small tolerance parameter which can be choosen to adjust the accuracy of the method. $e_{n+1}$ is an estimation of the local error. Indications on how it is computed are given below for each scheme.

For RK4 and GAU, which are multi-stage one-step methods, the approximate solution at $t=t_{n+1}$ can be written as follows:
\[ u_{n+1}=u_n+h_n(b_1k_1+\dots+b_sk_s) \]
where $s$ is the number of stages. The intermediate values $k_i$ and the coefficients $b_i$ are defined in equation~(1.8) and in Table~5.1 of \cite{book:hairer} for RK4, and in equation~(7.7) of \cite{book:hairer} and in page 71 of \cite{book:hairer2} for GAU. The estimated error $e_{n+1}$ is obtained from the difference between $u_{n+1}$ and a second estimation
\[ u^*_{n+1}=u_n+h_n(b^*_1k_1+\dots+b^*_sk_s) \]
of the approximate solution, having at least an order $k$ of constistency. The coefficients $b_i^*$ are provided in Table 5.1 of \cite{book:hairer} for RK4. They are defined in \cite{book:hairer2}, equation (8.16), and in \cite{swart97}, section 2, for GAU.

For BDF, the approximate solution is determined from a relation
\begin{equation}
	a_1^{n+1}u_{n+1}+\dots+a^{n+1}_{s}u_{n+1-s}=h_nf(t_{n+1},u_{n+1})
\end{equation}
where $s=4$ is the step number. The coefficients $a^{n+1}_i$, in the adaptive case, are described in appendix~G of \cite{tomas13}, subsection~G4.
The estimation of the local error is explained in pages 372-373, in  Theorem~6.2 and in Table~6.2 of \cite{book:hairer}.

Lastly, if equation (\ref{equation}) is a scalar ODE then the ETD approximate solution is defined in \cite{cox_2002}, equation (29). When equation (\ref{equation}) is not scalar, a pseudo-inversion and an exponentiation of a matrix is needed. They are carried out respectively with a singular value decomposition and a matrix Padé approximation. The local error is obtained from a perturbation of equation (29) of  \cite{cox_2002}.

All the schemes are implemented entirely in python with a fairly equal effort in optimization. The computations are run on a single processor. 
We focus on accuracy, the size of time step and computation time.

We now apply these schemes to some classes of differential problems.

\subsection{Lotka-Volterra equations}\label{lotka}

Consider a prey-predator system, dynamically governed by the Lotka-Volterra equations \cite{hofbauer88}:
\begin{equation}
\begin{cases}
	\td ut =  &\quad \alpha ~ u  - \beta ~ uv,\\[7pt]
	\td vt = &-\ \delta ~ v  + \gamma ~ uv,\\[5pt]
\end{cases}
\label{lotkavolterra}
\end{equation}
where $u$ and $v$ are respectively the number of preys and predators in the population, and $\alpha,\beta,\delta,\gamma$ are real positive constants. The reproduction parameter $\alpha$  is the natural (exponential) growth rate of preys in absence of predators whereas $\delta$ is the natural decline rate of predators in absence of preys. $\beta v$ is the mortality rate of prey depending on the the number $v$ of predators and $\delta u$ is the birth rate of predators depending on the number of prey eaten. It is straight forward to show that system~(\ref{lotkavolterra}) possesses the first integral:
\begin{equation}
	I(u,v)= \beta v + \gamma u - \alpha \ln v  - \delta \ln u .
	\label{first_integral}
\end{equation}

We first choose a set of coefficients for which the problem is not stiff.

\subsubsection{Non-stiff case}

Take an initial population which consists of two preys and one predator, that is $u_0=2$ and $v_0=1$. The reproduction/decline parameters are set to $\alpha=2/3$ and $\delta=2$ and the predation parameters to $\beta=4/3$ and $\gamma=2$. 

\begin{figure}[htp]
\centering
\begin{subfigure}{\figwidth}
	\includegraphics[scale=1.0,width =\textwidth]{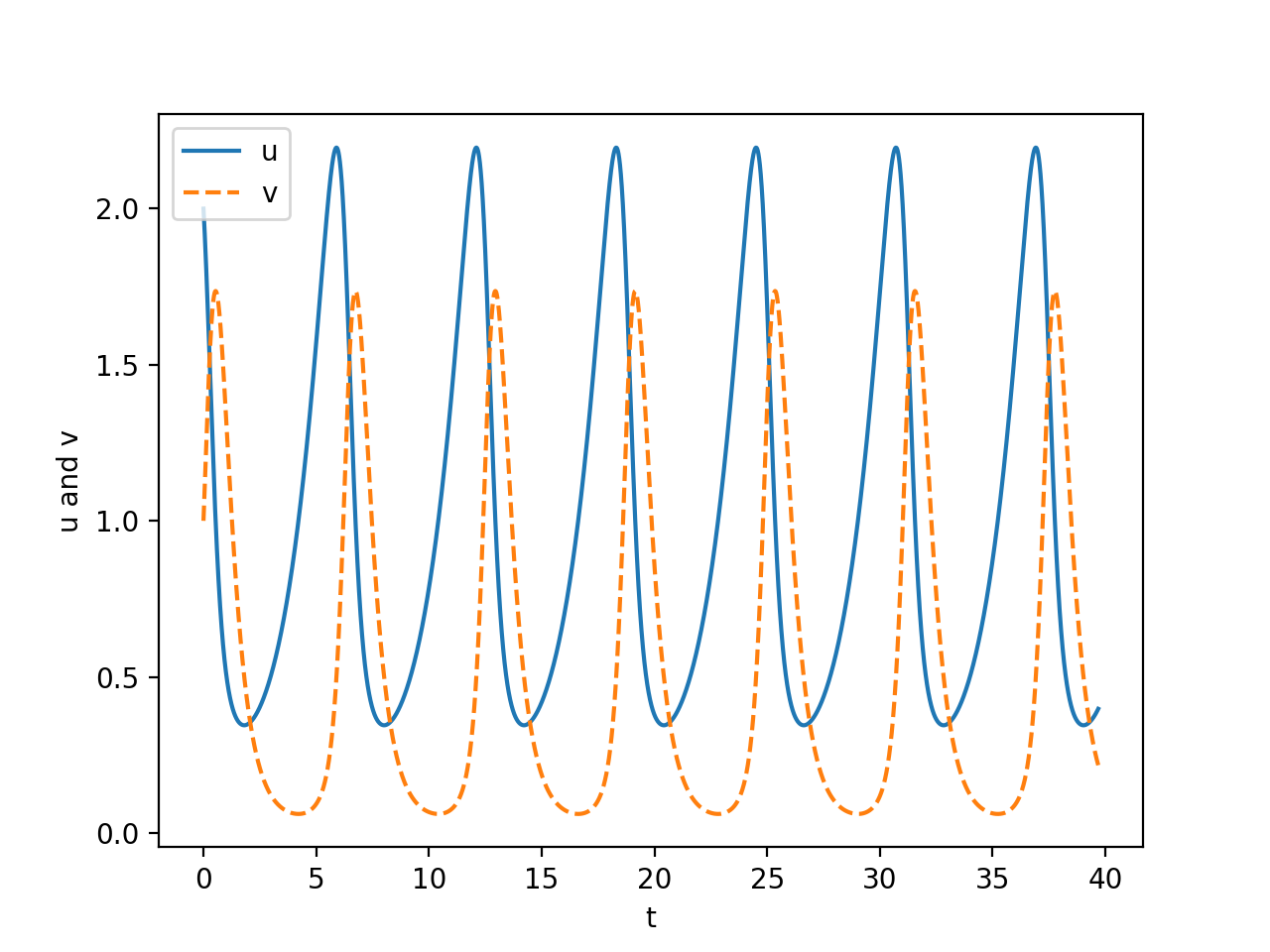}
	\caption{Time evolution}
\end{subfigure}
\begin{subfigure}{\figwidth}
	\includegraphics[scale=1.0,width =\textwidth]{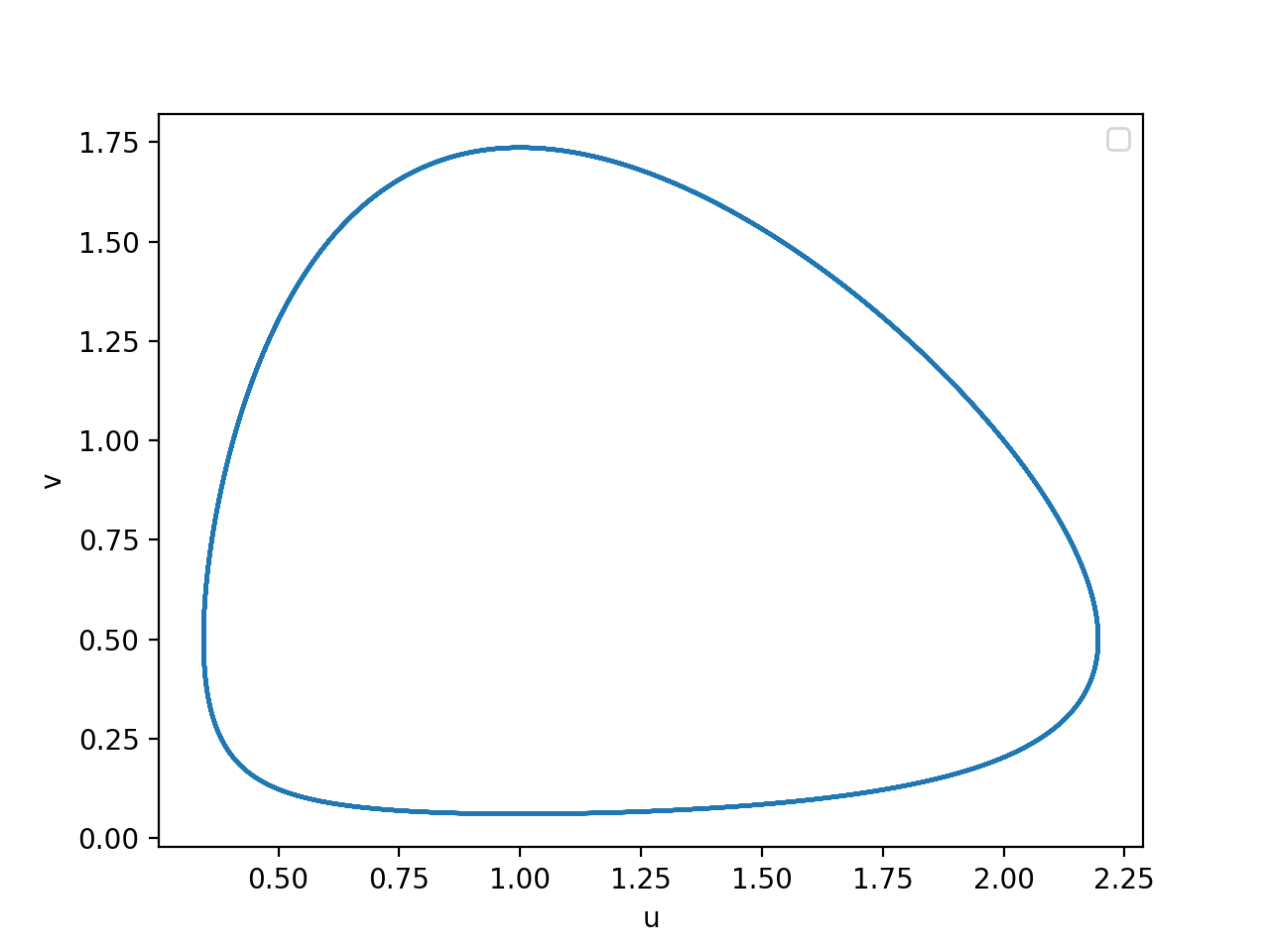}
	\caption{Trajectory in $(u,v)$ plane}
\end{subfigure}
\caption{Approximate solution with BPL}
\label{fig:sol}
\end{figure}

For BPL, the function $F_k$ which operates in the recurrence relation (\ref{ukp1}) is defined by
\begin{equation}F_k(u_0,\dots,u_k,v_0,\dots,v_k)=
	\begin{pmatrix}\displaystyle
		\ \alpha u_k+\beta\sum_{l=0}^ku_lv_{k-l}\\[10pt]\displaystyle
		-\delta v_k+\gamma\sum_{l=0}^ku_lv_{k-l}
	\end{pmatrix}.
\end{equation}
The parameter $ε$ of BPL which is used in the accuracy criterion (\ref{eps}) is set such that the mean error on the first integral (\ref{first_integral}) is about $1.35·10^{-7}$ over a simulation time $T=1000$. This mean or overall error is defined as an approximation of 
\begin{equation}
	\cfrac 1T\int_0^T\bigg|I\big(u(t),v(t)\big)-I\big(u(0),v(0)\big)\bigg|\mathrm{d}t.
\end{equation}
The approximate solution over 40 seconds is visualized in Figure \ref{fig:sol}. It necessitated 254 iterations. To obtain the smooth plots in Figure \ref{fig:sol}, not only the value of the solution at the discrete times $(t_i)_{i=0,\dots,254}$ but also at some intermediate times $t\in]t_i,t_{i+1}[$ are plotted. In contrast to many other schemes, no interpolation method is needed for this. Formula (\ref{gausslaguerre}) directly provides the approximate solution within each interval $]t_i,t_{i+1}[$.

The accuracy parameter $\tau$ are set for each step (BDF, ETD, GAU and RK4) such that their a posteriori accuracies are comparable to that of BPL. The mean errors are reported in Table \ref{tab:error}. They are around $4·10^{-7}$.
\begin{table}
	\centering\small
	\begin{tabular}{V{2.5}lcccccV{2.5}}
		\hlineB{2.5}
				&\textbf{BDF}	&\textbf{BPL}	&\textbf{ETD}	&\textbf{GAU}	&\textbf{RK4}\\\hline
		Mean error	&$5.56·10^{-7}$	&$1.35·10^{-7}$	&$7.22·10^{-7}$	&$4.6·10^{-7}$	&$2.38·10^{-7}$
		\\
		Mean time step	&$2.42·10^{-3}$	&$1.65·10^{-1}$	&$2.07·10^{-4}$	&$3.70·10^{-2}$	&$3.10·10^{-2}$
		\\
		CPU		&$5.96·10^2$	&$4.32$		&$9.34·10^2$	&$4.50·10^{1}$	&$4.96$
		\\\hlineB{2.5}
	\end{tabular}
	\caption{Error on the first integral and CPU time}
	\label{tab:error}
\end{table}

\begin{figure}
	\centering
\begin{subfigure}{\figwidth}
	\includegraphics[width=\figwidth]{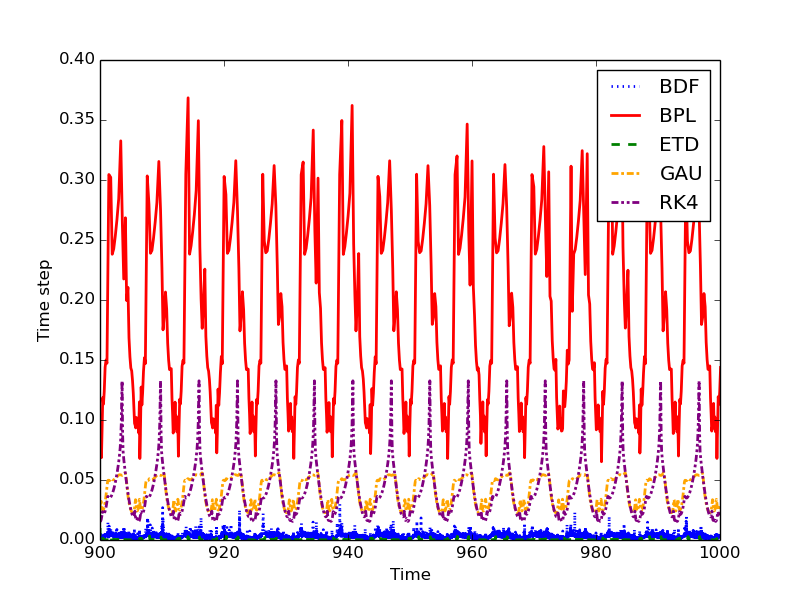}
	\caption{Linear scales}\label{timestep_comparable_precision:lin}
\end{subfigure}
\begin{subfigure}{\figwidth}
	\includegraphics[width=\figwidth]{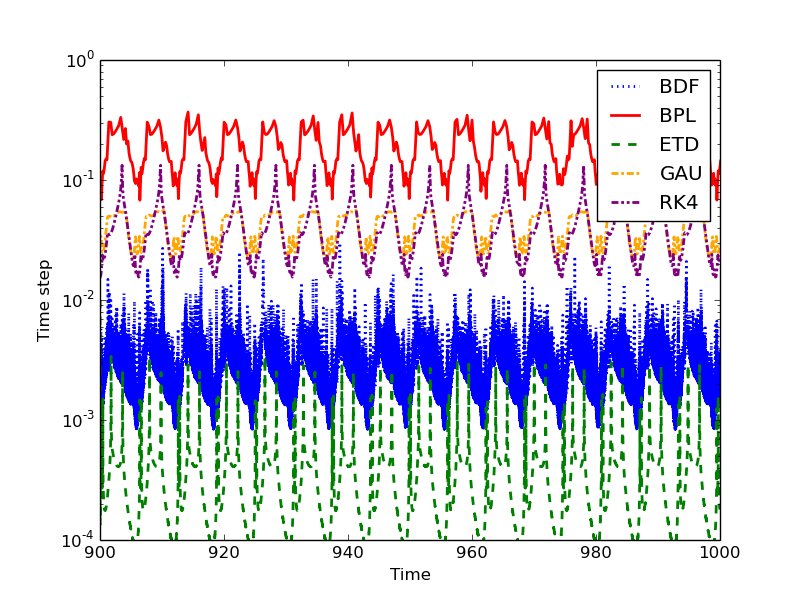}
	\caption{Semi-logarithmic scale}\label{timestep_comparable_precision:log}
\end{subfigure}
	\caption{Non-stiff Lotka-Volterra. Evolution of time step}
	\label{timestep_comparable_precision}
\end{figure}

Figure \ref{timestep_comparable_precision:lin} shows the evolution of the time steps of the different methods. The solution being periodic, only the evolution over the last 100 seconds are plotted. As can be seen, it is with BPL that the time step is the largest. Since we are in a non-stiff case, the classical Runge-Kutta method has a good performance and competes with the 10-th order Gauss scheme in terms of time step. Figure \ref{timestep_comparable_precision:log} represents the same data as Figure  \ref{timestep_comparable_precision:lin} but with a logarithmic scale in ordinate. It shows that the time step of BPL is about 60 times larger than that of BDF and about 80 times as large as that of ETD. It is confirmed in Table \ref{tab:error} which compares the mean values. As for CPU time, the two explicit integrators, BPL and RK4, have a comparable performance (see Table \ref{tab:error}). They need about 10 times less computation time than GAU and at least 100 times less than BDF and ETD.

\begin{figure}[tp]
	\centering
	\includegraphics[width=\figwidth]{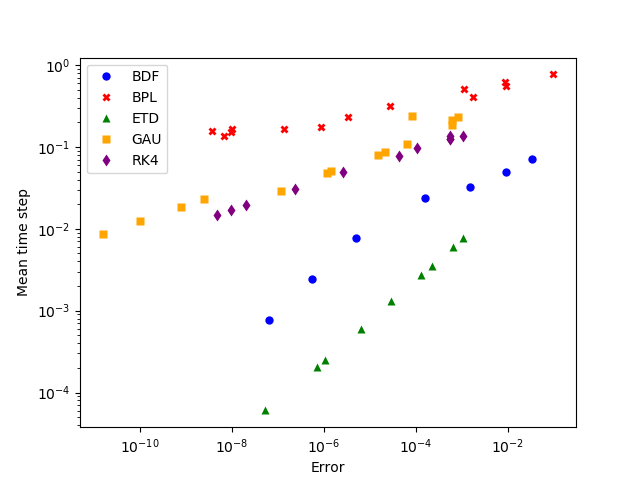}
	\caption{Non-stiff Lotka-Volterra. Evolution of the mean time step with the mean error
	\label{precision_timestep}}
	\centering
	\includegraphics[width=\figwidth]{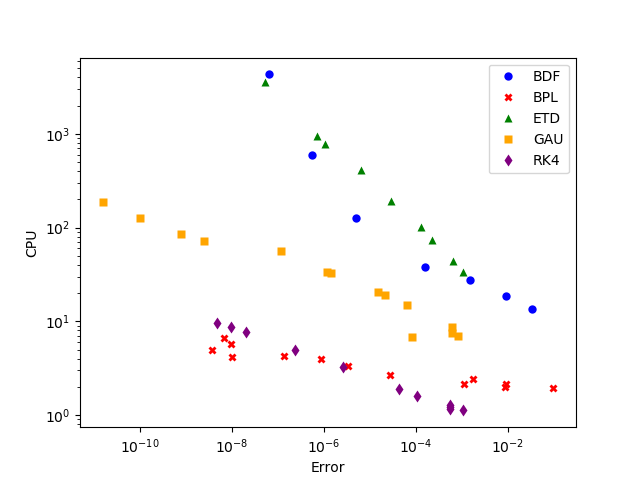}
	\caption{Non-stiff Lotka-Volterra. Evolution of CPU with the mean error
	\label{precision_cpu}}
\end{figure}

In a second test, each scheme is run with multiple values of the (residue or estimated error) tolerance. The mean time step is plotted in Figure \ref{precision_timestep} against the overall accuracy. This figure clearly shows that, amongst the considered schemes, BPL has always the largest mean time step, whatever the precision. This mean time step is about 6.6 times as large as that of the Gauss scheme with the same order, for an error around $3·10^{-9}$. This large time step results in a faster computation. Indeed, as can be noticed in Figure \ref{precision_cpu}, BPL requires much less CPU time than GAU, for comparable precisions. Only RK4 is faster than BPL for a medium or a low precision. But when a high precision is required, BPL tends to be more interesting.

\subsubsection{Increasing the stiffness ratio}

We now examine the behaviour of the schemes when the stiffness ratio varies. The stiffness ratio $r$ is defined as the spectral condition number of the linear part of equations (\ref{lotkavolterra}), that is
\begin{equation}
	r=\cfrac{\max(\alpha,\delta)}{\min(\alpha,\delta)}
	\label{stiffness}
\end{equation}
since $\alpha$ and $\delta$ are positive real numbers. For the sake of simplicity, and since it will be the case in the numerical experiments, assume that $\delta> \alpha$, such that $r=\delta/\alpha$.
In fact, $r$ appears naturally when equations (\ref{lotkavolterra}) are adimensionalized  with the variables
\begin{equation}
	u^*=\frac{\gamma u}{\delta},\quad v^*=\frac{\beta v}{\alpha},\quad t^*=\alpha t. 
\end{equation}
Indeed, equations (\ref{lotkavolterra}) can be written as follows:
\begin{equation}
\begin{cases}
	\td {u^*}{t^*} =  &\quad   u^*(1  -   v^*),\\[7pt]
	\td {v^*}{t^*} = &rv^*(-   1  +    u^*).\\[5pt]
\end{cases}
\end{equation}

All the parameters of the equations are kept at the same value as before, except $\delta$ which is increased. 
As before, simulations are run with multiple values of the (residue or estimated error) tolerance until 1000 seconds. The error on the first integral and the CPU time are recorded and plotted hereafter.

\begin{figure}[tp]
	\centering
	\includegraphics[width=.6\textwidth]{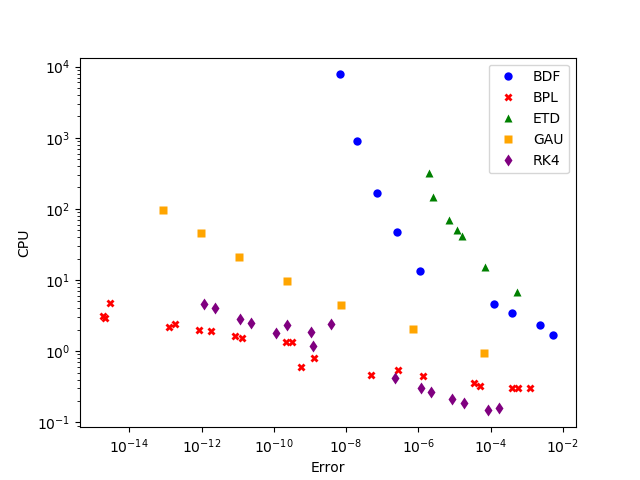}
	\caption{Lotka-Volterra. Stiffness ratio $r=8$\label{r8}}
%
%
	\includegraphics[width=.6\textwidth]{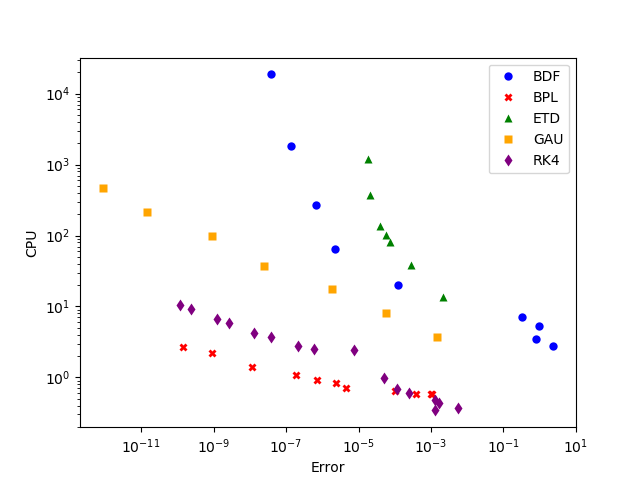}
	\caption{Lotka-Volterra. Stiffness ratio $r=16$\label{r16}}
	\includegraphics[width=.6\textwidth]{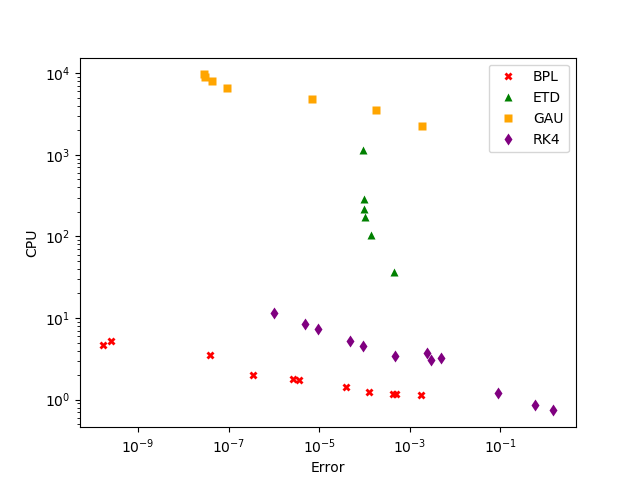}
	\caption{Lotka-Volterra. Stiffness ratio $r=32$\label{r32}}
\end{figure}

For a moderate stiffness ratio $r=8$, Figure \ref{r8} shows that RK4 and BPL compete in terms of CPU time, even if BPL indicates a slight advantage for high precision simulations. It can also be stated in this figure that GAU can provide a very accurate solution, due to its high order, but with a higher cost than BPL. BDF and ETD are much more expensive than the other schemes.

For $r=16$, we have approximately the same picture, except that BPL becomes more interesting than RK4 even for moderate precisions. This can be observed in Figure \ref{r16}.

When the stiffness ratio is set to a high value $r=32$, the situation changes significantly. First, as can be seen in Figure \ref{r32}, RK4 cannot reach very high precision any longer, compared to BPL. The precision that can be achieved with GAU is still very high but not as high as with $r=16$. ETD also looses precision. Only BPL is able to maintain the same precision as previously.

Concerning the numerical cost, the increase of CPU time needed by BPL and ETD is very small compared to that of GAU.

Lastly, BDF does not appear in Figure \ref{r32}. Indeed, although it is a popular method for stiff equations, it fails with $r=32$. It diverges as soon as $t$ reaches few seconds. This behaviour has also been observed with the optimized BDF solver of the python \texttt{scipy} package, with the optimized BDF solver of Scilab, and with the option \texttt{CVODE\_BDF} of the package \texttt{Sundials} of Julia language.

\begin{figure}
	\centering
	\includegraphics[width=\figwidth]{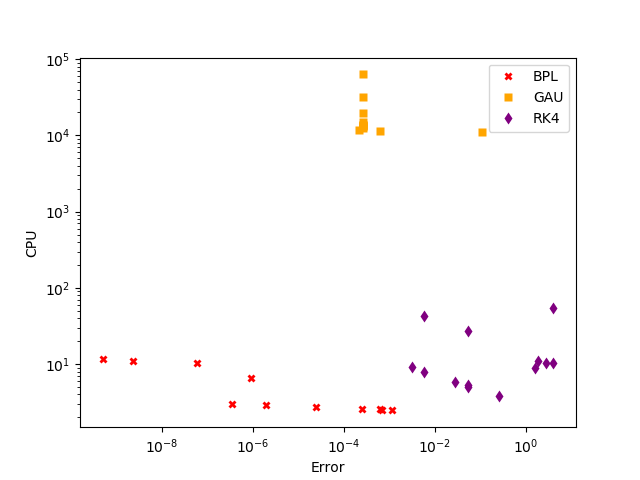}
	\caption{Lotka-Volterra. Stiffness ratio $r=64$}\label{r64}
	\includegraphics[width=\figwidth]{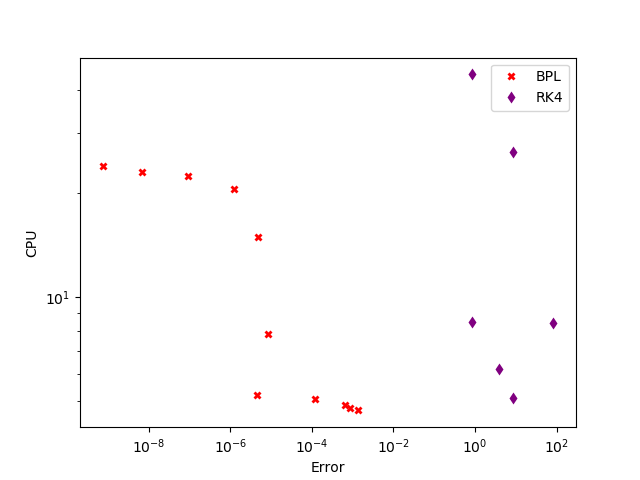}
	\caption{Lotka-Volterra. Stiffness ratio $r=128$}\label{r128}
\end{figure}

For $r=64$, ETD also fails. As remarked in Figure \ref{r64}, RK4 gives moderately accurate solutions, and even wrong solutions for some values of the predicted error tolerance. Indeed, even for very small value of the tolerance, the overall error may be larger than one. Figure \ref{r64} also shows that the precision of GAU seems to stagnate around $2.6·10^{-4}$. Only BPL can provide highly accurate solutions. Moreover, its CPU cost is very small compared to that of GAU.

At last, with $r=128$, the Gauss method also fails. This behaviour has as well been observed with the Gauss solver of the Matlab package \texttt{numeric::odesolve}. For this value of the stiffness ratio, RK4 cannot give an accurate solution any longer, whereas with BPL, the error can be as small as $7.4·10^{-10}$ (see Figure~\ref{r128}).

These numerical experiments shows that BPL is an interesting alternative method for stiff problems. First, its arbitrary high order allows to get highly accurate solutions. With Lotka-Volterra equations, it never fails for values of $r$ up to 128. Moreover, its cost is generally much smaller than that of the other methods, due to its explicit property. 

The previous tests show the performance of BPL for the resolution of stiff and non-stiff ODE's. In the next subsection, we examine its efficiency in solving partial differential equations.

\subsection{Korteweg-de-Vries equation}

In this subsection, we consider the Korteweg-de-Vries equation (KdV)
\begin{equation}
\label{kdv}
\cfrac{\partial u}{\partial t}  + c_0 \cfrac{\partial u}{\partial x} + \beta \cfrac{\partial^3 u}{\partial x^3} + \cfrac{\alpha}2 \cfrac{\partial u^2}{\partial x}= 0
\end{equation}
which models waves on shallow water surfaces \cite{korteweg95}. In this equation, the linear propagation velocity $c_0$, the non-linear coefficient $\alpha$ and the dispersion coefficient $\beta$ are positive constants, linked to the gravity acceleration $g$ and the mean depth $d$ of the water by:
\begin{equation}
	c_0=\sqrt{gd},\quad \alpha=\cfrac32\,\sqrt{\cfrac gd},\quad \beta=\cfrac{d^2c_0}6.
\end{equation}

In order to focus on the performance of the time integrators, we choose a high order scheme, namely a spectral method, for the space discretization. The solution is assumed to be periodic with period $X$ in space, and integrable. It is approximated by its truncated Fourier series:
\begin{equation}
	u(x,t)\simeq \sum_{|m|\leq M}\hat u^m(t)\e^{im\omega x},
	\label{fourier}
\end{equation}
where $M\in \mathbb N$ and $\omega=\frac{2\pi}X$.
The substitution of equation (\ref{fourier}) into (\ref{kdv}) leads to a $(2M+1)$-dimensional ODE
\begin{equation}
	\cfrac{\textrm d\hat u}{\textrm d t}=A\hat u+N(\hat u)
	\label{kdv_ode}
\end{equation}
where the array $\hat u$ contains the unknowns $\hat u^m$, $A$ is a diagonal matrix with diagonal entries 
\begin{equation}
	A^m_m=-c_0i\omega m+i\beta\omega^3 m^3
\end{equation}
and $N(\hat u)$ is a non-linear array containing convolution terms:
\begin{equation}
	N(\hat u)=-\cfrac12\,i\alpha m\omega \ \hat u*\hat u.
\end{equation}
Convolution operations are performed in physical space and the standard dealiasing 3/2 rule is applied.

With BPL, each component $\hat u^m(t)$ of the Fourier coefficient array $\hat u(t)$ is decomposed into its Taylor series
\begin{equation}
	\hat u^m(t)=\sum_{k=0}^K\hat u^m_kt^k.
	\label{bpl_fourier}
\end{equation}
The series coefficients are computed explicitely as follows:
\begin{equation}
	\hat u_{k+1}=\cfrac1{k+1}\left[(-c_0i\omega m+i\beta\omega^3 m^3)\hat u_k-\cfrac12\,i\alpha m\omega \ \sum_{l=0}^k\hat u_n*\hat u_{k-l}\right].
	\label{uk1}
\end{equation}
For the simulations, the initial condition is the periodic prolongation of the function
\begin{equation}
	u_0(x)=U\operatorname{sech}^2(\kappa x), \quad\quad\quad x\in\left[ -\frac X2,\frac X2 \right],
	\label{kdv_initial}
\end{equation}
$U$ being a constant and
$\kappa =\sqrt{\frac{3U}{4d^3}}.$
The corresponding exact solution is the traveling wave
\begin{equation}
	u(x,t)=u_0(x-ct).
	\label{kdv_solution}
\end{equation}
with
$c=c_0\left(1+\frac U{2d}\right).$
We take $X=24\pi$, $d=2$, $g=10$ and $U=\frac12$. The solution is periodic in time, with a period $T\simeq14.986$s. 

We use $D=2M$ to indicate the size of the system, instead of the dimension $2M+1$ of equation (\ref{kdv_ode}). The simulations are run over one period, for some values of $D$ between 64 and 512.

It is hard to calibrate the tolerance parameter $\tau$ of all the schemes to have the same (a posteriori) overall error at each value of $D$. So, this calibration has not been done. Instead, we require more accuracy to BPL than to the other methods, in order not to overestimate the performance of BPL. The overall error are recorded in Table \ref{tab:kdv_error}. The error reported in this table is an approximation of 
\begin{equation}
	\int_0^T\cfrac{\|u_{computed}(t)-u_{exact}(t)\|}{\|u_{exact}(t)\|}ｄt.
	\label{eq:kdv_error}
\end{equation}

\begin{table}
	\centering
	\begin{tabular}{V{2.5}cccccccV{2.5}}
		\hlineB{2.5}
			$D$	&&\textbf{BDF}	&\textbf{BPL}	&\textbf{ETD}	&\textbf{GAU}	&\textbf{RK4}\\\hline
			\textbf{64}	&&$1.09·10^{-1}$	&$3.71·10^{-4}$	&$2.92·10^{-3}$ &$5.11·10^{-3}$ &$1.83·10^{-3}$\\
			 \textbf{128}	&&$1.08·10^{-2}$	&$3.54·10^{-4}$ &$3.27·10^{-3}$ &$5.81·10^{-3}$ &$1.69·10^{-3}$\\
			 \textbf{256}	&&--		&$3.61·10^{-4}$ &$3.66·10^{-3}$ &$4.00·10^{-3}$ &$1.23·10^{-3}$\\
			 \textbf{512}	&&--		&$3.17·10^{-4}$ &$4.11·10^{-3}$ &$2.65·10^{-3}$ &$6.50·10^{-4}$
		\\\hlineB{2.5}
	\end{tabular}
	\caption{KdV. Overall error}
	\label{tab:kdv_error}
\end{table}

\begin{figure}
	\centering
	\includegraphics[width=\figwidth]{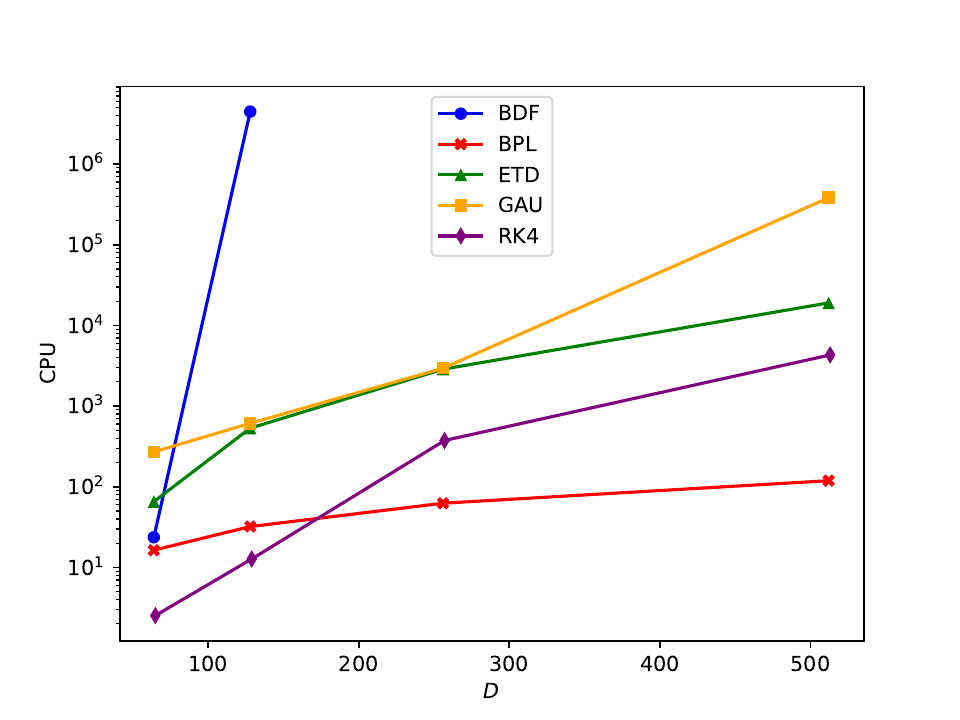}
	\caption{KdV. Evolution of the computation time with the size $D$ of the problem}
	\label{fig:kdv_cpu}
\end{figure}

The evolution of the computation time of each scheme is plotted in Figure \ref{fig:kdv_cpu}. It can be seen there that BDF requires a very high cost for $D=128$, despite the low precision (see second column of Table \ref{tab:kdv_error}). As a consequence, it has not been used for higher values of $D$.

Figure \ref{fig:kdv_cpu} also shows that, among the considered schemes, RK4 is the fastest for (non-stiff) small-sized problems, for the given precisions. But for high degrees of freedom, BPL becomes the most interesting in terms of computational time. BPL also has the smallest slope.

\begin{figure}
	\centering
	\includegraphics[width=\figwidth]{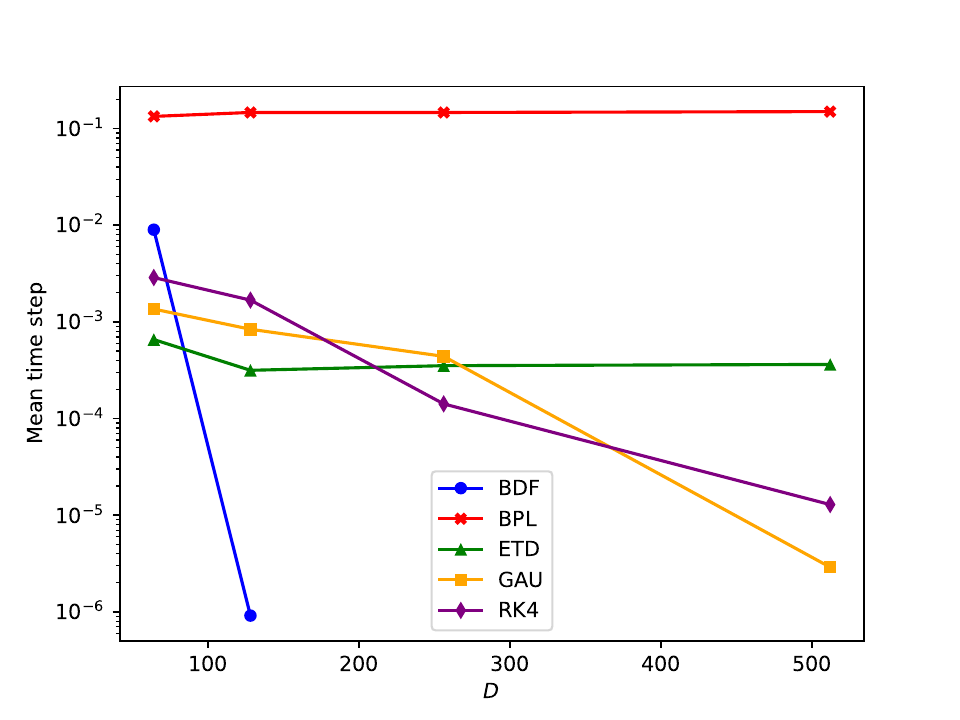}
	\caption{KdV. Evolution of the mean time step with the size $D$ of the problem}
	\label{fig:kdv_timestep}
\end{figure}

Figure \ref{fig:kdv_timestep} indicates that BPL has a very large mean time step compared to the other schemes, whatever the size of the problem is. It is also striking that the mean time step does not vary very much with the size of the problem. However, BPL is not the only scheme which presents this characteristics since ETD exhibits the same behavior, but with much smaller time steps.

\begin{figure}
	\centering
	\includegraphics[width=\figwidth]{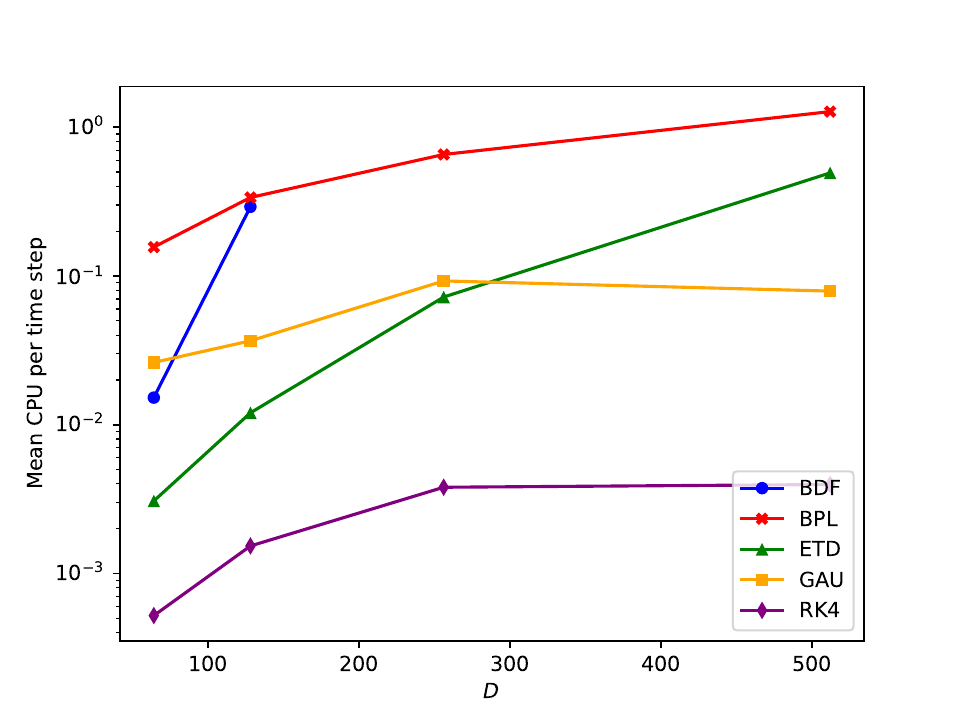}
	\caption{KdV. Mean CPU time per time step}
	\label{fig:kdv_cpu_timestep}
\end{figure}

The large time step of BPL is of a great importance in its performance. Indeed, the CPU time spent at each time step is very high with BPL in comparison to the other schemes, as can be stated in Figure \ref{fig:kdv_cpu_timestep}. One reason for this is the evaluation of the residue in step 6 of the algorithm presented in section \ref{sec:algorithm}. This evaluation is done multiple times at each time step to decide if the solution is still accurate enough. Another precision evaluation is desirable, but not available yet. Fortunately, this expensive precision evaluation is largely counter-balanced by large time steps.

\begin{figure}
	\centering
	\includegraphics[width=\figwidth]{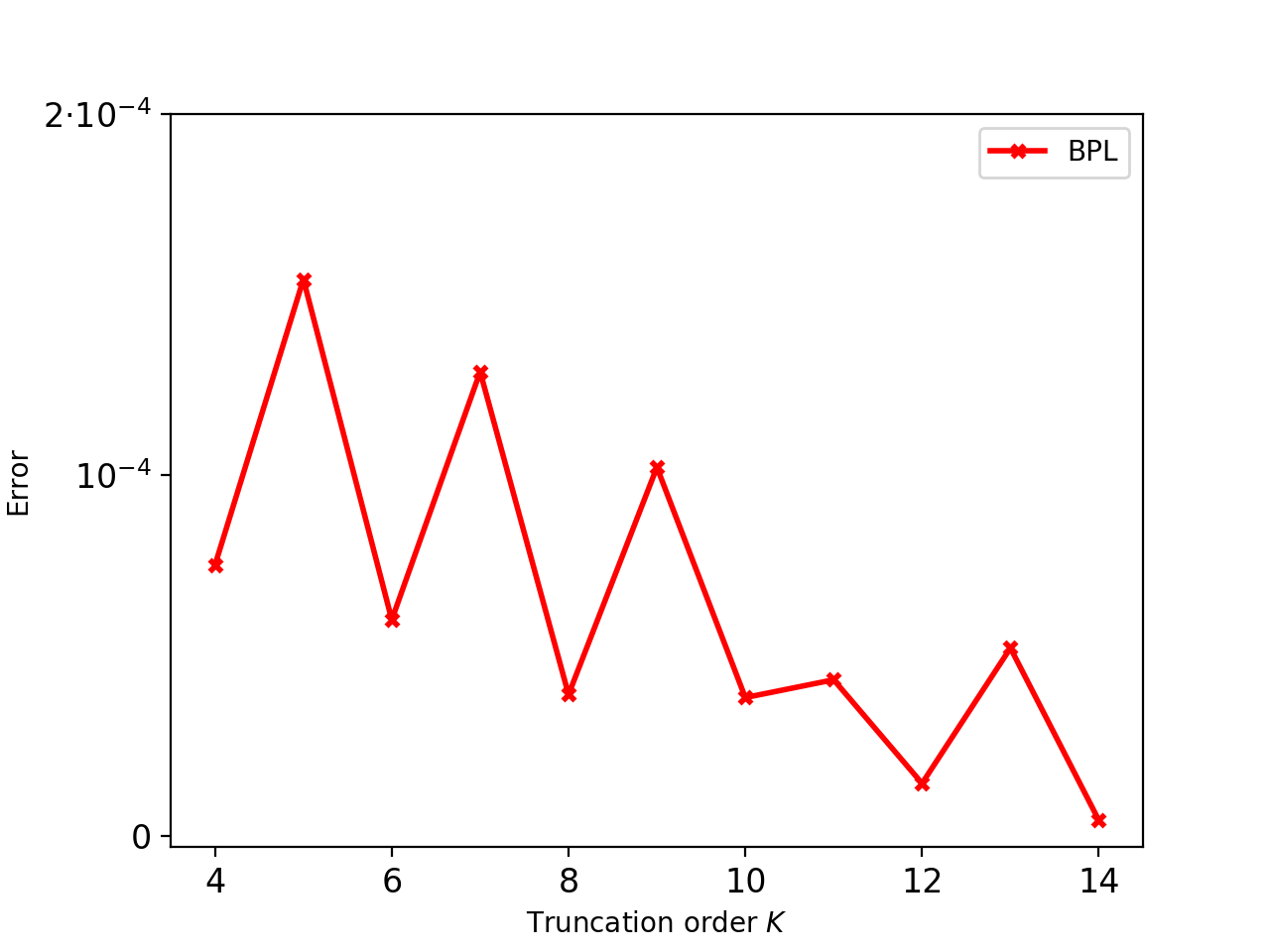}
	\caption{KdV. Evolution of the error with $K$\label{fig:kdv_order_error}}
\end{figure}

In all of the previous simulations, the order $K$ of the time series in BPL was set to 10. In our last test, the effect of $K$ on the performance of BPL is analysed. For this, the size of the problem is set to $D=128$. A residue tolerance $ε=1·10^{-4}$ is chosen. Figure \ref{fig:kdv_order_error} shows the $L^1$ relative error (defined in equation (\ref{eq:kdv_error})) over one period. This figure reveals a fluctuation of the error according to the parity of $K$. Note that such fluctuation is not uncommon when manipulating truncated series. Moreover, the parity of $K$ intervenes in the choice of the Padé approximants in Borel space. Indeed, when $K$ is odd, the numerator and the denominator of the Padé approximant have the same degree; and when $K$ is even, the denominator has a higher degree than the numerator (see choice in equation (\ref{ka})). Figure  \ref{fig:kdv_order_error} however tells us that globally, the accuracy increases with the order $K$ of the series, for a fixed value of the residue tolerance.

\begin{figure}
	\centering
	\includegraphics[width=\figwidth]{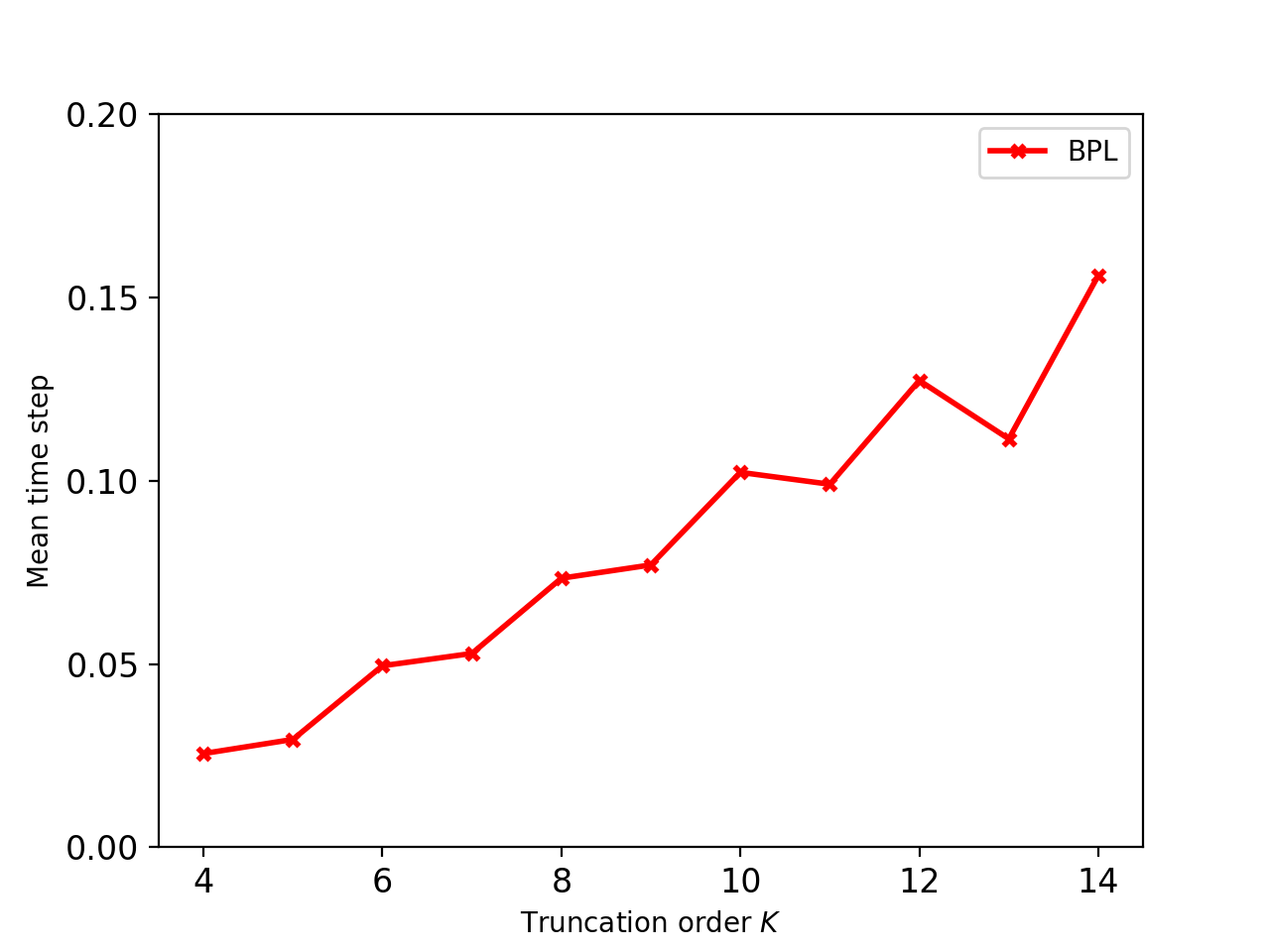}
	\caption{KdV. Evolution of the mean time step with $K$\label{fig:kdv_order_tim}}
\end{figure}

\begin{figure}
	\centering
	\includegraphics[width=\figwidth]{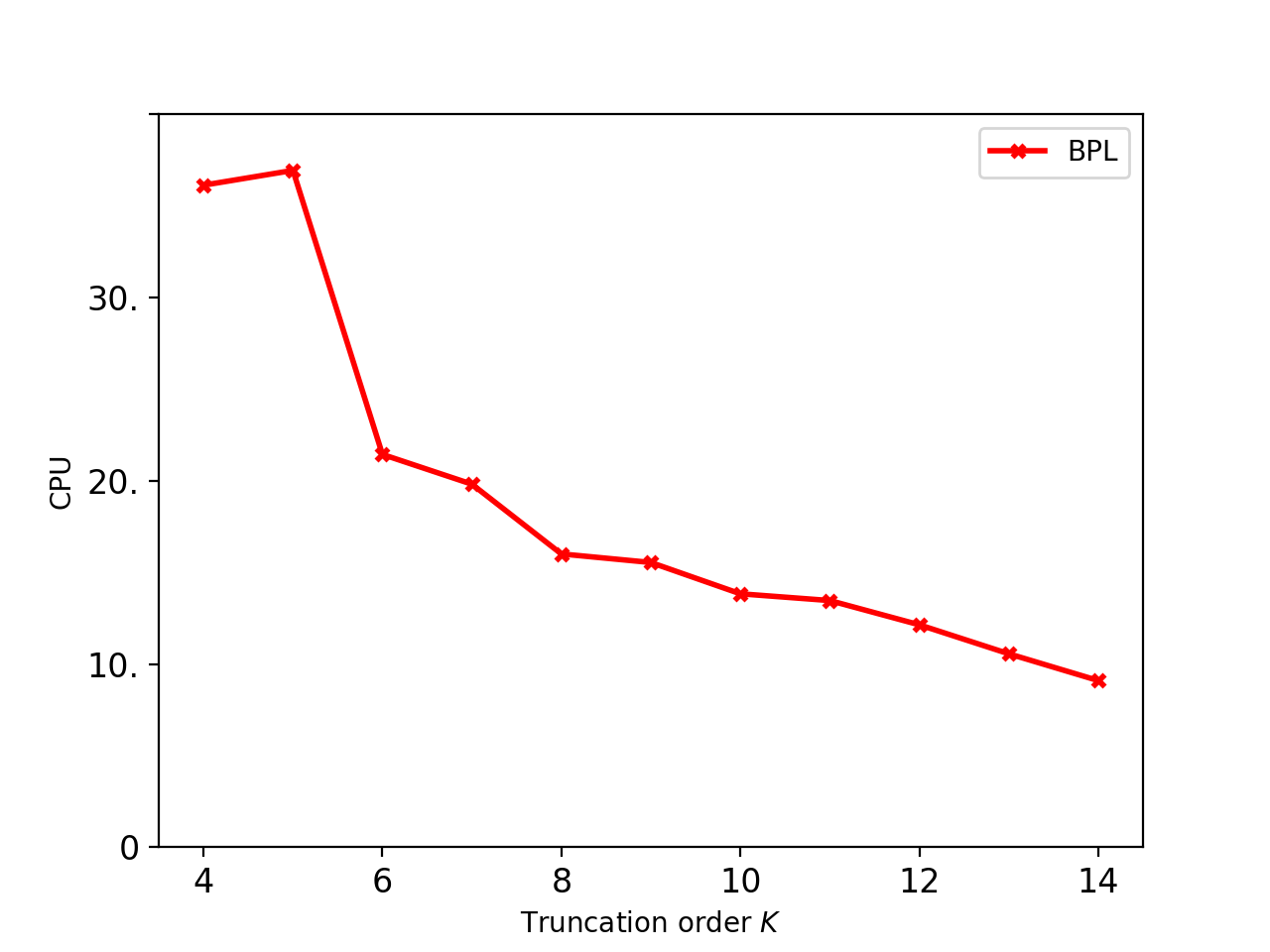}
	\caption{KdV. Evolution of the computation time with $K$\label{fig:kdv_order_cpu}}
\end{figure}

The mean time step has also a globally increasing tendency with $K$ as can be seen in Figure \ref{fig:kdv_order_tim}, passing from $Δt_{mean}=0.0256$ for $K=4$ to $Δt_{mean}=0.156$ when $K=14$. As a consequence, the CPU time decreases with $K$, as can be noted in Figure \ref{fig:kdv_order_cpu}. These results tend to indicate that high values of $K$ accelerate the computation. 


\section{Conclusion}

In this article, we studied the linear stability of the Borel-Padé-Laplace integrator. It has been shown that if the summation procedure is not applied, the scheme has the same linear stability domain as an explicit Runge-Kutta integrator. But when the summation is carried out, the linear stability domain enlarges very significantly, even when the Taylor series of the solution is convergent. It has been obsereved that the size of this domain increases with the truncation order of the series. We also saw that the choice of Padé approximants in Borel space has a substantial impact on the size of the linear stability domain.

Even if BPL is not $A$-stable, it has been shown that this scheme is more efficient than many explicit and implicit ones, in solving stiff problems. It runs without any particular difficulty for a wide range of values of the stiffness ratio. Due to its high order, it can reach very high precisions even when the stiffness number is high. Moreover, its explicit property makes it very fast compared to the other integrators.

Numerical tests on non-stiff Lotka-Volterra and on Korteweg-de-Vries equations showed that for small-size systems, the popular 4-th order Runge-Kutta method has a comparable speed than BPL when only a moderate precision is needed. But when high precision is required BPL becomes more interesting. It is even more true when the size of the system is large.

It is worth to notice that increasing the approximation order of BPL does not require any programming effort. One has simply to raise the cut-off parameter $K$ of the series, without changing anything else in the algorithm. As could be observed in the last part of the article, the higher this value is, the faster BPL is.

Despite its speed, one optimization should be brought to the algorithm of BPL. Indeed, a numerical test with Korteweg-de-Vries equation showed that a BPL time step is rather expensive, due among others to many evaluations of the residue. A more efficient accuracy estimation should be developed. This should increase the speed of the scheme.

To obtain the previous results, the computation was done on a single processor. Note however that BPL also presents some advantage regarding parallelization. Indeed, the summation algorithm can be done component-wise, letting the computation to be shared between many processors.

In this paper, only the computational aspects of BPL are discussed. The founding theory was skipped. Yet, some optimizations may be brought to the algorithm with help of theoretical considerations. For instance, a theoretical study of the equation may be helpful to determine the actual Gevrey order (which was set to one in this article). However, it is conceivable to evaluate numerically this Gevrey order from the coefficients of the series. A theoretical study of the equation may also help to find a better (than the real positive semi-line) integration direction in the Laplace transform.



\bibliographystyle{unsrt}

\begin{thebibliography}{10}
\expandafter\ifx\csname url\endcsname\relax
  \def\url#1{\texttt{#1}}\fi
\expandafter\ifx\csname urlprefix\endcsname\relax\def\urlprefix{URL }\fi
\expandafter\ifx\csname href\endcsname\relax
  \def\href#1#2{#2} \def\path#1{#1}\fi

\bibitem{book:hairer2}
E.~Hairer, G.~Wanner, Solving Ordinary Differential Equations II: Stiff and
  Differential-Algebraic Problems, {S}econd {R}evised, {C}orrected second
  printing Edition, Springer Series in Computational Mathematics, Springer,
  2002.

\bibitem{isereles96}
A.~Iserles, A first course in the numerical analysis of differential equations,
  Cambridge University Press, 1996.

\bibitem{butcher03}
J.~Butcher, Numerical Methods for Ordinary Differential Equations, J. Wiley \&
  Sons, Ltd., 2003.

\bibitem{lambert1991}
J.~D. Lambert, Numerical Methods for Ordinary Differential Systems: The Initial
  Value Problem, Wiley, 1991.

\bibitem{willoughby74}
R.~Willoughby (Ed.), Stiff differential systems, Plenum Press, 1974.

\bibitem{hackbusch14}
W.~Hackbusch, The Concept of Stability in Numerical Mathematics, Vol.~45 of
  Springer Series in Computational Mathematics, Springer, 2014.

\bibitem{curtiss52}
C.~F. Curtiss, J.~O. Hirschfelder, Integration of stiff equations, Proceedings
  of the National Academy of Sciences of the United States of America 38~(3)
  (1952) 235--243.

\bibitem{runge95}
C.~Runge, {\"U}ber die numerische {A}uflösung von {D}ifferentialgleichungen,
  Mathematische Annelen 46 (1895) 167--178.

\bibitem{butcher63}
J.~Butcher, Coefficients for the study of {R}unge-{K}utta integration
  processes, Journal of the Australian Mathematical Society 3~(2) (1963)
  185--201.

\bibitem{butcher96}
J.~Butcher, A history of {R}unge-{K}utta methods, Applied Numerical Mathematics
  20~(3) (1996) 247--260.

\bibitem{butcher16}
J.~Butcher, Numerical Methods for Ordinary Differential Equations, 3rd Edition,
  Wiley, 2016.

\bibitem{verwer96}
J.~Verwer, Explicit {R}unge-{K}utta methods for parabolic partial differential
  equations, Applied Numerical Mathematics 22~(1) (1996) 359 -- 379.

\bibitem{Friedli_1978}
A.~Friedli, Verallgemeinerte {R}unge-{K}utta {V}erfahren zur {L}\"osung steifer
  {D}ifferentialgleichungssysteme, in: Bulirsch, Grigorieff, Schröder (Eds.),
  Numerical Treatment of Differential Equations, Oberwolfach 1976, Springer
  Berlin Heidelberg, 1978, pp. 35--50.

\bibitem{Norsett_1969}
P.~Norsett, An a-stable modification of the {A}dams-{B}ashforth methods, in:
  J.~L. Morris (Ed.), Conference on the Numerical Solution of Differential
  Equations, Dundee/Scotland, Springer Berlin Heidelberg, 1969, pp. 214--219.

\bibitem{vanderhoven_1974}
P.~v.~d. Houwen, J.~Verwer, Generalized linear multistep methods, 1 :
  Development of algorithms with zero-parasitic roots, Stichting Mathematisch
  Centrum. Numerieke Wiskunde 74~(10) (1974) 1--16.

\bibitem{certaine1960}
J.~Certaine, The solution of ordinary differential equations with large time
  constants, Mathematical methods for digital computers (1960) 128--132.

\bibitem{hochbruck2005323}
M.~Hochbruck, A.~Ostermann, Exponential {R}unge–{K}utta methods for parabolic
  problems, Applied Numerical Mathematics 53~(2) (2005) 323 -- 339.

\bibitem{hochbruck2010exponential}
M.~Hochbruck, A.~Ostermann, Exponential integrators, Acta Numerica 19 (2010)
  209--286.

\bibitem{cox_2002}
S.~Cox, P.~Matthews, Exponential time differencing for stiff systems, Journal
  of Computational Physics 176~(2) (2002) 430 -- 455.

\bibitem{borel}
E.~Borel, Leçons sur les séries divergentes, Gauthier-Villars, 1901.

\bibitem{lutz1999}
D.~Lutz, M.~Miyake, R.~Schäfke, On the {B}orel summability of divergent
  solutions of the heat equation, Nagoya Mathematical Journal 154 (1999) 1--29.

\bibitem{lysik2009}
G.~Lysik, {Borel summable solutions of the Burgers equation}, Annales Polonici
  Mathematici 95 (2009) 187--197.

\bibitem{Costin_2006}
O.~Costin, S.~Tanveer, Borel summability of {N}avier-{S}tokes equation in
  $\mathbb{R}^3$ and small time existence, ArXiv Mathematics e-prints (dec
  2006).
\newblock \href {http://arxiv.org/abs/math/0612063}
  {\path{arXiv:math/0612063}}.

\bibitem{dyson_1952}
F.~Dyson, Divergence of perturbation theory in quantum electrodynamics,
  Physical Review 85 (1952) 631--632.

\bibitem{suslov05}
I.~Suslov, Divergent perturbation series, Journal of Experimental and
  Theoretical Physics 100~(6) (2005) 1188--1233.

\bibitem{kontopolous_2002}
G.~Kontopoulos, Order and chaos in dynamical astronomy, Astronomy and
  astrophysics library, Springer, Berlin, Heidelberg, New York, 2002.

\bibitem{thomann00}
J.~Thomann, Formal and numerical summation of formal power series solutions of
  {ODE}'s, Tech. rep., CIRM Luminy (2000).

\bibitem{jcp13}
D.~Razafindralandy, A.~Hamdouni, Time integration algorithm based on divergent
  series resummation, for ordinary and partial differential equations, Journal
  of Computational Physics 236 (2013) 56--73.

\bibitem{esaim14}
A.~Deeb, A.~Hamdouni, E.~Liberge, D.~Razafindralandy, Borel-{L}aplace summation
  method used as time integration scheme, ESAIM: Procedings and Surveys 45
  (2014) 318--327.

\bibitem{dcds16}
A.~Deeb, A.~Hamdouni, D.~Razafindralandy, Comparison between {B}orel-{P}ad\'e
  summation and factorial series, as time integration methods, Discrete and
  Continuous Dynamical Systems - Serie S 9~(2) (2016) 393--408.

\bibitem{ramis_poincare_1}
J.-P. Ramis, Poincaré et les développements asymptotiques ({P}remière
  partie), Gazettes des Mathématiques 133 (Juillet 2012).

\bibitem{ramis_poincare_2}
J.-P. Ramis, Les développements asymptotiques après {P}oincaré : continuité
  et... divergences, Gazettes des Mathématiques 134 (octobre 2012).

\bibitem{costin_2008_book}
O.~Costin, Asymptotics and {B}orel Summability, Monographs and Surveys in Pure
  and Applied Mathematics, CRC Press, 2008.

\bibitem{brezinski79}
C.~Brezinski, Rationnal approximation to formal power serie, Journal of
  Approximation Theory 25~(4) (1979) 295--317.

\bibitem{brezinski94}
C.~Brezinski, J.~Van~Iseghem, Pad\'e approximations, in: P.~G. Ciarlet, J.~L.
  Lions (Eds.), Handbook of Numerical Analysis, Vol.~3, Elsevier, 1994, pp. 47
  -- 222.

\bibitem{stroud66}
A.~Stroud, D.~Secrest, Gaussian quadrature formulas (without numerical tables),
  Prentice-Hall, 1966.

\bibitem{gonnet11}
P.~Gonnet, S.~G\"uttel, L.~Trefethen, Robust {P}ad\'e approximation via {SVD},
  SIAM Review 51~(1) (2013) 101--117.

\bibitem{cochelin1994}
B.~Cochelin, A path-following technique via an asymptotic-numerical method,
  Computers and Structures 53 (1994) 1181--1192.

\bibitem{zahrouni04}
H.~Zahrouni, W.~Aggoune, J.~Brunelot, M.~Potier-Ferry, Asymptotic numerical
  method for strong nonlinearities, Revue Europ\'eenne des El\'ements Finis
  13~(1-2) (2004) 97--118.

\bibitem{bucker06}
M.~B\"ucker, G.~Corliss, U.~Naumann, P.~Hovland, B.~Norris (Eds.), Automatic
  differentiation: applications, theory, and implementations, Vol.~50 of
  Lecture Notes in Computational Science and Engineering, Springer, 2006.

\bibitem{griewank08}
A.~Griewank, A.~Walther, Evaluating derivatives. {P}rinciples and techniques of
  algorithmic differentiation, 2nd Edition, Frontiers in Applied Mathematics,
  SIAM, 2008.

\bibitem{cochelin07}
B.~Cochelin, N.~Damil, M.~Potier-Ferry, M{\'e}thode asymptotique num{\'e}rique,
  Methodes num{\'e}riques, Hermes Lavoisier, 2007.

\bibitem{baker61}
G.~A. Baker~Jr., J.~Gammel, J.~Wills, An investigation of the applicability of
  the {P}ad\'e approximant method, Journal of Mathematical Analysis and
  Applications 2 (1961) 405--418.

\bibitem{nuttall_1970}
J.~Nuttall, The convergence of {P}adé approximants of meromorphic functions,
  Journal of Mathematical Analysis and Applications 31~(1) (1970) 147 -- 153.

\bibitem{baker00}
G.~Baker, Defects and the convergence of {P}adé approximants, Acta Applicandae
  Mathematica (2000).

\bibitem{baker75}
G.~Baker, Essentials of {P}ad{\'e} Approximants, Elsevier Science, 1975.

\bibitem{fehlberg70}
E.~Fehlberg, Klassische {R}unge-{K}utta-{F}ormeln vierter und niedrigerer
  {O}rdnung mit {S}chrittweiten-{K}ontrolle und ihre {A}nwendung auf
  {W}{\"a}rmeleitungsprobleme, Computing 6~(1) (1970) 61--71.

\bibitem{book:hairer}
E.~Hairer, S.~Nørsett, G.~Wanner, Solving Ordinary Differential Equations {I}:
  Nonstiff Problems, {S}econd {R}evised, {C}orrected third printing Edition,
  Springer Series in Computational Mathematics, Springer, 2008.

\bibitem{swart97}
J.~de~Swart, G.~S\"oderlind, On the construction of error estimators for
  implicit {R}unge-{K}utta methods, Journal of Computational and Applied
  Mathematics 86 (1997) 347--358.

\bibitem{tomas13}
T.~Co, Methods of {A}pplied {M}athematics for {E}ngineers {S}cientists,
  Michigan Technology University, Cambridge University Press, 2013.

\bibitem{hofbauer88}
J.~Hofbauer, K.~Sigmund, The Theory of Evolution and Dynamical Systems:
  Mathematical Aspects of Selection, London Mathematical Society Student Texts,
  Cambridge University Press, 1988.

\bibitem{korteweg95}
D.~Korteweg, G.~de~Vries, On the change of form of long waves advancing in a
  rectangular canal, and on a new type of long stationary waves, Philosophical
  Magazine 39~(240) (1895) 422--443.

\end{thebibliography}

\end{document}